\documentclass{dis}
\usepackage[dvips]{graphics}
\usepackage{epsf}
\usepackage{graphicx}
\usepackage{epsfig}
\usepackage{amsmath}
\usepackage[psamsfonts]{amssymb}
\usepackage[psamsfonts]{eucal}
\usepackage{amsthm}
\usepackage{textcomp}
\usepackage{mathptmx}
\theoremstyle{plain}
\newtheorem{theorem}[subsection]{Theorem}
\newtheorem{lemma}[subsection]{Lemma}

\newtheorem{prop}[subsection]{Proposition}
\theoremstyle{definition}
\newtheorem{defi}[subsection]{Definition}
\newtheorem{defis}[subsection]{Definitions}
\newtheorem{rmk}[subsection]{Remark}
\newtheorem{rmks}[subsection]{Remarks}
\newtheorem{ex}[subsection]{Example}
\newtheorem{exs}[subsection]{Examples}
\swapnumbers
\newcommand{\field}[1]{\mathbb{#1}}
\newcommand{\CC}{\field{C}}
\newcommand{\KK}{\field{K}}

\newcommand{\RR}{\field{R}}
\newcommand{\ZZ}{\field{Z}}
\def\id{\mathop{\rm id}\nolimits}

\def\ad{\mathop{\rm ad}\nolimits}

\def\toto{\mathop{\rightrightarrows}\limits}

\begin{document}
\mathclass{Primary 58H05; Secondary 53D17.}
\thanks{The author acknowledges the facilities offered by the team \lq\lq
Analyse alg\'ebrique\rq\rq\ of the French \lq\lq Institut de math\'ematiques de
Jussieu\lq\lq.}
\abbrevauthors{C.-M. Marle}
\abbrevtitle{Calculus on Lie algebroids}

\title{Calculus on Lie algebroids, Lie groupoids and Poisson manifolds
% To appear in Dissertationes Mathematicae
}

\author{Charles-Michel Marle}
\address{Universit{\'e} Pierre et Marie Curie\\4, place Jussieu\\
 75252 Paris cedex 05,  France\\ E-mail: cmm1934@orange.fr}

\maketitledis

\tableofcontents
\begin{abstract}
We begin with a short presentation of the basic concepts related to Lie
groupoids and Lie algebroids, but the main part of this paper deals with Lie
algebroids. A Lie algebroid over a manifold is a vector
bundle over that manifold whose properties are very similar to those of a
tangent bundle. Its dual bundle has properties very similar to those of a
cotangent bundle: in the graded algebra of sections of its
external powers, one can define an operator $d_E$ similar to the
exterior derivative. We present the theory of Lie derivatives,
Schouten-Nijenhuis brackets and exterior derivatives
in the general setting of a Lie algebroid, its dual bundle and
their exterior powers. All the results (which, for
their most part, are already known) are given with detailed proofs. In the
final sections, the results are applied to Poisson manifolds, whose links with
Lie algebroids are very close.
\end{abstract}
\makeabstract
To appear in {\it Dissertationes Mathematicae}
\vskip 1.0cm
\centerline{\large In memory of Paulette Libermann}

\chapter{Introduction}
\label{Introduction}
Groupoids are mathematical structures able to describe symmetry
properties more general than those described by groups. They were
introduced (and named) by H.~Brandt \cite{Brandt} in 1926. The reader will find
a general presentation of that important concept in \cite{Wein4} and
\cite{Cartier}. 

A groupoid with a topological structure (resp., a differtiable structure) is
called a \emph{topological groupoid} (resp., a \emph{Lie groupoid}). Around
1950, Charles Ehresmann \cite{Ehres} used these concepts as essential tools in
Topology and differential Geometry.

In recent years, Mickael Karasev \cite{Kara},  Alan
Weinstein \cite{Wein2,
Coste} and Stanis\l aw Zakrzewski \cite{Zakr} independently discovered that
Lie groupoids equipped with a symplectic structure can be used for
the construction of noncommutative deformations of the algebra of smooth
functions on a manifold, with potential applications to quantization. Poisson
groupoids were introduced by Alan Weinstein \cite{Wein3} as generalizations of
both Poisson Lie groups and symplectic groupoids.

The infinitesimal counterpart of the notion of a Lie groupoid is the notion of
a Lie algebroid, in the same way as the infinitesimal counterpart of the notion
of a Lie group is the notion of a Lie algebra. Lie algebroids were first
considered by Jean Pradines \cite{Pra}.

Lie groupoids and Lie algebroids  are now an active domain of research, with
applications in various parts of mathematics \cite{Wein3, CannW, BanKo, deLeon,
Lib1, Cartier}. More specifically, Lie algebroids have applications in
Mechanics \cite{Wein5, Lib2, GraGraUrb, Martinez} and provide a very natural setting in
which one can develop the theory of differential operators such as the
exterior derivative of forms and the Lie derivative with respect
to a vector field. In such a setting, slightly more general than
that of the tangent and cotangent bundles to a smooth manifold and
their exterior powers, the theory of Lie derivatives extends, in a
very natural way, into the theory of the Schouten-Nijenhuis
bracket (first introduced in differential Geometry by
J.A.~Schouten \cite{Schou} and developed by A.~Nijenhuis
\cite{Ni}). Other bidifferential operators such as the bracket of
exterior forms on a Poisson manifold, first discovered for Pfaff
forms by F.~Magri and C.~Morosi \cite{MaMo} and extended to forms
of all degrees by J.-L.~Koszul \cite{Ko2} appear in such a setting
as very natural: they are Schouten-Nijenhuis brackets for the Lie
algebroid structure of the cotangent bundle to a Poisson manifold.
\par\smallskip
We first present in this paper the basic concepts related to Lie groupoids and
Lie algebroids. Then we develop the theory of Lie derivatives,
Schouten-Nijenhuis brackets and exterior derivatives in the
general setting of a Lie algebroid, its dual bundle and their
exterior powers. All the results (which, for their most part, are
already known, see for example \cite{Xu, GraUrb1, GraUrb2}) are given with
detailed proofs. Most of these proofs are the same as the classical ones (when
the Lie algebroid is the tangent bundle to a smooth manifold); a few ones are
slightly more complicated because, contrary to the algebra of exterior
differential forms on a manifold, the algebra of sections of
exterior powers of the dual of a Lie algebroid is not locally
generated by its elements of degree $0$ and their differentials.
These results may even be extended for more general algebroids with no
assumption of skew-symmetry \cite{GraUrb3}, but here we will not discuss these
generalizations, nor will we discuss the Schouten bracket for symmetric tensors.
In the final section, the results are applied to
Poisson manifolds. We show that the cotangent space of a Poisson manifold
has a Lie algebroid structure and that the total space of the vector bundle dual
to a Lie algebroid has a natural Poisson structure, and we use these
properties for lifting to the tangent bundle Poisson structures and Lie
algebroid structures. 

\chapter{Lie groupoids}
\label{Lie groupoids}
\section{Definition and first properties}
Before stating the formal definition  of a groupoid, let us
explain, in an informal way, why it is a very natural concept. The
easiest way to understand that concept is to think of two sets,
$\Gamma$ and $\Gamma_0$. The first one, $\Gamma$, is called the
\emph{set of arrows} or \emph{total space} of the groupoid, and
the other one, $\Gamma_0$, the \emph{set of objects} or \emph{set
of units} of the groupoid. One may think of an element $x\in
\Gamma$ as an arrow going from an object (a point in $\Gamma_0$)
to another object (another point in $\Gamma_0$). The word \lq\lq
arrow\rq\rq\ is used here in a very general sense: it means a way
for going from a point in $\Gamma_0$ to another point in
$\Gamma_0$. One should not think of an arrow as a line drawn in
the set $\Gamma_0$ joining the starting point of the arrow to its
end point. Rather, one should think of an arrow as living outside
$\Gamma_0$, with only its starting point and its end point in $\Gamma_0$,
as shown on Figure~1.
\par\smallskip

\begin{figure}[h]
\centerline{\includegraphics{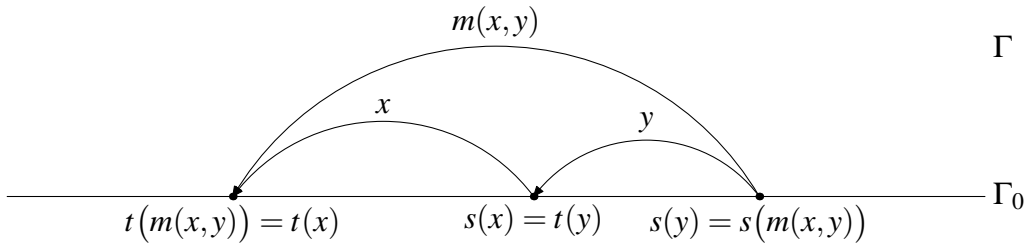}}
\caption{Two arrows $x$ and $y\in\Gamma$, with the target of $y$,
$t(y)\in \Gamma_0$, equal to the source of $x$,
$s(x)\in\Gamma_0$, and the composed arrow $m(x,y)$.}
\end{figure}

The following ingredients enter the definition of a groupoid.

\begin{description}

\item{(i)} Two maps $s:\Gamma\to\Gamma_0$ and $t:\Gamma\to
\Gamma_0$, called the {\it source map\/} and the {\it target
map\/} of the groupoid. If $x\in\Gamma$ is an arrow,
$s(x)\in\Gamma_0$ is its starting point and $t(x)\in\Gamma_0$
its end point.

\item{(ii)} A {\it composition law\/} on the set of arrows; we can
compose an arrow $y$ with another arrow $x$, and get an arrow
$m(x,y)$, by following first the arrow $y$, then the arrow $x$. Of
course, $m(x,y)$ is defined if and only if the target of $y$ is
equal to the source of $x$. The source of $m(x,y)$ is equal to the
source of $y$, and its target is equal to the target of $x$, as
illustrated on Figure~1. It is only by convention that we write
$m(x,y)$ rather than $m(y,x)$: the arrow which is followed first
is on the right, by analogy with the usual notation $f\circ g$ for
the composition of two maps $g$ and $f$. The composition of
arrows is associative.

\item{(iii)} An {\it embedding\/} $\varepsilon$ of the set $\Gamma_0$
into the set $\Gamma$, which associates a unit arrow
$\varepsilon(u)$ with each $u\in\Gamma_0$. That unit arrow is such
that both its source and its target are $u$, and it plays the role
of a unit when composed with another arrow, either on the right or
on the left: for any arrow $x$,
$m\Bigl(\varepsilon\bigl(t(x)\bigr),x\Bigr)=x$, and
$m\Bigl(x, \varepsilon\bigl(s(x)\bigr)\Bigr)=x$.

\item{(iv)} Finally, an {\it inverse map\/} $\iota$ from the set of
arrows onto itself. If $x\in\Gamma$ is an arrow, one may think of
$\iota(x)$ as the arrow $x$ followed in the reverse sense.
\end{description}

\par\smallskip
Now we are ready to state the formal definition of a groupoid.

\begin{defi}
A \emph{groupoid} is a pair of sets $(\Gamma,\Gamma_0)$ equipped
with the structure defined by the following data:
\begin{description}

\item{(i)} an injective map $\varepsilon:\Gamma_0\to\Gamma$, called
the \emph{unit section} of the groupoid;

\item{(ii)} two maps $s:\Gamma\to\Gamma_0$ and
$t:\Gamma\to\Gamma_0$, called, respectively, the \emph{source
map} and the \emph{target map\/}; they satisfy
\begin{equation*}
s\circ\varepsilon=t\circ\varepsilon=\id_{\Gamma_0}\,;
\end{equation*}

\item{(iii)} a composition law $m:\Gamma_2\to \Gamma$, called the
\emph{product}, defined on the subset $\Gamma_2$ of
$\Gamma\times\Gamma$, called the \emph{set of composable
elements},
\begin{equation*}
\Gamma_2=\bigl\{\,(x,y)\in\Gamma\times\Gamma;
s(x)=t(y)\,\bigr\}\,,
\end{equation*}
which is associative, in the sense that whenever one side of the
equality
 \begin{equation*}
 m\bigl(x,m(y,z)\bigr)=m\bigl(m(x,y),z\bigr)
 \end{equation*}
is defined, the other side is defined too, and the equality holds;
moreover, the composition law $m$ is such that for each
$x\in\Gamma$,
\begin{equation*}
 m\Bigl(\varepsilon\bigl(t(x)\bigr),x\Bigr)
 =m\Bigl(x,\varepsilon\bigl(s(x)\bigr)\Bigr)=x\,;
\end{equation*}

\item{(iv)} a map $\iota:\Gamma\to\Gamma$, called the
\emph{inverse}, such that, for every $x\in\Gamma$,
$\bigl(x,\iota(x)\bigr)\in\Gamma_2$, $\bigr(\iota(x),x\bigr)\in
\Gamma_2$ and
\begin{equation*}
m\bigl(x,\iota(x)\bigr)=\varepsilon\bigl(t(x)\bigr)\,,\quad
m\bigl(\iota(x),x\bigr)=\varepsilon\bigl(s(x)\bigr)\,.
\end{equation*}
\end{description}
The sets $\Gamma$ and $\Gamma_0$ are called, respectively, the
\emph{total space} and the \emph{set of units} of the groupoid,
which is itself denoted by $\Gamma\toto^t_s\Gamma_0$.
\end{defi}

\begin{rmk} The definition of a groupoid can be stated very briefly in the
language of category theory: a \emph{groupoid} is a small category all of
whose arrows are invertible. We recall that a category is said to be
\emph{small} if the collections of
its arrows and of its objects are sets.
\end{rmk} 

\subsection{Identification and notations} In what follows, by
means of the injective map $\varepsilon$, we will identify the set
of units $\Gamma_0$ with the subset $\varepsilon(\Gamma_0)$ of
$\Gamma$. Therefore $\varepsilon$ will be the canonical injection
in $\Gamma$ of its subset $\Gamma_0$.
\par
For $x$ and $y\in\Gamma$, we will sometimes write $x\circ y$, $x.y$, or even
simply $xy$ for $m(x,y)$, and $x^{-1}$ for $\iota(x)$. Also we
will write \lq\lq the groupoid $\Gamma$\rq\rq\ for \lq\lq the
groupoid $\Gamma\toto^t_s \Gamma_0$\rq\rq.

\section{Properties and comments}
The above definition has the following consequences.

\subsection{Involutivity of the inverse map} The inverse
map $\iota$ is involutive:
 \begin{equation*}
 \iota\circ\iota=\id_\Gamma\,.
 \end{equation*}
 We have indeed, for any
$x\in\Gamma$,
 \[
 \begin{split}
 \iota\circ\iota(x)&=m\bigl(\iota\circ\iota(x),
s\bigl(\iota\circ\iota(x)\bigr)\bigr)
        =m\bigl(\iota\circ\iota(x),s(x)\bigr)
        =m\bigl(\iota\circ\iota(x),m\bigl(\iota(x),x\bigr)\bigr)\\
       &=m\bigl(m\bigl(\iota\circ\iota(x),\iota(x)\bigr),x\bigr)
        =m\bigl(t(x),x\bigr)
         =x\,.\\
 \end{split}
 \]

\subsection{Unicity of the inverse} Let $x$ and $y\in
\Gamma$ be such that
\[m(x,y)=t(x)\quad\hbox{and}\quad m(y,x)=s(x)\,.\]
Then we have
 \[\begin{split}
   y&=m\bigl(y,s(y)\bigr)=m\bigl(y,t(x)\bigr)
     =m\bigl(y,m\bigl(x,\iota(x)\bigr)\bigr)
     =m\bigl(m(y,x),\iota(x)\bigr)\\
    &=m\bigl(s(x),\iota(x)\bigr)
     =m\bigl(t\bigl(\iota(x)\bigr),\iota(x)\bigr)
     =\iota(x)\,.\\
 \end{split}\]
Therefore for any $x\in \Gamma$, the unique $y\in\Gamma$ such that
$m(y,x)=s(x)$ and $m(x,y)=t(x)$ is $\iota(x)$.

\subsection{The fibres of the source and target maps and the
isotropy
groups}\label{isotropy} The target map $t$ (resp. the source
map $s$) of a groupoid $\Gamma\toto^t_s\Gamma_0$
determines an equivalence relation on $\Gamma$: two elements $x$
and $y\in\Gamma$ are said to be $t$-equivalent (resp.
$s$-equivalent) if $t(x)=t(y)$ (resp. if
$s(x)=s(y)$). The corresponding equivalence classes are
called the \emph{$t$-fibres} (resp. the
\emph{$s$-fibres}) of the groupoid. They are of the form
$t^{-1}(u)$ (resp. $s^{-1}(u)$), with $u\in\Gamma_0$.
\par
For each unit $u\in\Gamma_0$, the subset
 \begin{equation*}
 \Gamma_u=t^{-1}(u)\cap s^{-1}(u)
 =\bigl\{\,x\in\Gamma;s(x)=t(x)=u\,\bigr\}
 \end{equation*}
is called the {\it isotropy group\/} of $u$. It is indeed a group,
with the restrictions of $m$ and $\iota$ as composition law and
inverse map.

\begin{figure}[h]
\centerline{\includegraphics{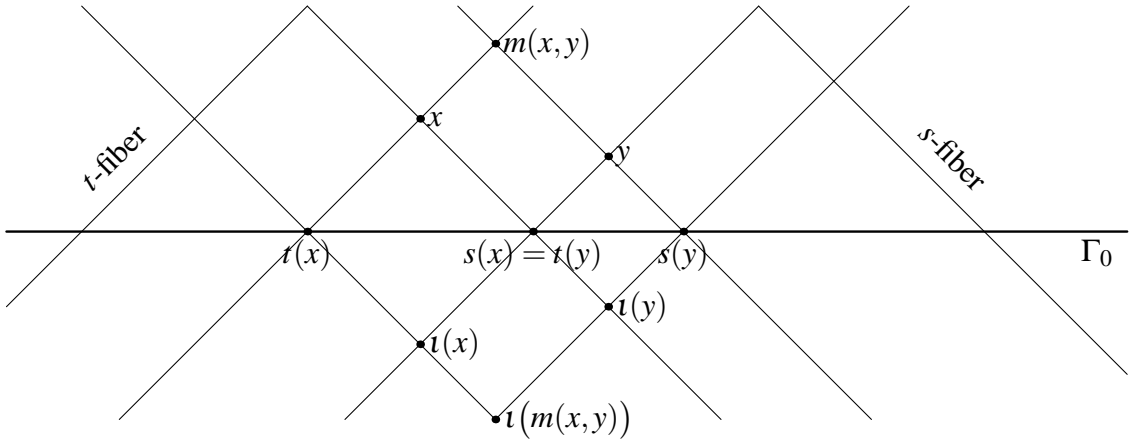}}
\caption{A way to visualize groupoids.}
\end{figure}

\subsection{A way to visualize groupoids} We have seen
(Figure~1) a way in which groupoids may be visualized, by using
arrows for elements in $\Gamma$ and points for elements in
$\Gamma_0$. There is another, very useful way to visualize
groupoids, shown on Figure~2. The total space $\Gamma$ of the groupoid is
represented as a plane, and the set $\Gamma_0$ of units as a straight line in
that plane. The $t$-fibres (resp. the $s$-fibres) are
represented as parallel straight lines, transverse to $\Gamma_0$.

Such a visualization should be used with care:
one may think, at first sight, that there is only one element in the groupoid
with a given source and a given target, which is not true in general.

\section{Simple examples of groupoids}
\subsection{The groupoid of pairs}
\label{groupoid of pairs}
Let $E$ be a nonempty set. Let $\Gamma=E\times E$, $\Gamma_0=E$,
$s:E\times E\to E$ be the projection on the right factor $s(x,y)=y$,
$t:E\times E\to E$ the projection on the left factor $t(x,y)=x$,
$\varepsilon:E\to E\times E$ be the diagonal map $x\mapsto (x,x)$. We define
the composition law $m:(E\times E)\times(E\times E)\to E\times E$ and the
inverse $\iota:E\times E\to E\times E$ by
 $$m\bigl((x,y),(y,z)\bigr)=(x,z)\, ,\quad \iota(x,y)=(y,x)\,.$$
Then $E\times E\toto^t_sE$ is a groupoid, called the \emph{groupoid of pairs} of
elements in $E$.
\subsection{Equivalence relations}
Let $E$ be a nonempty set with an equivalence relation $r$. Let
$\Gamma=\{\,(x,y)\in E\times E; x\,r\,y\,\}$ and $\Gamma_0=E$.
The source and target maps $s$ and $t$ are the restrictions to $\Gamma$ of the
source and target maps, above defined on $E\times E$ for the groupoid of pairs.
The composition law $m$, the injective map $\varepsilon$ and the inverse
$\iota$, are the same as for th egroupoid of pairs, suitably restricted. Then
$\Gamma\toto^t_sE$ is a groupoid, more precisely a subgroupoid of the groupoid
of pairs of elements in E.

\begin{rmk}
This example shows that equivalence relations may be considered as special
groupoids. Conversely, on the set of units $\Gamma_0$ of a general groupoid
$\Gamma\toto^t_s\Gamma_0$, there is a natural equivalence relation: $u_1$ and
$u_2\in \Gamma_0$ are said to be equivalent if there exists $x\in\Gamma$ such
that $s(x)=u_1$ and $t(x)=u_2$. But the groupoid structure generally carries
more information than that equivalence relation: there may be several $x\in
\Gamma$ such that $s(x)=u_1$ and $t(x)=u_2$, \emph{i.e.}, several ways in which
$u_1$ and $u_2$ are equivalent.
\end{rmk}

\section{Topological and Lie groupoids}
\begin{defis}
A \emph{topological groupoid} is a groupoid
$\Gamma\toto^t_s\Gamma_0$ for which $\Gamma$ is a (maybe
non Hausdorff) topological space, $\Gamma_0$ a Hausdorff
topological subspace of $\Gamma$, $t$ and $s$ surjective
continuous maps, $m:\Gamma_2\to\Gamma$ a continuous map and
$\iota:\Gamma\to\Gamma$ an homeomorphism.
\par\smallskip
A \emph{Lie groupoid} is a groupoid
$\Gamma\toto^t_s\Gamma_0$ for which $\Gamma$ is a smooth
(maybe non Hausdorff) manifold, $\Gamma_0$ a smooth Hausdorff
submanifold of $\Gamma$, $t$ and $s$ smooth surjective
submersions (which implies that $\Gamma_2$ is a smooth submanifold
of $\Gamma\times\Gamma$), $m:\Gamma_2\to\Gamma$ a smooth map and
$\iota:\Gamma\to\Gamma$ a smooth diffeomorphism.
\end{defis}
\goodbreak
\section{Examples of topological and Lie groupoids}
\nobreak
\subsection{Topological groups and Lie groups} A topological
group (resp. a Lie group) is a topological groupoid (resp. a Lie
groupoid) whose set of units has only one element $e$.

\subsection{Vector bundles}
\label{vector bundle groupoid} 
A smooth vector bundle $\tau:E\to M$
on a smooth manifold $M$ is a Lie groupoid, with the base $M$ as
set of units (identified with the image of the zero section); the
source and target maps both coincide with the projection $\tau$,
the product and the inverse maps are the addition $(x,y)\mapsto
x+y$ and the opposite map $x\mapsto -x$ in the fibres.

\subsection{The fundamental groupoid of a topological space}
\label{fundamental groupoid} 
Let $M$ be a topological space. A \emph{path} in $M$ is a continuous
map $\gamma:[0,1]\to M$. We denote by $[\gamma]$ the homotopy
class of a path $\gamma$ and by $\Pi(M)$ the set of homotopy
classes of paths in $M$ (with fixed endpoints). For $[\gamma]\in
\Pi(M)$, we set
 $t\bigl([\gamma]\bigr)=\gamma(1)$,
   $s\bigl([\gamma]\bigr)=\gamma(0)$,
where $\gamma$ is any representative of the class $[\gamma]$. The
concatenation of paths determines a well defined composition law
on $\Pi(M)$, for which $\Pi(M)\toto^t_s M$ is a
topological groupoid, called the \emph{fundamental groupoid} of
$M$. The inverse map is $[\gamma]\mapsto [\gamma^{-1}]$, where
$\gamma$ is any representative of $[\gamma]$ and $\gamma^{-1}$ is
the path $t\mapsto\gamma(1-t)$. The set of units is $M$, if we
identify a point in $M$ with the homotopy class of the constant
path equal to that point. 
\par\smallskip
Given a point $x\in M$, the isotropy group of the fundamental groupoid of $M$
at $x$ is the fundamental group at that point.
\par\smallskip
When $M$ is a smooth manifold, the same construction can be made
with piecewise smooth paths, and the fundamental groupoid
$\Pi(M)\toto^t_s M$ is a Lie groupoid.

\subsection{The gauge groupoid of a fibre bundle with structure
group}
The structure of a locally trivial topological bundle $(B,p,M)$ 
with standard fibre $F$ and structure group a
topological group $G$ of homeomorphisms of $F$, 
is usually determined {\it
via\/} an admssible fibred atlas $(U_i,\varphi_i)$, $i\in I$. The
$U_i$ are open subsets of $M$ such that $\bigcup_{i\in I}U_i=M$. For each $i\in
I$, $\varphi_i$ is a homeomorphism of $U_i\times F$ onto $p^{-1}(U_i)$ which,
for each $x\in U_i$, maps $\{x\}\times F$ onto $p^{-1}(x)$. For each pair
$(i,j)\in I^2$ such that $U_i\cap U_j\neq\emptyset$, each $x\in U_i\cap U_j$,
the homeomorphism $\varphi_j\circ\varphi_i^{-1}$ restricted to $\{x\}\times F$
is an element of $G$ ($F$ being identified with $\{x\}\times F$).
Elements of $G$ are called \emph{admissible homeomorphisms} of $F$. Another,
maybe more natural, way of describing that structure, is by looking at the set
$\Gamma$ of admissible homeomorphisms between two fibres of that fibre bundle,
$B_x=p^{-1}(x)$ and $B_y=p^{-1}(y)$, with $x$ and $y\in M$. The set $\Gamma$ has
a topological structure (in general not Hausdorff). For $\gamma\in \Gamma$
mapping $B_x$ onto $B_y$, we define $s(\gamma)=x$, $t(\gamma)=y$. Then 
$\Gamma\toto^t_s M$ is a topological groupoid,called the \emph{gauge groupoid}
of the fibre bundle $(B,p,M)$. When that bundle is smooth, its gauge groupoid
is a Lie groupoid. 

\section{Properties of Lie groupoids}

\subsection{Dimensions} Let $\Gamma\toto^t_s\Gamma_0$ be a Lie
groupoid. Since $t$ and $s$ are submersions, for any
$x\in\Gamma$, the $t$-fibre
$t^{-1}\bigl(t(x)\bigr)$ and the $s$-fibre
$s^{-1}\bigl(s(x)\bigr)$ are submanifolds of $\Gamma$,
both of dimension $\dim\Gamma-\dim\Gamma_0$. The inverse map
$\iota$, restricted to the $t$-fibre through $x$ (resp. the
$s$-fibre through $x$) is a diffeomorphism of that fibre onto
the $s$-fibre through $\iota(x)$ (resp. the $t$-fibre
through $\iota(x)$). The dimension of the submanifold $\Gamma_2$
of composable pairs in $\Gamma\times \Gamma$ is
$2\dim\Gamma-\dim\Gamma_0$.

\subsection{The tangent bundle of a Lie groupoid}
Let $\Gamma\toto^t_s\Gamma_0$ be a Lie
groupoid. Its tangent bundle $T\Gamma$ is a Lie groupoid, with
$T\Gamma_0$ as set of units, $Tt:T\Gamma\to T\Gamma_0$ and
$Ts:T\Gamma\to T\Gamma_0$ as target and source maps. Let us
denote by $\Gamma_2$ the set of composable pairs in
$\Gamma\times\Gamma$, by $m:\Gamma_2\to\Gamma$ the composition law
and by $\iota:\Gamma\to\Gamma$ the inverse. Then the set of
composable pairs in $T\Gamma\times T\Gamma$ is simply $T\Gamma_2$,
the composition law on $T\Gamma$ is $Tm:T\Gamma_2\to T\Gamma$ and
the inverse is $T\iota:T\Gamma\to T\Gamma$.

When the groupoid $\Gamma$ is a Lie group $G$, the Lie groupoid
$TG$ is a Lie group too.

\begin{rmk}
The cotangent bundle of a Lie groupoid is a Lie groupoid, and more precisely a
symplectic groupoid \cite{CannW, Coste, Wein2, Albert, Daz}. Remarkably, the
cotangent bundle
of a
non-Abelian Lie group is not a Lie group: it is a Lie groupoid. This fact may
be considered as a justification of the current interest in Lie groupoids: as
soon as one is interested in Lie groups, by looking at their cotangent bundles, 
one has to deal with Lie groupoids! 
\end{rmk}

\subsection{Isotropy groups} For each unit $u\in\Gamma_0$ of a Lie groupoid,
the
isotropy group $\Gamma_u$ (defined in \ref{isotropy}) is a Lie
group.

\chapter{Lie algebroids}
\label{Lie algebroids}
\setcounter{equation}{0} \ The concept of a Lie algebroid was
first introduced by J.~Pradines \cite{Pra}, in relation with Lie
groupoids.
\section{Definition and examples}
 A Lie algebroid over a manifold is a vector bundle
based on that manifold, whose properties are very similar to those
of the tangent bundle. Let us give its formal definition.

 \begin{defi}
\label{Lie algebroid}
Let $M$ be a smooth manifold and $(E,\tau,M)$ be a vector bundle
with base $M$. A {\it Lie algebroid structure\/} on that bundle is
the structure defined by the following data:
 \begin{description}
 \item{(1)} a composition law $(s_1,s_2)\mapsto\{s_1,s_2\}$ on the
space $\Gamma(\tau)$ of smooth sections of that bundle, for which
$\Gamma(\tau)$ becomes a Lie algebra,

 \item{(2)} a smooth vector bundle map $\rho:E\to TM$, where $TM$
is the tangent bundle of $M$, such that for every pair
$(s_1,s_2)$ of smooth sections of $\tau$, and every smooth function
$f:M\to\RR$, we have the Leibniz-type formula,
 \begin{equation*}\{s_1,fs_2\}=f\{s_1,s_2\}+\bigl({\cal L}(\rho\circ
s_1)f\bigr)s_2\,.
 \end{equation*}
We have denoted by ${\cal L}(\rho\circ
s_1)f$ the Lie derivative of $f$ with respect to the vector field 
$\rho\circ s_1$:
\begin{equation*}
 {\cal L}(\rho\circ
s_1)f=i(\rho\circ s_1)df\,.
\end{equation*}
\end{description}
The vector bundle $(E,\tau,M)$ equipped with its Lie algebroid
structure will be called a {\it Lie algebroid\/} and denoted by
$(E,\tau,M,\rho)$; the composition law
$(s_1,s_2)\mapsto\{s_1,s_2\}$ will be called the {\it bracket\/}
and the map $\rho:E\to TM$ the {\it anchor\/} of the Lie algebroid
$(E,\tau,M,\rho)$.
\end{defi}

\begin{prop}
\label{property of anchor}
Let $(E,\tau,M,\rho)$ be a Lie algebroid. The map $s\mapsto\rho\circ s$, which
associates to a smooth section $s$ of $\tau$ the smooth vector field
$\rho\circ s$ on $M$, is a Lie algebras homomorphism. In other words, for each
pair $(s_1,s_2)$ of smooth sections of $\tau$,
 \begin{equation*}[\rho\circ s_1,\rho\circ s_2]=\rho\circ\{s_1,s_2\}\,.
 \end{equation*}
\end{prop}

\begin{proof}
 Let $s_1$, $s_2$ and $s_3$ be three smooth sections of $\tau$ and $f$ be a
smooth function on $M$.
By the Jacobi identity for the Lie algebroid bracket,
\begin{equation*}
 \bigl\{\{s_1,s_2\},fs_3\bigr\}=\bigl\{s_1,\{s_2,fs_3\}\bigr\}-
 \bigl\{s_2,\{s_1,fs_3\}\bigr\}\,.
\end{equation*}
Using the property of the anchor, the right hand side becomes
\begin{equation*}
 \begin{split}
 \bigl\{s_1,\{s_2,fs_3\}\bigr\}-
 \bigl\{s_2,\{s_1,fs_3\}\bigr\}
 &=f\bigl(\bigl\{s_1,\{s_2,s_3\}\bigr\}
 -\bigl\{s_2,\{s_1,s_3\}\bigr\}\bigr)\\
 &\quad +\Bigl(\bigl({\cal L}(\rho\circ s_1)\circ{\cal L}(\rho\circ s_2)
  -{\cal L}(\rho\circ s_2)\circ{\cal L}(\rho\circ s_1)\bigr)f\Bigr)s_3\,.
 \end{split}
\end{equation*}
Similarly, the left hand side becomes
\begin{equation*}
 \bigl\{\{s_1,s_2\},fs_3\bigr\}=f\bigl\{\{s_1,s_2\},s_3\bigr\}
 +\Bigl({\cal L}\bigl(\rho\circ\{s_1,s_2\}\bigr)f\Bigr)s_3\,.
\end{equation*}
Using again the Jacobi identity for the Lie algebroid bracket, we obtain
 \begin{equation*}
 \Bigl(\Bigl({\cal L}\bigl(\rho\circ\{s_1,s_2\}\bigr)
 -\bigl({\cal L}(\rho\circ s_1)\circ{\cal L}(\rho\circ s_2)
  -{\cal L}(\rho\circ s_2)\circ{\cal L}(\rho\circ
   s_1)\bigr)\Bigr)f\Bigr)s_3=0\,. 
\end{equation*}
But we have
\begin{equation*}
 {\cal L}(\rho\circ s_1)\circ{\cal L}(\rho\circ s_2)
  -{\cal L}(\rho\circ s_2)\circ{\cal L}(\rho\circ
   s_1)={\cal L}\bigl([\rho\circ s_1,\rho\circ s_2]\bigr)\,.
\end{equation*}
Finally
\begin{equation*}
 \Bigl({\cal L}\bigl(\rho\circ\{s_1,s_2\}-[\rho\circ s_1,\rho\circ
s_2]\bigr)f\Bigr)s_3=0\,.
\end{equation*}
This result, which holds for any smooth function $f$ on $M$ and any smooth
section $s_3$ of $\tau$, proves that $s\mapsto\rho\circ s$ is a Lie algebras
homomorphism. 
\end{proof}

\begin{rmks} Let $(E,\tau,M,\rho)$ be a Lie algebroid.
\par\nobreak\smallskip\noindent
{\rm(i)}\quad{\it Lie algebras homomorphisms.}\quad
For each smooth vector field $X$ on $M$, the Lie derivative  ${\cal L}(X)$
with respect to $X$ is a derivation of $C^\infty(M,\RR)$:
for every pair $(f,g)$ of smooth functions on $M$,
 \begin{equation*}{\cal L}(X)(fg)=\bigl({\cal L}(X)f\bigr)g+f\bigl({\cal
 L}(X)g\bigr)\,.
 \end{equation*}
The map $X\mapsto{\cal L}(X)$ is a Lie algebras
homomorphism from the Lie algebra $A^1(M)$ of smooth vector fields
on $M$, into the Lie algebra
 $\mathop{\rm Der}\bigl(C^\infty(M,\RR)\bigr)$
of derivations of $C^\infty(M,\RR)$, equipped with the commutator
 \begin{equation*}(D_1,D_2)\mapsto[D_1,D_2]=D_1\circ
 D_2-D_2\circ D_1
 \end{equation*}
as composition law. These facts were used in the proof of
Proposition~\ref{property of anchor}.
\par\smallskip
The map $s\mapsto {\cal L}(\rho\circ s)$, obtained by composition
of two Lie algebras homomorphisms, is a Lie algebras homomorphism,
from the Lie algebra $\Gamma(\tau)$ of smooth sections of the Lie
algebroid $(E,\tau,M,\rho)$, into the Lie algebra of derivations of
$C^\infty(M,\RR)$.
\par\smallskip\noindent
{\rm(ii)}\quad{\it Leibniz-type formulae\/.}\quad According to
Definition \ref{Lie algebroid} we have, for any pair $(s_1,s_2)$ of smooth
sections of $\tau$ and any smooth function $f$ on $M$,
 \begin{equation*}\{s_1,fs_2\}=f\{s_1,s_2\}+\bigl(i(\rho\circ
 s_1)df\bigr)\,s_2\,.
 \end{equation*}
As an easy consequence of the definition, we also have
 \begin{equation*}\{fs_1,s_2\}=f\{s_1,s_2\}-\bigl(i(\rho\circ
 s_2)df\bigr)\,s_1\,.
 \end{equation*}
More generally, for any pair $(s_1,s_2)$ of smooth sections of
$\tau$ and any pair $(f_1,f_2)$ of smooth functions  on $M$, we
have
 \begin{equation*}\{f_1s_1,f_2s_2\}=f_1f_2\{s_1,s_2\}+f_1\bigl(i(\rho\circ
 s_1)df_2\bigr)s_2-f_2\bigl(i(\rho\circ
 s_2)df_1\bigr)s_1\,.
 \end{equation*}
Using the Lie derivative operators, that formula may also be
written as
 \begin{equation*}\{f_1s_1,f_2s_2\}=f_1f_2\{s_1,s_2\}
 +f_1\bigl({\cal L}(\rho\circ
 s_1)f_2\bigr)s_2-f_2\bigl({\cal L}(\rho\circ
 s_2)f_1\bigr)s_1\,.
 \end{equation*}
\end{rmks}

\subsection{Simple examples of Lie algebroids}\hfill
\par\smallskip\noindent
{\rm(i)}\quad{\it The tangent bundle\/.}\quad The tangent bundle
$(TM,\tau_M,M)$ of a smooth manifold $M$, equipped with the usual
bracket of vector fields as composition law and with the identity
map $\id_{TM}$ as anchor, is a Lie algebroid.
\par\smallskip\noindent
{\rm(ii)}\quad{\it An involutive distribution\/.}\quad Let $V$ be
a smooth distribution on a smooth manifold $M$, {\it i.e.}, a
smooth vector subbundle of the tangent bundle $TM$. We assume that
$V$ is involutive, {\it i.e.}, such that the space of its smooth
sections is stable under the bracket operation. The vector bundle
$(V,\tau_M|_V,M)$, with the usual bracket of vector fields as
composition law and with the canonical injection $i_V:V\to TM$ as
anchor, is a Lie algebroid. We have denoted by $\tau_M:TM\to M$
the canonical projection of the tangent bundle and by $\tau_M|_V$
its restriction to the subbundle $V$.
\par\smallskip\noindent
{\rm(iii)}\quad{\it A sheaf of Lie algebras\/.}\quad Let
$(E,\tau,M)$ be a vector bundle over the smooth manifold $M$ and
$(z_1,z_2)\mapsto[z_1,z_2]$ be a smooth, skew-symmetric bilinear
bundle map defined on the fibred product $E\times_M E$, with
values in $E$, such that for each $x\in M$, the fibre
$E_x=\tau^{-1}(x)$, equipped with the bracket $(z_1,z_2)\mapsto
[z_1,z_2]$, is a Lie algebra. We define the bracket of two smooth
sections $s_1$ and $s_2$ of $\tau$ as the section $\{s_1,s_2\}$
such that, for each $x\in M$,
$\{s_1,s_2\}(x)=\bigl[s_1(x),s_2(x)\bigr]$. For the anchor, we
take the zero vector bundle map from $E$ to $TM$. Then $(E,\tau,M)$
is a Lie algebroid of particular type, called a {\it sheaf of Lie
algebras\/} over the manifold $M$.

\par\smallskip\noindent
{\rm(iv)}\quad{\it A finite-dimensional Lie algebra\/.}\quad In
particular, a finite-dimensional Lie algebra can be considered as
a Lie algebroid over a base reduced to a single point, with the
zero map as anchor.

%\end{exs}

\section{The Lie algebroid of a Lie groupoid} 
We now describe the most
important example of Lie algebroid: to every Lie groupoid, there is an
associated Lie algebroid, much like to every Lie group there is an associated
Lie algebra. It is in that setting that Pradines \cite{Pra} introduced Lie
algebroids for the first time. For more informations about Lie groupoids and
their associated Lie algebroids, the reader is referred to
\cite{Mack, Mack2, Coste, DazSon, Albert}.

In the following propositions and definitions, $\Gamma\toto^t_s\Gamma_0$ is
a Lie groupoid.

\begin{prop}
For each $x\in \Gamma$, the
maps
$$y\mapsto L_x(y)=xy\quad\hbox{and}\quad z\mapsto R_x(z)=zx$$
are smooth diffeomorphisms, respectively from $t^{-1}\bigl(s(x)\bigr)$ onto
$t^{-1}\bigl(t(x)\bigr)$ and from
$s^{-1}\bigl(t(x)\bigr)$ onto
$s^{-1}\bigl(s(x)\bigr)$. These maps are called the \emph{left translation}
and the \emph{right translation} by $x$, respectively.
\end{prop}

\begin{proof}
 The smoothness of the groupoid composition law $m:(x,y)\mapsto xy$ implies the
smoothness of $L_x$ and $R_x$. These maps are diffeomorphisms whose inverses are
 $$(L_x)^{-1}=L_{x^{-1}}\,,\quad (R_x)^{-1}=R_{x^{-1}}\,,$$
so the proof is complete.
\end{proof}

\begin{defi}
A vector field $Y$ and a vector field $Z$, defined on open subsets of
$\Gamma$, are said to be, respectively, \emph{left invariant} and \emph{right
invariant} if they satisfy the two properties
\begin{description}
 \item{(i)} the projections on $\Gamma_0$  of $Y$ by the target map $t$, and
of $Z$ by the source map $s$, vanish:
$$Tt(Y)=0\,,\quad Ts(Z)=0\,;$$

\item{(ii)} for each $y$ in the domain of definition of $Y$ and each $x\in
s^{-1}\bigl(t(y)\bigr)$, $xy$ is in the domain of definition of $Y$ and
 $$Y(xy)= TL_x\bigl(Y(y)\bigr)\,;
$$ 
similarly, for each $z$ in the domain of definition of $Z$ and each $x\in
t^{-1}\bigl(s(z)\bigr)$, $zx$ is in the domain of definition of $Z$ and
 $$Z(zx)= TR_x\bigl(Z(z)\bigr)\,.
$$ 
\end{description}
\end{defi}

\begin{prop} Let
$A(\Gamma)$ be the intersection of $\ker Tt$ and
$T_{\Gamma_0}\Gamma$ (the tangent bundle $T\Gamma$ restricted to
the submanifold $\Gamma_0$). Then $A(\Gamma)$ is the total
space of a vector bundle $\tau:A(\Gamma)\to\Gamma_0$, with base
$\Gamma_0$, the canonical projection $\tau$ being the map which
associates a point $u\in \Gamma_0$ to every vector in $\ker
T_u t$. That vector bundle has a natural Lie algebroid structure and is called
the \emph{Lie algebroid} of the Lie groupoid $\Gamma$. Its composition law
is so defined. Let $w_1$ and $w_2$ be two smooth sections of that bundle over an
open subset $U$ of $\Gamma_0$. Let $\widehat w_1$ and $\widehat w_2$ be the two
left invariant vector fields, defined on $s^{-1}(U)$, whose restrictions to $U$
are equal to $w_1$ and $w_2$ respectively. Then for each $u\in U$,
 $$\{w_1,w_2\}(u)=[\widehat w_1, \widehat w_2](u)\,.$$
The anchor $\rho$ of that Lie algebroid is the map $Ts$ restricted to
$A(\Gamma)$.
\end{prop}

\begin{proof}
The correspondence which associates, to each smooth section $w$ of the vector
bundle $\tau:A(\Gamma)\to\Gamma_0$, the
prolongation of that section by a left invariant vector field $\widehat w$, is 
a vector space isomorphism. Therefore, by setting 
$$\{w_1,w_2\}(u)=[\widehat w_1, \widehat w_2](u)\,,$$
we obtain a Lie algebra structure on the space of smooth sections of
$\tau:A(\Gamma)\to\Gamma_0$. Let $f$ be a smooth function, defined on the open
subset $U$ of $\Gamma_0$ on which $w_1$ and $w_2$ are defined. For each $u\in
U$,
\begin{equation*}
 \begin{split}
\{w_1, fw_2\}(u)&=[\widehat w_1,\widehat{fw}_2](u)=[\widehat w_1, (f\circ
s)\widehat w_2](u)\\
 &=f(u)[\widehat w_1, \widehat w_2(u)]+ \bigl(i(\widehat w_1)
d(f\circ s)\bigr)(u)\widehat w_2(u)\\
 &= f(u)[\widehat w_1, \widehat w_2(u)]+\bigl\langle df(u),
Ts\bigl(\widehat w_1(u)\bigr)\bigr\rangle \widehat w_2(u)\\
&= f(u)\{w_1,w_2\}(u)+\bigl\langle df(u), Ts\bigl(w_1(u)\bigr)\bigr\rangle
w_2(u)\,,
 \end{split}
\end{equation*}
which proves that $Ts$ has the properties of an anchor.
\end{proof}

\begin{rmk}
We could exchange the roles of $t$ and $s$ and use
\emph{right invariant} vector fields instead of left invariant
vector fields. The Lie algebroid obtained remains the same, up to
an isomorphism.
\end{rmk}

\begin{exs}
\par\nobreak\noindent
{\rm(i)}\quad
When the Lie groupoid $\Gamma\toto^t_s$ is a Lie group,
its Lie algebroid is simply its Lie algebra.
\par\nobreak\smallskip\noindent
{\rm(ii)}\quad
We have seen (\ref{vector bundle groupoid}) that a vector bundle $(E,\tau,M)$,
 with addition in the fibres as composition law, can be considered as a Lie
groupoid. Its Lie algebroid is the same vector bundle, with the zero bracket on
its space of sections, and the zero map as anchor. 
\par\nobreak\smallskip\noindent
{\rm(ii)}\quad
Let $M$ be a smooth manifold. The groupoid of pairs $M\times M\toto^t_sM$
(\ref{groupoid of pairs}) is a Lie groupoid whose Lie algebroid is isomorphic
to the tangent bundle $(TM, \tau_M,M)$ with the identity map as anchor.
\par\nobreak\smallskip\noindent
{\rm(iii)}\quad
The fundamental groupoid (\ref{fundamental groupoid}) of a
smooth connected manifold $M$ is a Lie groupoid. Its total space is the simply
connected covering space of $M\times M$ and, as in the previous example, its
Lie algebroid is isomorphic to the tangent bundle $(TM, \tau_M,M)$.
\end{exs} 

\subsection{Integration of Lie algebroids}
According to Lie's third theorem, for any given finite-dimensional
Lie algebra, there exists a Lie group whose Lie algebra is
isomorphic to that Lie algebra. The same property is not true for
Lie algebroids and Lie groupoids. The problem of finding necessary
and sufficient conditions under which a given Lie algebroid is
isomorphic to the Lie algebroid of a Lie groupoid remained open
for more than 30 years. Partial results were obtained by
J.~Pradines \cite{Pra2}, K.~Mackenzie \cite{Mack}, P.~Dazord
\cite{Daz}, P.~Dazord ans G.~Hector \cite{DazHect}. An important breakthrough
was made by Cattaneo and Felder~\cite{CattaneoFelder} who, starting from
a Poisson manifold, built a groupoid (today called the \emph{Weinstein
groupoid}) which, when its total space is regular, has a dimension twice that
of the Poisson manifold, has a symplectic structure and has as Lie algebroid the
cotangent space to the Poisson manifold. That groupoid was obtained by
symplectic reduction of an infinite dimensional manifold. That method may in 
fact be used for any Lie algebroid, as shown by Cattaneo~\cite{Cattaneo}.
A complete solution the integration problem of Lie algebroids was obtained by
M.~Crainic and R.L.~Fernandes \cite{crainicfernandes1}. They have shown that
with each given Lie algebroid, one can associate a topological groupoid with
connected and simply connected $t$-fibres, now called the \emph{Weinstein
groupoid} of that Lie algebroid. That groupoid, when the Lie algebroid is the
cotangent bundle to a Poisson manifold, is the same as that previously
obtained by Cattaneo and Felder by another method. When that topological
groupoid is in fact a Lie groupoid, {\it
i.e.}, when it is smooth, its Lie algebroid is isomorphic to the given Lie
algebroid. Crainic and Fernandes have obtained computable necessary and
sufficient conditions under which the Weinstein groupoid of a Lie algebroid is
smooth. In \cite{crainicfernandes2} they have used these results for the
integration of Poisson manifolds, {\it i.e.}, for the construction of a
symplectic groupoid whose set of units is a given Poisson manifold.

\section{Locality of the bracket}

We will prove that the value, at any point $x\in M$, of the
bracket of two smooth sections $s_1$ and $s_2$ of the Lie
algebroid $(E,\tau,M,\rho)$, depends only on the jets of order $1$
of $s_1$ and $s_2$ at $x$. We will need the following lemma.

\begin{lemma}
\label{locality}
 Let $(E,\tau,M,\rho)$ be a Lie algebroid, $s_1:M\to E$ a smooth
 section of $\tau$, and $U$ an open subset of $M$ on which $s_1$
 vanishes. Then for any other smooth section $s_2$ of
 $\tau$, $\{s_1,s_2\}$ vanishes on $U$.

\end{lemma}

\begin{proof} Let $x$ be a point in $U$. There exists a smooth
function $f:M\to \RR$, whose support is contained in $U$ such that
$f(x)=1$. The section $fs_1$ vanishes identically, since $s_1$
vanishes on $U$ while $f$ vanishes outside of $U$. Therefore, for
any other smooth section $s_2$ of $\tau$,
 \begin{equation*}0=\{fs_1,s_2\}=-\{s_2,fs_1\}=-f\{s_2,s_1\}-\bigl(i(\rho\circ
 s_2)df\bigr)s_1\,.\end{equation*}
So at $x$ we have
 \begin{equation*}f(x) \{s_1,s_2\}(x)=\bigl(i(\rho\circ
 s_2)df\bigr)(x)s_1(x)=0\,.\end{equation*}
Since $f(x)=1$, we obtain $\{s_1,s_2\}(x)=0$.
\end{proof}

\begin{prop}
 Let $(E,\tau,M,\rho)$ be a Lie algebroid. The value $\{s_1,s_2\}(x)$
 of the bracket of two smooth sections $s_1$ and $s_2$ of $\tau$,
 at a point $x\in M$, depends only on the jets of order $1$ of
 $s_1$ and $s-_2$ at $x$. Moreover, if $s_1(x)=0$ and $s_2(x)=0$, then
 $\{s_1,s_2\}(x)=0$.

\end{prop}

\begin{proof} Let $U$ be an open neighbourhood of $x$ in $M$ on
which there exists a local basis $(\sigma_1,\ldots,\sigma_k)$ of
of smooth sections of  $\tau$. For any point $y\in U$,
$\bigl(\sigma_1(y),\ldots,\sigma_k(y)\bigr)$ is a basis of the
fibre $E_y=\tau^{-1}(y)$. Let $s_1$ and $s_2$ be two smooth
sections of $\tau$. On the open subset $U$, these two sections can
be expressed, in a unique way, as
 \begin{equation*}s_1=\sum_{i=1}^k f_i\sigma_i\,,\qquad
  s_2=\sum_{j=1}^k g_j\sigma_j\,,\end{equation*}
where the $f_i$ and $g_j$ are smooth functions on $U$.
\par
By Lemma \ref{locality}, the values of $\{s_1,s_2\}$ in $U$ depend only on
the values of $s_1$ and $s_2$ in $U$. Therefore we have in $U$
 \begin{equation*}\{s_1,s_2\}=\sum_{i,j}\Bigl(f_ig_j\{\sigma_i,\sigma_j\}
 +f_i\bigl({\cal L}(\rho\circ \sigma_i)g_j\bigr)\sigma_j
 -g_j\bigl({\cal L}(\rho\circ
\sigma_j)f_i\bigr)\sigma_i\Bigr)\,.\end{equation*} This expression
proves that the value of $\{s_1,s_2\}$ at $x$ depends only on the
$f_i(x)$, $df_i(x)$, $g_j(x)$ and $dg_j(x)$, that means on the
jets of order $1$ of $s_1$ and $s_2$ at $x$.
\par
If $s_1(x)=0$, we have, for all $i\in\{1,\ldots,k\}$,  $f_i(x)=0$,
and similarly if $s_2(x)=0$, we have, for all
$j\in\{1,\ldots,k\}$, $g_j(x)=0$. The above expression shows then
that $\{s_1,s_2\}(x)=0$.
\end{proof}
\

\chapter{Exterior powers of vector bundles}
\label{Exterior powers}
We recall in this section some definitions and general properties
related to vector bundles, their dual bundles and exterior powers.
In a first subsection we recall some properties of graded
algebras, graded Lie algebras and their derivations. The second
subsection applies these properties to the graded algebra of
sections of the exterior powers of a vector bundle. For more
details the reader may look at the book by Greub, Halperin and
Vanstone \cite{GrHV}. The reader already familiar with this
material may skip this section, or just look briefly at the sign
conventions we are using.

\section{Graded vector spaces and graded algebras}

\begin{defis}
\label{graded definitions}
\par\nobreak\smallskip\noindent
{\rm(i)}\quad An {\it algebra\/} is a vector space $A$ on the
field $\KK=\RR$ or $\CC$, endowed with a $\KK$-bilinear map called
the {\it composition law\/},
 \begin{equation*}A\times A\to A\,,\quad (x,y)\mapsto xy\,,\quad \hbox{where}\
 x\ \hbox{and}\ y\in A\,.\end{equation*}
\par\smallskip\noindent
{\rm(ii)}\quad An algebra $A$ is said to be {\it associative\/} if
its composition law is associative, {\it i.e.}, if for all $x$,
$y$ and $z\in A$,
 \begin{equation*}x(yz)=(xy)z\,.\end{equation*}
\par\smallskip\noindent
{\rm(iii)}\quad A vector space $E$ on the field $\KK=\RR$ or $\CC$
is said to be {\it $\ZZ$-graded\/} if one has chosen a family
$(E^p\,,\ p\in\ZZ)$ of vector subspaces of $E$, such that
 \begin{equation*}E=\bigoplus_{p\in\ZZ}E^p\,.\end{equation*}
For each $p\in\ZZ$, an element $x\in E$ is said to be {\it
homogeneous of degree $p$ \/} if $x\in E^p$.
\par\smallskip\noindent
{\rm(iv)}\quad Let $E=\bigoplus_{p\in\ZZ}E^p$ and
$F=\bigoplus_{p\in\ZZ}F^p$ be two $\ZZ$-graded vector spaces on the
same field $\KK$. A $\KK$-linear map $f:E\to F$ is said to be {\it
homogeneous of degree $d$\/} (with $d\in\ZZ$) if for each
$p\in\ZZ$,
  \begin{equation*}f(E^p)\subset F^{p+d}\,.\end{equation*}
\par\smallskip\noindent
{\rm(v)}\quad An algebra $A$ is said to be $\ZZ$-graded if
$A=\bigoplus_{p\in\ZZ}A^p$ is $\ZZ$-graded as a vector space and if
in addition, for all $p$ and $q\in \ZZ$, $x\in A^p$ and $y\in
A^q$,
  \begin{equation*}xy\in A^{p+q}\,.\end{equation*}
\par\smallskip\noindent
{\rm(vi)}\quad A $\ZZ$-graded algebra $A=\bigoplus_{p\in\ZZ}A^p$ is
said to be {\it $\ZZ_2$-commutative\/}
if for all $p$ and $q\in\ZZ$, $x\in
A^p$ and $y\in A^q$,
 \begin{equation*}xy=(-1)^{pq}yx\,.
\end{equation*}
It is said to be {\it
$\ZZ_2$-anticommutative\/})
if for all $p$ and $q\in\ZZ$, $x\in
A^p$ and $y\in A^q$,
 \begin{equation*}
 xy=-(-1)^{pq}yx\,\hbox{)}\,.
\end{equation*}
\end{defis}

\subsection{Some~properties~and~examples}\hfill
\label{graded spaces properties}
\par\smallskip\nobreak\noindent
{\rm(i)}\quad{\it Composition of homogeneous linear maps\/.}\quad
We consider three $\ZZ$-graded vector spaces,
$E=\bigoplus_{p\in\ZZ}E^p$, $F=\bigoplus_{p\in\ZZ}F^p$ and
$G=\bigoplus_{p\in\ZZ}G^p$, on the same field $\KK$. Let $f:E\to F$
and $g:F\to G$ be two linear maps, the first one $f$ being
homogeneous of degree $d_1$ and the second one $g$ homogeneous of
degree $d_2$. Then $g\circ f:E\to G$ is homogeneous of degree
$d_1+d_2$.
\par\smallskip\nobreak\noindent
{\rm(ii)}\quad {\it The algebra of linear endomorphisms of a
vector space\/.}\quad Let $E$ be a vector space and ${\cal
L}(E,E)$ be the space of linear endomorphisms of $E$. We take as
composition law on that space the usual composition of maps,
 \begin{equation*}(f,g)\mapsto f\circ g\,,\quad\hbox{with}\quad f\circ
 g(x)=f\bigl(g(x)\bigr)\,,\ x\in E\,.\end{equation*}
With that composition law, ${\cal L}(E,E)$ is an associative
algebra.
\par\smallskip\nobreak\noindent
{\rm(iii)}\quad{\it The graded algebra of graded linear
endomorphisms\/.}\quad We assume now that $E=\bigoplus_{p\in\ZZ}E^p$
is a $\ZZ$-graded vector space. For each $d\in\ZZ$, let $A^d$ be
the vector subspace of ${\cal L}(E,E)$ whose elements are the
linear endomorphisms $f:E\to E$ which are homogeneous of degree
$d$, {\it i.e.}, such that for all $p\in\ZZ$, $f(E^p)\subset
E^{p+d}$. Let $A=\bigoplus_{d\in\ZZ}A^d$. By using
property~\ref{graded spaces properties}~(i), we see that with the usual
composition of maps
as composition law, $A$ is a $\ZZ$-graded associative algebra.
\par\bigskip
Let us use property~\ref{graded spaces properties}~(i) with $E=F=G$, in the
following
definition.

\begin{defi}
\label{graded bracket} 
Let $E=\bigoplus_{p\in\ZZ}E^p$ be a $\ZZ$-graded
vector space, $f$ and $g\in {\cal L}(E,E)$ be two homogeneous
linear endomorphisms of $E$ of degrees $d_1$ and $d_2$,
respectively. The linear endomorphism $[f,g]$ of $E$ defined by
 \begin{equation*}[f,g]=f\circ g-(-1)^{d_1d_2}g\circ f\,,\end{equation*}
which, by \ref{graded spaces properties}~(i), is homogeneous of degree
$d_1+d_2$, is called
the \emph{graded bracket} of $f$ and $g$.
\end{defi}

\begin{defi}
\label{def graded derivation} 
Let $A=\bigoplus_{p\in\ZZ}A^p$ be a $\ZZ$-graded
algebra. Let $\theta:A\to A$ be a linear endomorphism of the
graded vector space $A$. Let $d\in\ZZ$. The linear endomorphism
$\theta$ is said to be a {\it derivation of degree $d$\/} of the
graded algebra $A$ if
\begin{description}
\item{\rm(i)} as a linear endomorphism of a graded vector space,
$\theta$ is homogeneous of degree $d$,

\item{\rm(ii)} for all $p\in\ZZ$, $x\in A^p$ and $y\in
A$,
 \begin{equation*}\theta(xy)=
 \bigl(\theta(x)\bigr)y+(-1)^{dp}x\bigl(\theta(y)\bigr)\,.\end{equation*}
\end{description}
\end{defi}

\begin{rmk} More generally, as shown by  Koszul \cite{Ko1}, for an
algebra $A$ equip\-ped with an involutive automorphism, one can
define two types of remarkable linear endomorphisms of $A$, the
{\it derivations\/} and the {\it antiderivations\/}. When
$A=\bigoplus_{p\in\ZZ}A^p$ is a $\ZZ$-graded algebra, and when the
involutive automorphism used is that which maps each $x\in A^p$
onto $(-1)^px$, it turns out that all nonzero graded derivations
are of even degree, that all nonzero graded antiderivations are of
odd degree, and that both derivations and antiderivations can be
defined as done in Definition~\ref{def graded derivation}. 
For simplicity we have chosen
to call {\it derivations\/} both the derivations and
antiderivations.
\end{rmk}

\subsection{Some~properties~of~derivations}
\label{properties of derivations}
Let $A=\bigoplus_{p\in\ZZ}A^p$ be a $\ZZ$-graded algebra.
\par\smallskip\nobreak\noindent
{\rm(i)}\quad{\it A derivation of degree $0$\/.}\quad
 For every
$p\in\ZZ$ and $x\in A^p$, we set
 \begin{equation*}\mu(x)=px\,.\end{equation*}
The map $\mu$, defined for homogeneous elements of $A$, can be
extended in a unique way as a linear endomorphism of $A$, still
denoted by $\mu$. This endomorphism is a derivation of degree $0$
of $A$.
\par\smallskip\nobreak\noindent
{\rm(ii)}\quad{\it The graded bracket of two derivations\,.}\quad
Let $\theta_1:A\to A$ and $\theta_2:A\to A$ be two derivations of
$A$, of degree $d_1$ and $d_2$, respectively. Their graded
bracket (Definition \ref{graded bracket})
 \begin{equation*}[\theta_1,\theta_2]
 =\theta_1\circ\theta_2-(-1)^{d_1d_2}
 \theta_2\circ\theta_1\,,\end{equation*} 
is a derivation of degree $d_1+d_2$.

\begin{defi}
\label{def graded Lie algebra}
A {\it $\ZZ$-graded Lie algebra\/} is a $\ZZ$-graded algebra
$A=\bigoplus_{p\in\ZZ}A^p$ (in the sense of \ref{graded definitions}~(v)), whose
composition law, often denoted by $(x,y)\mapsto [x,y]$ and called
the {\it graded bracket\/}, satisfies the following two
properties:
\par\smallskip
\begin{description}
\item{\rm(i)} it is $\ZZ_2$-anticommutative in the sense of \ref{graded
definitions}~(vi),
{\it i.e.}, for all $p$ and $q\in\ZZ$, $P\in A^p$ and $Q\in A^q$,
 \begin{equation*}[P,Q]=-(-1)^{pq}[Q,P]\,,\end{equation*}

\item{\rm(ii)} it satisfies the {\it $\ZZ$-graded Jacobi identity\/}, {\it
i.e.}, for $p$, $q$ and $r\in\ZZ$, $P\in A^p$, $Q\in A^q$ and
$R\in A^r$,
 \begin{equation*}(-1)^{pr}\bigl[P,[Q,R]\bigr]+(-1)^{qp}\bigl[Q,[R,P]\bigr]
 +(-1)^{rq}\bigl[R,[P,Q]\bigr]=0\,.\end{equation*}
\end{description}
\end{defi}

\subsection{Examples~and~remarks}
\label{examples graded Lie algebras}
\par\smallskip\nobreak\noindent
{\rm(i)}\quad{\it Lie algebras and $\ZZ$-graded Lie
algebras\/.}\quad A $\ZZ$-graded Lie algebra
$A=\bigoplus_{p\in\ZZ}A^p$ is not a Lie algebra in the usual sense,
unless $A^p=\{0\}$ for all $p\neq 0$. However, its subspace $A^0$
of homogeneous elements of degree $0$ is a Lie algebra in that
usual sense: it is stable under the bracket operation and when
restricted to elements in $A^0$, the bracket is skew-symmetric and
satisfies the usual Jacobi identity.

\par\smallskip\nobreak\noindent
{\rm(ii)}\quad{\it The graded Lie algebra associated to a graded
associative algebra\/.}\quad Let $A=\bigoplus_{p\in\ZZ}A^p$ be a
$\ZZ$-graded associative algebra, whose composition law is denoted
by $(P,Q)\mapsto PQ$. We define another composition law, denoted
by $(P,Q)\mapsto[P,Q]$ and called the {\it graded commutator\/};
we first define it for homogeneous elements in $A$ by setting, for
all $p$ and $q\in\ZZ$, $P\in A^p$ and $Q\in A^q$,
 \begin{equation*}[P,Q]=PQ-(-1)^{pq}QP\,;\end{equation*}
then we extend the definition of that composition law to all pairs
of elements in $A$ by bilinearity. The reader will easily verify
that with this composition law, $A$ is a graded Lie algebra. When
$A^p=\{0\}$ for all $p\neq 0$, we recover the well known way in
which one can associate a Lie algebra to any associative algebra.

\par\smallskip\nobreak\noindent
{\rm(iii)}\quad{\it The graded Lie algebra of graded
endomorphisms\/.}\quad Let $E=\bigoplus_{p\in\ZZ}E^p$ be a graded
vector space. For each $p\in\ZZ$, let $A^p\subset{\cal L}(E,E)$ be
the space of linear endomorphisms of $E$ which are homogeneous of
degree $p$, and let $A=\bigoplus_{p\in \ZZ}A^p$. As we have seen
in~\ref{graded spaces properties}~(iii), when equipped with the composition of
applications
as composition law, $A$ is a $\ZZ$-graded associative algebra. Let
us define another composition law on $A$, called the {\it graded
commutator\/}; we first define it for homogeneous elements in $A$
by setting, for all $p$ and $q\in\ZZ$, $P\in A^p$ and $Q\in A^q$,
 \begin{equation*}[P,Q]=PQ-(-1)^{pq}QP\,;\end{equation*}
then we extend the definition of that composition law to all pairs
of elements in $A$ by bilinearity. By using \ref{examples graded Lie
algebras}~(ii), we see
that $A$, with this composition law, is a $\ZZ$-graded Lie
algebra.

\par\smallskip\nobreak\noindent
{\rm(iv)}\quad{\it Various interpretations of the graded Jacobi
identity\/.}\quad Let $A=\bigoplus_{p\in\ZZ}A^p$ be a $\ZZ$-graded
Lie algebra. The $\ZZ$-graded Jacobi identity indicated in
Definition~\ref{def graded Lie algebra} can be cast into other forms, which
better
indicate its meaning. Let us set, for all $P$ and $Q\in A$,
 \begin{equation*}\ad_PQ=[P,Q]\,.\end{equation*}
For each $p\in\ZZ$ and $P\in A^p$, $\ad_P:A\to A$ is a graded
endomorphism of $A$, homogeneous of degree $p$. By taking into
account the $\ZZ_2$-anticommutativity of the bracket, the reader
will easily see that the graded Jacobi identity can be written
under the following two forms:
\par\smallskip
{\it First form\/.}\quad For all $p$, $q$ and $r\in\ZZ$, $P\in
A^p$, $Q\in A^q$ and $R\in A^r$,
 \begin{equation*}\ad_P\bigl([Q,R]\bigr)
 =[\ad_PQ,R]+(-1)^{pq}[Q,\ad_PR]\,.
 \end{equation*}
This equality means that for all $p\in \ZZ$ and $P\in A^p$, the
linear endomorphism $\ad_P:A\to A$ is  a derivation of degree $p$
of the graded Lie algebra $A$, in the sense of \ref{def graded derivation}.
\par\smallskip
{\it Second form\/.}\quad For all $p$, $q$ and $r\in \ZZ$, $P\in
A^p$, $Q\in A^q$ and $R\in A^r$,
 \begin{equation*}\ad_{[P,Q]}R=\ad_P\circ\ad_QR-(-1)^{pq}\ad_Q\circ\ad_PR=[\ad_P
 ,\ad_Q]R\,.\end{equation*} 
This equality means that for all $p$
and $q\in\ZZ$, $P\in A^p$ and $Q\in A^q$, the endomorphism
$\ad_{[P,Q]}:A\to A$ is the graded bracket (in the sense of \ref{graded
bracket})
of the two endomorphisms $\ad_P:A\to A$ and $\ad_Q:A\to A$. In
other words, the map $P\mapsto \ad_P$ is a $\ZZ$-graded Lie
algebras homomorphism from the $\ZZ$-graded Lie algebra $A$ into
the $\ZZ$-graded Lie algebra of sums of linear homogeneous
endomorphisms of $A$, with the graded bracket as composition law
(example \ref{examples graded Lie algebras}~(iii)).
\par\smallskip
When $A^p=\{0\}$ for all $p\neq 0$, we recover the well known
interpretations of the usual Jacobi identity.

\section{Exterior powers of a vector bundle and of its dual}

In what follows all the vector bundles will be assumed to be
locally trivial and of finite rank; therefore we will write simply
{\it vector bundle\/} for {\it locally trivial vector bundle\/}.

\subsection{The~dual~of~a~vector~bundle} Let $(E,\tau,M)$ be a
vector bundle on the field $\KK=\RR$ or $\CC$. We will denote its
{\it dual bundle\/} by $(E^*,\pi,M)$. Let us recall that it is
a vector bundle over the same base manifold $M$, whose fibre
$E^*_x=\pi^{-1}(x)$ over each point $x\in M$ is the dual vector
space of the corresponding fibre $E_x=\tau^{-1}(x)$ of $(E,\tau,M)$,
{\it i.e.}, the space of linear forms on $E_x$ ({\it i.e.}, linear
functions defined on $E_x$ and taking their values in the field
$\KK$).
\par\smallskip
For each $x\in M$, the duality coupling $E^*_x\times E_x\to\KK$
will be denoted by
 \begin{equation*}(\eta,v)\mapsto\langle\eta,v\rangle\,.\end{equation*}

\subsection{The~exterior~powers~of~a~vector~bundle} Let
$(E,\tau,M)$ be a vector bundle of rank $k$. For each integer
$p>0$, we will denote by $(\bigwedge^pE,\tau,M)$ the $p$-th
external power of $(E,\tau,M)$. It is a vector bundle over $M$
whose fibre $\bigwedge^pE_x$, over each point $x\in M$, is the
$p$-th external power of the corresponding fibre $E_x=\tau^{-1}(x)$
of $(E,\tau,M)$. We recall that $\bigwedge^pE_x$ can be canonically
identified with the vector space of $p$-multilinear skew-symmetric
forms on the dual $E^*_x$ of $E_x$.
\par\smallskip
Similarly, for any integer $p>0$, we will denote by
$(\bigwedge^pE^*,\pi,M)$ the $p$-th external power of the
bundle $(E^*,\pi,M)$, dual of $(E,\tau,M)$.
\par\smallskip
For $p=1$, $(\bigwedge^1E,\tau,M)$ is simply the bundle
$(E,\tau,M)$, and similarly $(\bigwedge^1E^*,\pi,M)$ is simply
the bundle $(E^*,\pi,M)$. For $p$ strictly larger than the rank
$k$ of $(E,\tau,M)$, $(\bigwedge^pE,\tau,M)$ and
$(\bigwedge^pE^*,\pi,M)$ are the trivial bundle over $M$,
$(M\times\{0\},p_1,M)$, whose fibres are zero-dimensional
($p_1:M\times\{0\}\to M$ being the projection onto the first
factor).
\par\smallskip
For $p=0$, we will set
 $(\bigwedge^0E,\tau,M)=(\bigwedge^0E^*,\pi,M)=(M\times\KK,p_1,M)$,
where $p_1:M\times \KK\to M$ is the projection onto the first
factor.
\par\smallskip
Finally, we will consider that for $p<0$, $(\bigwedge^pE,\tau,M)$
and $(\bigwedge^pE^*,\pi,M)$ are the trivial bundle over $M$,
$(M\times\{0\},p_1,M)$. With these conventions,
$(\bigwedge^pE,\tau,M)$ and $(\bigwedge^pE^*,\pi,M)$ are defined
for all $p\in\ZZ$.
\par\medskip

\subsection{Operations in the graded vector spaces $\bigwedge
E_x$ and $\bigwedge E^*_x$}  
\label{operations in exterior powers}
Let $(E,\tau,M)$ be a vector bundle of
rank $k$, $(E^*,\pi,M)$ its dual and, for each $p\in\ZZ$,
$(\bigwedge^pE,\tau,M)$ and $(\bigwedge^pE^*,\pi,M)$ their
$p$-th external powers. We recall in this sections some operations
which can be made for each point $x\in M$, in the vector spaces
$\bigwedge^pE_x$ and $\bigwedge^pE^*_x$.
\par\smallskip
For each $x\in M$, let us consider the $\ZZ$-graded vector spaces
 \begin{equation*}\bigwedge
 E_x=\bigoplus_{p\in\ZZ}\bigwedge^pE_x\quad\hbox{and}\quad
 \bigwedge E^*_x=\bigoplus_{p\in\ZZ}\bigwedge^pE^*_x\,.\end{equation*}
We will say that elements in $\bigwedge E^*_x$ are (multilinear)
{\it forms\/} at $x$, and that elements in $\bigwedge E_x$ are
{\it multivectors\/} at $x$.
\par\smallskip\nobreak\noindent
{\rm(i)}\quad{\it The exterior product\/.}\quad Let us recall that
for each $x\in M$, $p$ and $q\in\ZZ$, $P\in\bigwedge^pE_x$ and
$Q\in\bigwedge^qE_x$, there exists $P\wedge
Q\in\bigwedge^{p+q}E_x$, called the {\it exterior product\/} of
$P$ and $Q$, defined by the following formulae.
\par\smallskip
\begin{description}
\item{--} If $p<0$, then $P=0$, therefore, for any
$Q\in\bigwedge^qE_x$, $P\wedge Q=0$. Similarly, if $q<0$, then
$Q=0$, therefore, for any $P\in\bigwedge^pE_x$, $P\wedge Q=0$.
\par\smallskip
\item{--} If $p=0$, then $P$ is a scalar ($P\in\KK$), and
therefore, for any $Q\in\bigwedge^qE_x$, $P\wedge Q=PQ$, the usual
product of $Q$ by the scalar $P$. Similarly, for $q=0$, then $Q$
is a scalar $(Q\in\KK)$, and therefore, for any
$P\in\bigwedge^pE_x$, $P\wedge Q=QP$, the usual product of $P$ by
the scalar $Q$.
\par\smallskip
\item{--} If $p\geq 1$ and $q\geq 1$, $P\wedge Q$, considered
as a $(p+q)$-multilinear form on $E^*_x$, is given by the formula,
where $\eta_1,\ldots,\eta_{p+q}\in E^*_x$,
\end{description}
 \begin{equation*}P\wedge Q(\eta_1,\ldots,\eta_{p+q})=\sum_{\sigma\in{\cal
 S}_{(p,q)}}\varepsilon(\sigma)
 P(\eta_{\sigma(1)},\ldots,\eta_{\sigma(p)})Q(\eta_{\sigma(p+1)},\ldots,
 \eta_{\sigma(p+q)})\,.\end{equation*}
We have denoted by ${\cal S}_{(p,q)}$ the set of permutations
$\sigma$ of $\{\,1,\ldots,p+q\,\}$ which satisfy
 \begin{equation*}\sigma(1)<\sigma(2)<\cdots<\sigma(p)\quad\hbox{and}\quad
 \sigma(p+1)<\sigma(p+2)<\cdots<\sigma(p+q)\,,\end{equation*}
and set
 \begin{equation*}\epsilon(\sigma)=\begin{cases}
                    1&\text{if $\sigma$ is even},\\
                    -1&\text{if $\sigma$ is odd}.
                    \end{cases}\end{equation*}
\par\smallskip
Similarly, let us recall that for each $x\in M$, $p$ and
$q\in\ZZ$, $\xi\in\bigwedge^pE^*_x$ and $\eta\in\bigwedge^qE^*_x$,
there exists $\xi\wedge \eta\in\bigwedge^{p+q}E^*_x$, called the
{\it exterior product\/} of $\xi$ and $\eta$. It is defined by the
formulae given above, the only change being the exchange of the
roles of $E_x$ and $E^*_x$.
\par\smallskip
The exterior product is associative and $\ZZ_2$-commutative: for
all $x\in M$, $p$, $q$ and $r\in\ZZ$, $P\in \bigwedge^pE_x$, $Q\in
\bigwedge^qE_x$ and $R\in \bigwedge^rE_x$,
 \begin{equation*}P\wedge(Q\wedge R)=(P\wedge Q)\wedge R\,,\qquad Q\wedge
 P=(-1)^{pq}P\wedge Q\,,\end{equation*}
 and similarly, for $\xi\in\bigwedge^pE^*_x$, $\eta\in
\bigwedge^qE^*_x$ and $\zeta\in \bigwedge^rE^*_x$,
 \begin{equation*}\xi\wedge(\eta\wedge \zeta)=(\xi\wedge \eta)\wedge
\zeta\,,\qquad \eta\wedge
 \xi=(-1)^{pq}\xi\wedge \eta\,.\end{equation*}
\par\medskip
For all $x\in M$, the exterior product extends, by bilinearity, as
a composition law in each of the graded vector spaces $\bigwedge
E_x$ and $\bigwedge E^*_x$. With these composition laws, these
vector spaces become $\ZZ$-graded associative and
$\ZZ_2$-commutative algebras.
\par\smallskip\nobreak\noindent
{\rm(ii)}\quad{\it The interior product of a form by a
vector\/.}\quad Let us recall that for each $x\in M$, $v\in E_x$,
$p\in\ZZ$, $\eta\in\bigwedge^pE^*_x$, there exists $i(v)\eta\in
\bigwedge^{p-1}E^*_x$, called the {\it interior product\/} of
$\eta$ by $v$, defined by the following formulae.
\par\smallskip
\begin{description}
\item{--} For $p\leq 0$, $i(v)\eta=0$, since
$\bigwedge^{p-1}E^*_x=\{0\}$.
\par\smallskip
\item{--} For $p=1$,
 \begin{equation*}i(v)\eta=\langle\eta,v\rangle\in\KK\,.\end{equation*}

\par\smallskip
\item{--} For $p>1$, $i(v)\eta$ is the $(p-1)$-multilinear form
on $E_x$ such that, for all $v_1,\ldots,v_{p-1}\in E_x$,
 \begin{equation*}i(v)\eta(v_1,\ldots,v_{p-1})=\eta(v,v_1,\ldots,v_{p-1})\,.
\end{equation*}
\end{description}
\par\medskip
For each $x\in M$ and $v\in E_x$, the map $\eta\mapsto i(v)\eta$
extends, by linearity, as a graded endomorphism of degree $-1$ of
the graded vector space $\bigwedge E^*_x$. Moreover, that
endomorphism is in fact a derivation of degree $-1$ of the
exterior algebra of $E^*_x$, {\it i.e.}, for all $p$ and
$q\in\ZZ$, $\zeta\in \bigwedge^pE^*_x$, $\eta\in\bigwedge^qE^*_x$,
 \begin{equation*}i(v)(\zeta\wedge\eta)
 =\bigl(i(v)\zeta\bigr)\wedge\eta+(-1)^p
 \zeta\wedge\bigl(i(v)\eta\bigr)\,.
 \end{equation*}
\par\smallskip\nobreak\noindent
{\rm(iii)}\quad{\it The pairing between $\bigwedge E_x$ and
$\bigwedge E^*_x$\/.}\quad Let $x\in M$, $p$ and $q\in \ZZ$,
$\eta\in\bigwedge^pE^*_x$ and $v\in\bigwedge^qE_x$. We set
 \begin{equation*}\langle\eta,v\rangle=\begin{cases}
                        0&\text{if $p\neq q$, or if $p<0$, or if
                        $q<0$},\\
                        \eta v&\text{if $p=q=0$}.
                        \end{cases}
 \end{equation*}
In order to define $\langle\eta,v\rangle$ when $p=q\geq 1$, let us
first assume that $\eta$ and $v$ are decomposable, {\it i.e.},
that they can be written as
 \begin{equation*}\eta=\eta_1\wedge\cdots\wedge\eta_p\,,\qquad
 v=v_1\wedge\cdots\wedge v_p\,,\end{equation*}
where $\eta_i\in E^*_x$, $v_j\in E_x$, $1\leq i,j\leq p$. Then we
set
 \begin{equation*}\langle\eta,v\rangle=\det\bigl(\langle\eta_i,v_j\rangle\bigr)\
,.\end{equation*} One may see that $\langle\eta,v\rangle$ depends
only on $\eta$ and $v$, not on the way in which they are expressed
as exterior products of elements of degree $1$. The map
$(\eta,v)\mapsto\langle\eta,v\rangle$ extends, in a unique way as
a bilinear map
 \begin{equation*}\bigwedge E^*_x\times\bigwedge E_x\to\KK,\quad\hbox{still
 denoted by}\quad
 (\eta,v)\mapsto\langle\eta,v\rangle\,,\end{equation*}
called the {\it pairing\/}. That map allows us to consider each
one of the two graded vector spaces $\bigwedge E^*_x$ and
$\bigwedge E_x$ as the dual of the other one.
\par\smallskip
Let $\eta\in\bigwedge^pE^*_x$ and  $v_1,\ldots,v_p$ be elements of
$E_x$. The pairing $\langle\eta,v_1\wedge\cdots\wedge v_p\rangle$
is related, in a very simple way, to the value taken by $\eta$,
considered as a $p$-multilinear form on $E_x$,  on the set
$(v_1,\ldots,v_p)$. We have
 \begin{equation*}\langle\eta,v_1\wedge\cdots\wedge
 v_p\rangle=\eta(v_1,\ldots,v_p)\,.\end{equation*}
\par\smallskip\nobreak\noindent
{\rm(iv)}\quad{\it The interior product of a form by a
multivector\/.}\quad For each $x\in M$ and $v\in E_x$,  we have
defined in \ref{operations in exterior powers}~(ii) the interior product
$i(v)$ as a derivation
of degree $-1$ of the exterior algebra $\bigwedge E^*_x$ of forms
at $x$. Let us now define, for each multivector $P\in\bigwedge
E_x$, the interior product $i(P)$. Let us first assume that $P$ is
homogeneous of degree $p$, {\it i.e.}, that $P\in\bigwedge^pE_x$.
\par\smallskip
\begin{description}
\item{--} For $p<0$, $\bigwedge^pE_x=\{0\}$, therefore $i(P)=0$.
\par\smallskip
\item{--} For $p=0$, $\bigwedge^0E_x=\KK$, therefore $P$ is a
scalar and we set, for all $\eta\in\bigwedge E^*_x$,
 \begin{equation*}i(P)\eta=P\eta\,.\end{equation*}
\par\smallskip
\item{--} For $p\geq 1$ and $P\in\bigwedge^pE_x$ decomposable,
{\it i.e.},
 \begin{equation*}P=P_1\wedge \cdots\wedge P_p\,,\quad \hbox{with}\quad
P_i\in E_x\,,
 \quad 1\leq i\leq p\,,\end{equation*}
we set
 \begin{equation*}i(P_1\wedge\cdots\wedge P_p)=i(P_1)\circ\cdots\circ
i(P_p)\,.\end{equation*} We see easily that $i(P)$ depends only of
$P$, not of the way in which it is expressed as an exterior
product of elements of degree $1$.
\par\smallskip
\item{--} We extend by linearity the definition of
$i(P)$ for all $P\in\bigwedge^pE_x$, and we see that $i(P)$ is a
graded endomorphism of degree $-p$ of the graded vector space
$\bigwedge E^*_x$. Observe that for $p\neq 1$, $i(P)$ is not in
general a derivation of the exterior algebra $\bigwedge E^*_x$.
\end{description}
\par\smallskip
Finally, we extend by linearity the definition of $i(P)$ to all
elements $P\in\bigwedge E_x$.
\par\smallskip\nobreak\noindent
{\rm(v)}\quad{\it The interior product by an exterior
product\/.}\quad It is easy to see that for all $P$ and
$Q\in\bigwedge E_x$,
 \begin{equation*}i(P\wedge Q)=i(P)\circ i(Q)\,.\end{equation*}
\par\smallskip\nobreak\noindent
{\rm(vi)}\quad{\it Interior product and pairing\/.}\quad For
$p\in\ZZ$, $\eta\in\bigwedge^{p}E^*_x$ and $P\in\bigwedge^pE_x$,
we have
 \begin{equation*}i(P)\eta=(-1)^{(p-1)p/2}\langle\eta,P\rangle\,.\end{equation*}
More generally, for $p$ and $q\in\ZZ$, $P\in\bigwedge^p(E_x)$,
$Q\in\bigwedge^q(E_x)$ and  $\eta\in\bigwedge^{p+q}(E^*_x)$,
 \begin{equation*}\bigl\langle
i(P)\eta,Q\bigr\rangle=(-1)^{(p-1)p/2}\langle\eta,
 P\wedge Q\rangle\,.
 \end{equation*}
This formula shows that the interior product by
$P\in\bigwedge^pE_x$ is $(-1)^{(p-1)p/2}$ times the transpose,
with respect to the pairing, of the exterior product by $P$ on the
left.

\subsection{The~exterior~algebra~of~sections}
\label{The exterior algebra of sections} 
Let $(E,\tau,M)$ be
a vector bundle of rank $k$ on the field $\KK=\RR$ or $\CC$, over
a smooth manifold $M$, $(E^*,\pi,M)$ be its dual bundle and,
for each integer $p\geq 1$, let $(\bigwedge^pE,\tau,M)$ and
$(\bigwedge^pE^*,\pi,M)$ be their respective $p$-th exterior
powers.
\par\smallskip
For each $p\in\ZZ$, we will denote by $A^p(M,E)$ the space of
smooth sections of $(\bigwedge^pE,\tau,M)$, {\it i.e.}, the space
of smooth maps $Z:M\to\bigwedge^pE$ which satisfy
 \begin{equation*}\tau\circ Z=\id_M\,.\end{equation*}
Similarly, for each $p\in\ZZ$, we will denote by $\Omega^p(M,E)$
the space of smooth sections of the vector bundle $(\bigwedge^pE^*,\pi,M)$, {\it
i.e.}, the space of smooth maps $\eta:M\to\bigwedge^pE^*$ which
satisfy
 \begin{equation*}\pi\circ \eta=\id_M\,.\end{equation*}
Let us observe that $\Omega^p(M,E)=A^p(M,E^*)$.
\par
We will denote by $A(M,E)$ and $\Omega(M,E)$ the direct sums
 \begin{equation*}A(M,E)=\bigoplus_{p\in\ZZ}A^p(M,E)\,,\qquad
   \Omega(M,E)=\bigoplus_{p\in\ZZ}\Omega^p(M,E)\,.\end{equation*}
These direct sums, taken for all $p\in\ZZ$, are in fact taken for
all integers $p$ which satisfy $0\leq p\leq k$, where $k$ is the
rank of the vector bundle $(E,\tau,M)$, since we have
$A^p(M,E)=\Omega^p(M,E)=\{0\}$ for $p<0$ as well as for $p>k$.
\par\smallskip
For $p=0$, $A^0(M,E)$ and $\Omega^0(M,E)$ both coincide with the
space $C^\infty(M,\KK)$ of smooth functions defined on $M$ which
take their values in the field $\KK$.
\par\smallskip
Operations such as the exterior product, the interior product and
the pairing, defined for each point $x\in M$ 
in~\ref{operations in exterior powers}, can be
extended to elements in $A(M,E)$ and $\Omega(M,E)$.
\par\smallskip\nobreak\noindent
{\rm(i)}\quad{\it The exterior product of two sections\/.}\quad
For example, the exterior product of two sections $P$ and $Q\in
A(M,E)$ is the section
 \begin{equation*}x\in M\,,\quad x\mapsto (P\wedge Q)(x)=P(x)\wedge
Q(x)\,.\end{equation*} The exterior product of two sections $\eta$
and $\zeta\in\Omega(M,E)$ is similarly defined.
\par\smallskip
With the exterior product as composition law, $A(M,E)$ and
$\Omega(M,E)$ are $\ZZ$-graded associative and $\ZZ_2$-commutative
algebras, called the {\it algebra of multivectors\/} and the {\it
algebra of forms\/} associated to the vector bundle $(E,\tau,M)$.
Their subspaces $A^0(M,E)$ and $\Omega^0(M,E)$ of homogeneous
elements of degree $0$ both coincide with the usual algebra
$C^\infty(M,\KK)$ of smooth $\KK$-valued functions on $M$, with
the usual product of functions as composition law. We observe that
$A(M,E)$ and $\Omega(M,E)$ are $\ZZ$-graded modules over the ring
of functions $C^\infty(M,\KK)$.
\par\smallskip\nobreak\noindent
{\rm(ii)}\quad{\it The interior product by a section of
$A(M,E)$\/.}\quad For each $P\in A(M,E)$, the {\it interior
product\/} $i(P)$ is an endomorphism of the graded vector space
$\Omega(M,E)$. If $p\in\ZZ$ and $P\in A^p(M,E)$, the endomorphism
$i(P)$ is homogeneous of degree $-p$. For $p=1$, $i(P)$ is a
derivation of degree $-1$ of the algebra $\Omega(M,E)$.
\par\smallskip\nobreak\noindent
{\rm(iii)}\quad{\it The pairing between $A(M,E)$ and
$\Omega(M,E)$\/.}\quad The {\it pairing\/}
 \begin{equation*}(\eta,P)\mapsto\langle\eta,P\rangle\,,\quad \eta\in
 \Omega(M,E)\,,\quad P\in A(M,E)\,,\end{equation*}
is a $C^\infty(M,\KK)$-bilinear map, defined on $\Omega(M,E)\times
A(M,E)$, which takes its values in $C^\infty(M,\KK)$. \

\chapter{Exterior powers of a Lie algebroid and of its dual}
\label{Exterior powers 2}
We consider now a Lie algebroid
$(E,\tau,M,\rho)$ over a smooth manifold $M$. We denote by
$(E^*,\pi,M)$ its dual vector bundle, and use all the notations
defined in Section~ \ref{Exterior powers}. We will assume that the base field
$\KK$ is
$\RR$, but most results remain valid for $\KK=\CC$. We will prove
that differential operators such as the Lie derivative and the
exterior derivative, which are well known for sections of the
exterior powers of a tangent bundle and of its dual, still exist
in this more general setting.

\section{Lie derivatives with respect to sections of a Lie
algebroid}
\label{Lie derivative by section of algebroid}
We prove in this section that for each smooth section $V$ of
the Lie algebroid $(E,\tau,M,\rho)$, there exists a derivation of
degree $0$ of the exterior algebra $\Omega(M,E)$, called the {\it
Lie derivative\/} with respect to $V$ and denoted by ${\cal
L}_\rho(V)$. When the Lie algebroid is the tangent bundle
$(TM,\tau_M,M,\id_{TM})$, we will recover the usual Lie derivative
of differential forms with respect to a vector field.

\begin{prop}
\label{Lie derivative by section of algebroid 2}
Let $(E,\tau,M,\rho)$ be a Lie algebroid on a smooth
manifold $M$. For each smooth section $V\in A^1(M,E)$ of the
vector bundle $(E,\tau,M)$, there exists a unique graded
endomorphism of degree $0$ of the graded algebra of exterior forms
$\Omega(M,E)$, called the {\it Lie derivative\/} with respect to
$V$ and denoted by ${\cal L}_\rho(V)$, which satisfies the
following properties:
\begin{description}
\par\smallskip
\item{\rm(i)} For a smooth function
$f\in\Omega^0(M,E)=C^\infty(M,\RR)$,
 \begin{equation*}{\cal L}_\rho(V)f=i(\rho\circ V)df={\cal L}(\rho\circ
V)f\,,\end{equation*} where ${\cal L}(\rho\circ V)$ denotes the
usual Lie derivative with respect to the vector field $\rho\circ
V$;
\par\smallskip
\item{\rm (ii)} For a form $\eta\in\Omega^p(M,E)$ of degree $p>0$,
${\cal L}_\rho(V)\eta$ is the form defined by the formula, where
$V_1,\ldots,V_p$ are smooth sections of $(E,\tau,M)$,
 \begin{equation*}
 \begin{split}
 \bigl({\cal L}_\rho(V)\eta\bigr)(V_1,\ldots,V_p)
 &={\cal L}_\rho(V)\bigl(\eta(V_1,\ldots,V_p)\bigr)\\
 &\quad-\sum_{i=1}^p\eta(V_1,\ldots,V_{i-1},\{V,V_i\},V_{i+1},\ldots,V_p)\,.
 \end{split}
 \end{equation*}
\par\smallskip
\end{description}
\end{prop}

\begin{proof} Clearly (i) defines a function ${\cal L}_\rho(V)f\in
\Omega^0(M,E)=C^\infty(M,\RR)$. We see immediately that for $f$
and $g\in C^\infty(M,\RR)$,
 \begin{equation*}{\cal L}_\rho(V)(fg)=\bigl({\cal L}_\rho(V)f\bigr)g
 +f\bigl({\cal L}_\rho(V)g\bigr)\,.\eqno(*)\end{equation*}
Now (ii) defines a map $(V_1,\ldots V_p)\mapsto\bigl({\cal
L}_\rho(V)\eta\bigr)(V_1,\ldots,V_p)$ on $\bigl(A^1(M,E)\bigr)^p$,
with values in $C^\infty(M,\RR)$. In order to prove that this map
defines an element ${\cal L}_\rho(V)\eta$ in $\Omega^p(M,E)$, it
is enough to prove that it is skew-symmetric and
$C^\infty(M,\RR)$-linear in each argument. The skew-symmetry and
the $\RR$-linearity in each argument are easily verified. There
remains only to prove that for each function $f\in
C^\infty(M,\RR)$,
 \begin{equation*}\bigl({\cal L}_\rho(V)\eta\bigr)(fV_1,V_2,\ldots,V_p)=f
 \bigl({\cal
L}_\rho(V)\eta\bigr)(V_1,V_2,\ldots,V_p)\,.\eqno(**)\end{equation*}
We have
 \begin{equation*}
 \begin{split}
 \bigl({\cal L}_\rho(V)\eta\bigr)(fV_1,V_2,\ldots,V_p)
 &={\cal L}_\rho(V)\bigl(\eta(fV_1,V_2,\ldots,V_p)\bigr)\\
 &\quad-\eta(\{V,fV_1\},V_2,\ldots,V_p)\\
 &\quad-\sum_{i=2}^p\eta(fV_1,V_2,\ldots,V_{i-1},\{V,V_i\},V_{i+1},\ldots
 ,V_p)\,.
 \end{split}
 \end{equation*}
By using $(*)$, we may write
 \begin{equation*}
 \begin{split}
 {\cal L}_\rho(V)\bigl(\eta(fV_1,V_2,\ldots,V_p)\bigr)
 &={\cal L}_\rho(V)\bigl(f\eta(V_1,V_2,\ldots,V_p)\bigr)\\
 &=\bigl({\cal L}_\rho(V)f\bigr)\eta(V_1,V_2,\ldots,V_p)\\
 &\quad+f{\cal
  L}_\rho(V)\bigl(\eta(V_1,V_2,\ldots,V_p)\bigr)\,.
 \end{split}
  \end{equation*}
Using the property of the anchor, we also have
 \begin{equation*}\{V,fV_1\}=\bigl(i(\rho\circ
V)df\bigr)V_1+f\{V,V_1\}=\bigl({\cal
 L}_\rho(V)f\bigr)V_1+f\{V,V_1\}\,.\end{equation*}
Equality $(**)$ follows immediately.
\par
The endomorphism ${\cal L}_\rho(V)$, defined on the subspaces of
homogeneous forms, can then be extended, in a unique way, to
$\Omega(M,E)$, by imposing the $\RR$-linearity of the  map
$\eta\mapsto {\cal L}_\rho(V)\eta$.
\end{proof}
\par\bigskip

Let us now introduce the $\Omega(M,E)$-valued exterior derivative
of a function. In the next section, that definition will be
extended to all elements in $\Omega(M,E)$.

\begin{defi}
\label{exterior derivative of functions} 
Let $(E,\tau,M,\rho)$ be a Lie algebroid
on a
smooth manifold $M$. For each function
$f\in\Omega^0(M,E)=C^\infty(M,\RR)$, we call {\it
$\Omega(M,E)$-valued exterior derivative\/} of $f$, and denote by
$d_\rho f$, the unique element in $\Omega^1(M,E)$ such that, for
each section $V\in A^1(M,E)$,
 \begin{equation*}\langle d_\rho f,V\rangle=\langle df,\rho\circ
V\rangle\,.\end{equation*}

\end{defi}

\begin{rmk} Let us observe that the transpose of the anchor
$\rho:E\to TM$ is a vector bundle map ${}^t\!\rho:T^*M\to E^*$. By
composition of that map with the usual differential of functions,
we obtain the $\Omega(M,E)$-valued exterior differential $d_\rho$.
We have indeed
 \begin{equation*}d_\rho f={}^t\!\rho\circ df\,.\end{equation*}
\end{rmk}

\begin{prop}
\label{Lie derivative by section of algebroid 3}
Under the assumptions of
Proposition~\ref{Lie derivative by section of algebroid 2}, 
the Lie derivative has the following properties.
\par\smallskip
{\rm 1.} For each $V\in A^1(M,E)$ and $f\in C^\infty(M,\RR)$,
 \begin{equation*}{\cal L}_\rho(V)(d_\rho f)=d_\rho\bigl({\cal
 L}_\rho(V)f\bigr)\,.\end{equation*}
\par\smallskip
{\rm 2.} For each $V$ and $W\in A^1(M,E)$, $\eta\in\Omega(M,E)$,
 \begin{equation*}i\bigl(\{V,W\}\bigr)\eta=\bigl({\cal L}_\rho(V)\circ i(W)-
 i(W)\circ{\cal L}_\rho(V)\bigr)\eta\,.\end{equation*}
\par\smallskip
{\rm 3.} For each $V\in A^1(M,E)$, ${\cal L}_\rho(V)$ is a
derivation of degree $0$ of the exterior algebra $\Omega(M,E)$.
That means that for all $\eta$ and $\zeta\in\Omega(M,E)$,
 \begin{equation*}{\cal L}_\rho(V)(\eta\wedge\zeta)=\bigl({\cal
 L}_\rho(V)\eta\bigr)\wedge\zeta+\eta\wedge\bigl({\cal
 L}_\rho(V)\zeta\bigr)\,.\end{equation*}
\par\smallskip
{\rm 4.}  For each $V$ and $W\in A^1(M,E)$, $\eta\in\Omega(M,E)$,
 \begin{equation*}{\cal L}_\rho\bigl(\{V,W\}\bigr)\eta=\bigl(
 {\cal L}_\rho(V)\circ {\cal L}_\rho(W)-{\cal L}_\rho(W)\circ
 {\cal L}_\rho(V)\bigr)\eta\,.\end{equation*}
\par\smallskip
{\rm 5.}  For each $V\in A^1(M,E)$, $f\in C^\infty(M,\RR)$ and
$\eta\in\Omega(M,E)$,
 \begin{equation*}{\cal L}_\rho(fV)\eta=f{\cal
 L}_\rho(V)\eta+d_\rho f\wedge i(V)\eta\,.\end{equation*}
\end{prop}

\begin{proof}
{\rm 1.} Let $W\in A^1(M,E)$. Then
 \begin{equation*}
 \begin{split}
 \bigl\langle{\cal L}_\rho(V)(d_\rho f),W\bigr\rangle
 &={\cal L}_\rho(V)\langle d_\rho f,W\rangle-\langle d_\rho
 f,\{V,W\}\rangle\\
 &={\cal L}(\rho\circ V)\circ{\cal L}(\rho\circ W)f
 -{\cal L}\bigl(\rho\circ\{V,W\}\bigr)f\\
 &={\cal L}(\rho\circ W)\circ{\cal L}(\rho\circ V)f\\
 &=\bigl\langle d_\rho\bigl({\cal
 L}_\rho(V)f\bigr),W\bigr\rangle\,,
 \end{split}
 \end{equation*}
so Property~1 is proven.
\par\medskip
{\rm 2.} Let $V$ and $W\in A^1(M,E)$, $\eta\in \Omega^p(M,E)$,
$V_1,\ldots,V_{p-1}\in A^1(M,E)$. We may write
 \begin{equation*}
 \begin{split}
 \bigl({\cal L}_\rho(V)\circ
 i(W)\eta\bigr)(V_1,
 &\ldots,V_{p-1})\\
 &={\cal L}_\rho(V)\bigl(\eta(W,V_1,\ldots,V_{p-1})\bigr)\\
 &\quad-\sum_{k=1}^{p-1}\eta\bigl(W,V_1,\ldots,V_{k-1},\{V,V_k\},
 V_{k+1},\ldots,V_{p-1}\bigr)\\
 &=\bigl({\cal L}_\rho(V)\eta\bigr)(W,V_1,\ldots,V_{p-1})
 +\eta\bigl(\{V,W\},V_1,\ldots, V_{p-1}\bigr)\\
 &=\Bigl(\bigl(i(W)\circ{\cal
 L}_\rho(V)+i(\{V,W\})\bigr)\eta\Bigr)
 (V_1,\ldots,V_{p-1})\,,
 \end{split}
 \end{equation*}
so Property~2 is proven.
\par\medskip
{\rm 3.} Let $V\in A^1(M,E)$, $\eta\in \Omega^p(M,E)$ and
$\zeta\in \Omega^q(M,E)$. For $p<0$, as well as for $q<0$, both
sides of the equality stated in Property~3 vanish, so that
equality is trivially satisfied. For $p=q=0$, that equality is
also satisfied, as shown by Equality $(*)$ in the proof of
Proposition~\ref{Lie derivative by section of algebroid 2}. We still have to
prove that equality for $p>0$
and (or) $q>0$. We will do that by induction on the value of
$p+q$. Let $r\geq 1$ be an integer such that the equality stated
in Property~3 holds for $p+q\leq r-1$. Such an integer exists, for
example $r=1$. We assume now that $p\geq 0$ and $q\geq 0$ are such
that $p+q=r$. Let $W\in A^1(M,E)$. By using Property~2, we may
write
 \begin{equation*}
 \begin{split}
 i(W)\circ{\cal L}_\rho(V)(\eta_\wedge\zeta)
 &={\cal L}_\rho(V)\circ
 i(W)(\eta\wedge\zeta)-i(\{V,W\})(\eta\wedge\zeta)\\
 &={\cal L}_\rho(V)\bigl(
 i(W)\eta\wedge\zeta+(-1)^p\eta\wedge i(W)\zeta\bigr)\\
 &\quad -i(\{V,W\})\eta\wedge\zeta-(-1)^p\eta\wedge
 i(\{V,W\})\zeta\,.
 \end{split}
 \end{equation*}
Since $i(W)\eta\in\Omega^{p-1}(M,E)$ and
$i(W)\zeta\in\Omega^{q-1}(M,E)$, the induction assumption allows
us to use Property~3 to transform the first terms of the right
hand side. We obtain
 \begin{equation*}
 \begin{split}
 i(W)\circ{\cal L}_\rho(V)(\eta_\wedge\zeta)
 &=\bigl({\cal L}_\rho(V)\circ
 i(W)\eta\bigr)\wedge\zeta+i(W)\eta\wedge{\cal L}_\rho(V)\zeta\\
 &\quad+(-1)^p\bigl({\cal L}_\rho(V)\eta\bigr)\wedge
 i(W)\zeta
 +(-1)^p\eta\wedge\bigl({\cal L}_\rho(V)\circ
 i(W)\zeta\bigr)\\
 &\quad-i(\{V,W\})\eta\wedge\zeta-(-1)^p\eta\wedge
 i(\{V,W\})\zeta\,.
 \end{split}
 \end{equation*}
By rearranging the terms, we obtain
 \begin{equation*}
 \begin{split}
 i(W)\circ{\cal L}_\rho(V)(\eta_\wedge\zeta)
 &=\bigl({\cal L}_\rho(V)\circ
 i(W)\eta-i(\{V,W\})\eta\bigr)\wedge\zeta\\
 &\quad+(-1)^p\eta\wedge\bigl({\cal L}_\rho(V)\circ i(W)\zeta
 -i(\{V,W\})\zeta\bigr)\\
 &\quad+i(W)\eta\wedge{\cal L}_\rho(V)\zeta
 +(-1)^p\bigl({\cal L}_\rho(V)\eta\bigr)\wedge i(W)\zeta\,.
 \end{split}
 \end{equation*}
By using again Property~2 we get
 \begin{equation*}
 \begin{split}
 i(W)\circ{\cal L}_\rho(V)(\eta_\wedge\zeta)
 &=\bigl(i(W)\circ{\cal L}_\rho(V)\eta\bigr)\wedge\zeta
   +(-1)^p\eta\wedge\bigl(i(W)\circ{\cal L}_\rho(V)\zeta\bigr)\\
 &\quad+i(W)\eta\wedge{\cal L}_\rho(V)\zeta
   +(-1)^p{\cal L}_\rho(V)\eta\wedge i(W)\zeta\\
 &=i(W)\bigl({\cal L}_\rho(V)\eta\wedge\zeta+\eta\wedge{\cal
 L}_\rho(V)\zeta\bigr)\,.
 \end{split}
 \end{equation*}
Since that last equality holds for all $W\in A^1(M,E)$, it follows
that Property~3 holds for $\eta\in\Omega^p(M,E)$ and
$\zeta\in\Omega^q(M,E)$, with $p\geq 0$, $q\geq 0$ and $p+q=r$. We
have thus proven by induction that Property~3 holds for all $p$
and $q\in\ZZ$, $\eta\in\Omega^p(M,E)$, $\zeta\in\Omega^q(M,E)$.
The same equality holds, by bilinearity, for all $\eta$ and
$\zeta\in\Omega(M,E)$.
\par\medskip
{\rm 4.} Let $V$ and $W\in A^1(M,E)$. Then $\{V,W\}\in A^1(M,E)$
and, by Property~3, ${\cal L}_\rho(V)$, ${\cal L}_\rho(W)$ and
${\cal L}_\rho(\{V,W\})$ are derivations of degree $0$ of the
graded algebra $\Omega(M,E)$. By \ref{properties of derivations}~(ii), the
graded bracket
 \begin{equation*}\bigl[{\cal L}_\rho(V),{\cal L}_\rho(W)\bigr]
 ={\cal L}_\rho(V)\circ {\cal L}_\rho(W)-{\cal L}_\rho(W)\circ
 {\cal L}_\rho(V)\end{equation*}
is also a derivation of degree $0$ of $\Omega(M,E)$. Property~4
means that the derivations ${\cal L}_\rho(\{V,W\})$ and
$\bigl[{\cal L}_\rho(V),{\cal L}_\rho(W)\bigr]$ are equal. In
order to prove that equality, it is enough to prove that it holds
true for $\eta\in \Omega^0(M,E)$ and for $\eta\in\Omega^1(M,E)$,
since the graded algebra $\Omega(M,E)$ is generated by its
homogeneous elements of degrees $0$ and $1$.
\par\smallskip
Let $f\in\Omega^0(M,E)=C^\infty(M,\RR)$. We have
 \begin{equation*}
 \begin{split}
 {\cal L}_\rho(\{V,W\})f
 &={\cal L}\bigl(\rho\circ\{V,W\}\bigr)f\\
 &={\cal L}\bigl([\rho\circ V,\rho\circ W]\bigr)f\\
 &=\bigl[{\cal L}(\rho\circ V),{\cal L}(\rho\circ W)\bigr]f\\
 &=\bigl[{\cal L}_\rho(V),{\cal L}_\rho(W)\bigr]f\,,
 \end{split}
 \end{equation*}
therefore Property~4 holds for $\eta=f\in\Omega^0(M,E)$.
\par\smallskip
Now let $\eta\in\Omega^1(M,E)$ and $Z\in A^1(M,E)$. By using
Property~2, then Property~4 for elements $\eta\in\Omega^0(M,E)$,
we may write
\goodbreak
 \begin{equation*}
 \begin{split}
 i(Z)\circ{\cal L}_\rho(\{V,W\})\eta
 &={\cal
 L}_\rho(\{V,W\})\bigl(i(Z)\eta\bigr)-i\bigl(\bigl\{\{V,W\},Z\bigr\}
 \bigr)\eta\\
 &=\bigl({\cal L}_\rho(V)\circ {\cal L}_\rho(W)
   -{\cal L}_\rho(W)\circ {\cal
   L}_\rho(V)\bigr)\bigl(i(Z)\eta\bigr)\\
  &\quad-i\bigl(\bigl\{\{V,W\},Z\bigr\}\bigr)\eta\,.
 \end{split}
  \end{equation*}
By using Property~2 and the Jacobi identity, we obtain
 \begin{equation*}
 \begin{split}
 i(Z)\circ{\cal L}_\rho(\{V,W\})\eta
 &={\cal L}_\rho(V)\bigl(i(\{W,Z\})\eta+i(Z)\circ{\cal
 L}_\rho(W)\eta\bigr)\\
 &\quad-{\cal L}_\rho(W)\bigl(i(\{V,Z\})\eta+i(Z)\circ{\cal
 L}_\rho(V)\eta\bigr)\\
 &\quad-i\bigl(\bigl\{\{V,W\},Z\bigr\}\bigr)\eta\\
 &=i(\{W,Z\}){\cal L}_\rho(V)\eta + i(\{V,Z\}){\cal
 L}_\rho(W)\eta\\
 &\quad-i(\{V,Z\}){\cal L}_\rho(W)\eta
 - i(\{W,Z\}){\cal L}_\rho(V)\eta\\
 &\quad+i(Z)\circ\bigl({\cal L}_\rho(V)\circ{\cal L}_\rho(W)
 -{\cal L}_\rho(W)\circ{\cal L}_\rho(V)\bigr)\eta\\
 &\quad+i\Bigl(\bigl\{V,\{W,Z\}\bigr\}-\bigl\{W,\{V,Z\}\bigr\}
 -\bigl\{\{V,W\},Z\bigr\}\Bigr)\eta\\
 &=i(Z)\circ\bigl({\cal L}_\rho(V)\circ{\cal L}_\rho(W)
 -{\cal L}_\rho(W)\circ{\cal L}_\rho(V)\bigr)\eta\,.
 \end{split}
 \end{equation*}
Since that last equality holds for all $Z\in A^1(M,E)$, Property~4
holds for all $\eta\in\Omega^1(M,E)$, and therefore for all
$\eta\in\Omega(M,E)$.
\par\medskip
{\rm 5.} Let $V\in A^1(M,E)$ and $f\in C^\infty(M,\RR)$. We have
seen (Property~4) that ${\cal L}_\rho(fV)$ is a derivation of
degree $0$ of $\Omega(M,E)$. We easily verify that
 \begin{equation*}\eta\mapsto f{\cal L}_\rho(V)\eta+d_\rho f\wedge
i(V)\eta\end{equation*} is too a derivation of degree $0$ of
$\Omega(M,E)$. Property~5 means that these two derivations are
equal. As above, it is enough to prove that Property~5 holds for
$\eta\in\Omega^0(M,E)$ and for $\eta\in\Omega^1(M,E)$.
\par\smallskip
Let $g\in\Omega^0(M,E)=C^\infty(M,\RR)$. We may write
 \begin{equation*}{\cal L}_\rho(fV)g=i(fV)d_\rho g=f{\cal
L}_\rho(V)g\,,\end{equation*} which shows that Property~5 holds
for $\eta=g\in\Omega^0(M,E)$.
\par\smallskip
Let $\eta\in\Omega^1(M,E)$, and $W\in A^1(M,E)$. We have
 \begin{equation*}
 \begin{split}
 \bigl\langle {\cal L}_\rho(fV)\eta,W\bigr\rangle
 &={\cal L}_\rho(fV)\bigl(\langle\eta,W\rangle\bigr)
  -\bigl\langle\eta,\{fV,W\}\bigr\rangle\\
 &=f{\cal L}_\rho(V)\bigl(\langle\eta,W\rangle\bigr)
  -f\bigl\langle\eta,\{V,W\}\bigr\rangle
  +\bigl\langle\eta,\bigl(i(W)d_\rho
  f\bigr)V\bigr\rangle\\
 &=\bigl\langle f{\cal L}_\rho(V)\eta,W\bigr\rangle
  +\bigl(i(W)d_\rho f\bigr)i(V)\eta\\
 &=\bigl\langle f{\cal L}_\rho(V)\eta +d_\rho f\wedge i(V)\eta,W\bigr\rangle
 \,,
 \end{split}
 \end{equation*}
since, $\eta$ being in $\Omega^1(M,E)$, $i(W)\circ i(V)\eta=0$.
The last equality being satisfied for all $W\in A^1(M,E)$, the
result follows.
\end{proof}
\par\bigskip
The next Proposition shows that for each $V\in A^1(M,E)$, the Lie
derivative ${\cal L}_\rho(V)$, already defined as a derivation of
degree $0$ of the graded algebra $\Omega(M,E)$, can also be
extended into a derivation of degree $0$ of the graded algebra
$A(M,E)$, with very nice properties. As we will soon see, the
Schouten-Nijenhuis bracket will appear as a very natural further
extension of the Lie derivative.

\begin{prop}
\label{Lie derivative of multivectors} 
Let $(E,\tau,M,\rho)$ be a Lie algebroid on a
smooth manifold $M$. For each smooth section $V\in A^1(M,E)$ of
the vector bundle $(E,\tau,M)$, there exists a unique graded
endomorphism of degree $0$ of the graded algebra of multivectors
$A(M,E)$, called the {\it Lie derivative\/} with respect to $V$
and denoted by ${\cal L}_\rho(V)$, which satisfies the following
properties:

\begin{description}
\item{\rm(i)} For a smooth function
$f\in A^0(M,E)=C^\infty(M,\RR)$,
 \begin{equation*}{\cal L}_\rho(V)f=i(\rho\circ V)df={\cal L}(\rho\circ
V)f\,,\end{equation*} where ${\cal L}(\rho\circ V)$ denotes the
usual Lie derivative with respect to the vector field $\rho\circ
V$;

\item{\rm (ii)} For an integer $p\geq 1$ and a multivector
$P\in A^p(M,E)$, ${\cal L}_\rho(V)P$ is the unique element in
$A^p(M,E)$ such that, for all $\eta\in \Omega^p(M,E)$,
 \begin{equation*}\bigl\langle\eta,{\cal L}_\rho(V)P\bigr\rangle
 ={\cal L}_\rho(V)\bigl(\langle\eta,P\rangle\bigr)-\bigl\langle
 {\cal L}_\rho(V)\eta,P\bigr\rangle\,.\end{equation*}
\end{description}
\end{prop}

\begin{proof} Let us first observe that
$A^0(M,E)=\Omega^0(M,E)=C^\infty(M,\RR)$, and that for $f\in
A^0(M,E)$, the definition of ${\cal L}_\rho(V)f$ given above is
the same as that given in Proposition 
\ref{Lie derivative by section of algebroid 2}.
\par\smallskip
Now let $p\geq 1$ and $P\in A^p(M,E)$. The map
 \begin{equation*}\eta\mapsto K(\eta)=
 {\cal L}_\rho(V)\bigl(\langle\eta,P\rangle\bigr)-\bigl\langle
 {\cal L}_\rho(V)\eta,P\bigr\rangle\end{equation*}
is clearly an $\RR$-linear map defined on $\Omega^p(M,E)$, with
values in $C^\infty(M,\RR)$. Let $f\in C^\infty(M,\RR)$. We have
 \begin{equation*}
 \begin{split}
 K(f\eta)
 &={\cal L}_\rho(V)\bigl(\langle f\eta,P\rangle\bigr)-\bigl\langle
 {\cal L}_\rho(V)(f\eta),P\bigr\rangle\\
 &=f\Bigl({\cal L}_\rho(V)\bigl(\langle\eta,P\rangle\bigr)-\bigl\langle
 {\cal L}_\rho(V)\eta,P\bigr\rangle\Bigr)\\
 &\quad+\bigl({\cal L}_\rho(V)f\bigr)\langle\eta,P\rangle
        -\bigl({\cal L}_\rho(V)f\bigr)\langle\eta,P\rangle\\
 &=fK(\eta)\,.
 \end{split}
 \end{equation*}
This proves that the map $K$ is $C^\infty(M,\RR)$-linear. Since
the pairing allows us to consider the vector bundle
$(\bigwedge^pE,\tau,M)$ as the dual of $(\bigwedge^pE^*,\pi,M)$,
we see that there exists a unique element ${\cal L}_\rho(V)P\in
A^p(M,E)$ such that, for all $\eta\in \Omega^p(M,E)$,
 \begin{equation*}K(\eta)=
 {\cal L}_\rho(V)\bigl(\langle\eta,P\rangle\bigr)-\bigl\langle
 {\cal L}_\rho(V)\eta,P\bigr\rangle=\bigl\langle\eta,{\cal
 L}_\rho(V)P\bigr\rangle\,,\end{equation*}
and that ends the proof.
\end{proof}

\begin{prop}
\label{Lie derivative of multivectors 2}
Under the assumptions of 
Proposition~\ref{Lie derivative of multivectors}, 
the Lie derivative has the following properties.
\par\smallskip
{\rm 1.} For each $V$ and $W\in A^1(M,E)$,
 \begin{equation*}{\cal L}_\rho(V)(W)=\{V,W\}\,.\end{equation*}
\par\smallskip
{\rm 2.} For $V,V_1,\ldots,V_p\in A^1(M,E)$,
 \begin{equation*}{\cal L}_\rho(V)(V_1\wedge\cdots\wedge V_p)
 =\sum_{i=1}^pV_1\wedge\cdots\wedge V_{i-1}\wedge\{V,V_i\}\wedge V_{i+1}
 \wedge\cdots\wedge V_p\,.\end{equation*}
\par\smallskip
{\rm 3.} For each $V\in A^1(M,E)$, ${\cal L}_\rho(V)$ is a
derivation of degree $0$ of the exterior algebra $A(M,E)$. That
means that for all $P$ and $Q\in A(M,E)$,
 \begin{equation*}{\cal L}_\rho(V)(P\wedge Q)=
   \bigl({\cal L}_\rho(V)P\bigr)\wedge Q+
   P\wedge{\cal L}_\rho(V)Q\,.\end{equation*}
\par\smallskip
{\rm 4.} For each $V\in A^1(M,E)$, $P\in A(M,E)$ and
$\eta\in\Omega(M,E)$,
 \begin{equation*}{\cal L}_\rho(V)\bigl(i(P)\eta\bigr)=i\bigl({\cal
 L}_\rho(V)P\bigr)\eta
 +i(P)\bigl({\cal L}_\rho(V)\eta\bigr)\,.\end{equation*}
\par\smallskip
{\rm 5.} Similarly, for each $V\in A^1(M,E)$, $P\in A(M,E)$ and
$\eta\in\Omega(M,E)$,
 \begin{equation*}{\cal L}_\rho(V)\bigl(\langle\eta,P\rangle\bigr)
 =\bigl\langle{\cal L}_\rho(V)\eta,P\bigr\rangle
 +\bigl\langle\eta,{\cal L}_\rho(V)P\bigr\rangle\,.\end{equation*}
\par\smallskip
{\rm 6.}  For each $V$ and $W\in A^1(M,E)$, $P\in A(M,E)$,
 \begin{equation*}{\cal L}_\rho\bigl(\{V,W\}\bigr)P=\bigl(
 {\cal L}_\rho(V)\circ {\cal L}_\rho(W)-{\cal L}_\rho(W)\circ
 {\cal L}_\rho(V)\bigr)P\,.\end{equation*}
\par\smallskip
{\rm 7.}  For each $V\in A^1(M,E)$, $f\in C^\infty(M,\RR)$, $P\in
A(M,E)$ and $\eta\in \Omega(M,E)$,
 \begin{equation*}\bigl\langle\eta,{\cal L}_\rho(fV)P\bigr\rangle
 =f\bigl\langle{\cal L}_\rho(V)P,\eta\bigr\rangle
 +\bigl\langle d_\rho f\wedge i(V)\eta,P\bigr\rangle\,.\end{equation*}
\end{prop}

\begin{proof}
{\rm 1.} Let $V$ and $W\in A^1(M,E)$, $\eta\in\Omega(M,E)$. We may
write
 \begin{equation*}
 \begin{split}
 \bigl\langle \eta,{\cal L}_\rho(V)W\bigr\rangle
 &={\cal L}_\rho(V)\bigl(\langle\eta,W\rangle\bigr)-\bigl\langle{\cal
 L}_\rho(V)\eta,W\bigr\rangle\\
 &=\bigl\langle\eta,\{V,W\}\bigr\rangle\,.
 \end{split}
 \end{equation*}
We have proven Property~1.
\par\smallskip
{\rm 2.} The proof is similar to that of Property~1.
\par\smallskip
{\rm 3.} When $P=V_1\wedge\cdots\wedge V_p$ and $Q=W_1\wedge\cdots
\wedge W_q$ are decomposable homogeneous elements in $A(M,E)$,
Property~3 is an easy consequence of 2. The validity of Property~3
for all $P$ and $Q\in A(M,E)$ follows by linearity.
\par\smallskip
{\rm 4.}  When $P=V_1\wedge\cdots\wedge V_p$ is a decomposable
homogeneous element in $A(M,E)$, Property~4 is an easy consequence
of Property~2. The validity of Property~4 for all $P$ and $Q\in
A(M,E)$ follows by linearity.
\par\smallskip
{\rm 5.} This is an immediate consequence of Property~4.
\par\smallskip
{\rm 6.} This is an immediate consequence of Property~4 of this
Proposition and of Property~4 of Proposition~\ref{Lie derivative by section of
algebroid 3}.
\par\smallskip
{\rm 7.} This is an immediate consequence of Property~4 of this
Proposition and of Property~5 of 
Proposition~\ref{Lie derivative by section of algebroid 3}.
\end{proof}

\section{The $\Omega(M,E)$-valued exterior derivative}
\label{exterior derivative}
We have introduced above 
(Definition~\ref{exterior derivative of functions}) 
the $\Omega(M,E)$-valued exterior derivative of a function $f\in
\Omega^0(M,E)=C^\infty(M,\RR)$. The next proposition shows that
the $\Omega(M,E)$-valued exterior derivative extends as a graded
endomorphism of degree $1$ of the graded algebra $\Omega(M,E)$. We
will see later (Proposition~\ref{exterior derivative of forms 2})
that the $\Omega(M,E)$-valued
exterior derivative is in fact a derivation of degree $1$ of
$\Omega(M,E)$.

\begin{prop}
\label{exterior derivative of forms 1} 
Let $(E,\tau,M,\rho)$ be a Lie algebroid
over a
smooth manifold $M$. There exists a unique graded endomorphism of
degree $1$ of the exterior algebra of forms $\Omega(M,E)$, called
the {\it $\Omega(M,E)$-valued exterior derivative\/} (or, in
brief, the {\it exterior derivative\/}) and denoted by $d_\rho$,
which satisfies the following properties:
\par\smallskip
\item{\rm(i)} For $f\in\Omega^0(M,E)=C^\infty(M,\RR)$, $d_\rho f$ is the
unique element in $\Omega^1(M,E)$, already defined 
(Definition~\ref{exterior derivative of functions}), 
such that, for each $V\in A^1(M,E)$,
 \begin{equation*}\langle d_\rho f,V\rangle={\cal L}_\rho(V)f=\langle
df,\rho\circ
 V\rangle=\langle{}^t\!\rho\circ df,V\rangle\,,\end{equation*}
where $d$ stands for the usual exterior derivative of smooth
functions on $M$, and ${}^t\!\rho:E^*\to T^*M$ is the transpose of
the anchor $\rho$.
\par\smallskip
\item{\rm(ii)} For $p\geq 1$ and $\eta\in \Omega^p(M,E)$, $d_\rho\eta$ is
the unique element in $\Omega^{p+1}(M,E)$ such that, for all
$V_0,\ldots, V_p\in A^1(M,E)$,
 \begin{equation*}
 \begin{split}
 d_\rho\eta(V_0,\ldots,V_p)
 &=\sum_{i=0}^p(-1)^i{\cal L}_\rho(V_i)\bigl(\eta(V_0,\ldots,\widehat{V}_i,
 \ldots,V_p)\bigr)\\
 &\quad+\sum_{0\leq i<j\leq p}(-1)^{i+j}\eta\bigl(\{V_i,V_j\},V_0,\ldots,
 \widehat{V}_i,\ldots,\widehat{V}_j,\ldots,V_p\bigr)\,,
 \end{split}
 \end{equation*}
where the symbol~~$\widehat{\ }$~~over the terms $V_i$ and $V_j$
means that these terms are omitted.

\end{prop}

\begin{proof} For $f\in\Omega^0(M,E)$, $d_\rho f$ is clearly an element
in $\Omega^1(M,E)$.
\par
Let $p\geq 1$ and $\eta\in\Omega^p(M,E)$. As defined in {\rm(ii)},
$d_\rho\eta$ is a map, defined on $\bigl(A^1(M,E)\bigr)^p$, with
values in $C^\infty(M,\RR)$. The reader will immediately see that
this map is skew-symmetric and $\RR$-linear in each of its
arguments. In order to prove that $d_\rho\eta$ is an element in
$\Omega^{p+1}(M,E)$, it remains only to verify that as a map,
$d_\rho\eta$ is $C^\infty(M,\RR)$-linear in each of its arguments,
or simply in its first argument, since the skew-symmetry will
imply the same property for all other arguments. Let $f\in
C^\infty(M,\RR)$. We have
 \begin{equation*}
 \begin{split}
 d_\rho\eta(fV_0,&V_1,\ldots,V_p)\\
 &={\cal L}_\rho(fV_0)\bigl(\eta(V_1,\ldots,V_p)\bigr)\\
 &\quad+\sum_{i=1}^p(-1)^i{\cal
 L}_\rho(V_i)\bigl(f\eta(V_0,\ldots,\widehat{V}_i,\ldots,V_p)\bigr)\\
 &\quad+\sum_{1\leq j\leq
 p}(-1)^j\eta\bigl(\{fV_0,V_j\},V_1,\ldots,\widehat{V}_j,\ldots,V_p\bigr)\\
 &\quad+\sum_{1\leq i<j\leq p}(-1)^{i+j}\eta\bigl(\{V_i,V_j\},fV_0,V_1,
 \ldots,\widehat{V}_i,\ldots,\widehat{V}_j,\ldots,V_p\bigr)\,.
 \end{split}
 \end{equation*}
By a rearrangement of the terms in the right hand side, and by
using the formulae
 \begin{equation*}{\cal L}_\rho(V_i)\bigl(f\eta(\ldots)\bigr)=\bigl({\cal
 L}_\rho(V_i)f\bigr)\eta(\ldots) + f{\cal
 L}_\rho(V_i)\bigl(\eta(\ldots)\bigr)\,,\end{equation*}
and
 \begin{equation*}\{fV_0,V_j\}=f\{V_0,V_j\}-\bigl({\cal
L}_\rho(V_j)f\bigr)V_0\,,\end{equation*} we obtain
 \begin{equation*}d_\rho\eta(fV_0,V_1,\ldots,V_p)=fd_\rho\eta(V_0,V_1,\ldots,V_p
)\,.\end{equation*} We have shown that
$d_\rho\eta\in\Omega^{p+1}(M,E)$. The $\Omega(M,E)$-valued
exterior derivative so defined on $\Omega^p(M,E)$ for all
$p\in\ZZ$ extends, in a unique way, into a graded endomorphism of
degree $1$ of $\Omega(M,E)$. \end{proof}

\begin{rmk} Let $p\geq 1$, $\eta\in\Omega^p(M,E)$ and
$V_0,\ldots,V_p\in A^1(M,E)$. The formula for $d_\rho\eta$ given
in Proposition~\ref{exterior derivative of forms 1} can be cast into
another form, often useful:
 \begin{equation*}
 \begin{split}
 d_\rho\eta(V_0,\ldots,V_p)
 &=\sum_{i=0}^p(-1)^i\bigl({\cal
 L}_\rho(V_i)\eta\bigr)(V_0,\ldots,\widehat{V_i},
 \ldots,V_p)\\
 &\quad-\sum_{0\leq i<j\leq p}(-1)^{i+j}\eta\bigl(\{V_i,V_j\},V_0,\ldots,
 \widehat{V}_i,\ldots,\widehat{V}_j,\ldots,V_p\bigr)\,.
 \end{split}
 \end{equation*}
\par\smallskip
For example, for $p=1$,
 \begin{equation*}
 \begin{split}
 d_\rho\eta(V_0,V_1)
 &={\cal L}_\rho(V_0)\bigl(\eta(V_1)\bigr)
   -{\cal L}_\rho(V_1)\bigl(\eta(V_0)\bigr)
   -\eta\bigl(\{V_0,V_1\}\bigr)\\
 &=\bigl\langle{\cal L}_\rho(V_0)\eta,V_1\bigr\rangle
   -\bigl\langle{\cal L}_\rho(V_1)\eta,V_0\bigr\rangle
   +\eta\bigl(\{V_0,V_1\}\bigr)\,.
 \end{split}
   \end{equation*}
\end{rmk}

\begin{prop}
\label{exterior derivative of forms 2}
Under the assumptions of 
Proposition~\ref{exterior derivative of forms 1},
the $\Omega(M,E)$-valued exterior derivative has the following
properties:
\par\smallskip
{\rm 1.} Let $V\in A^1(M,E)$. The Lie derivative ${\cal
L}_\rho(V)$, the exterior derivative $d_\rho$ and the interior
product $i(V)$ are related by the formula
 \begin{equation*}{\cal L}_\rho(V)=i(V)\circ d_\rho+d_\rho\circ
i(V)\,.\end{equation*}
\par\smallskip
{\rm 2.} The exterior derivative $d_\rho$ is a derivation of
degree $1$ of the exterior algebra $\Omega(M,E)$. That means that
for each $p\in\ZZ$, $\eta\in\Omega^p(M,E)$ and
$\zeta\in\Omega(M,E)$,
 \begin{equation*}d_\rho(\eta\wedge\zeta)=d_\rho\eta\wedge\zeta+(-1)^p\eta\wedge
 d_\rho\zeta\,.\end{equation*}
\par\smallskip
{\rm 3.} Let $V\in A^1(M,E)$. Then
 \begin{equation*}{\cal L}_\rho(V)\circ d_\rho=d_\rho\circ{\cal
L}_\rho(V)\,.\end{equation*}
\par\smallskip
{\rm 4.} The square of $d_\rho$ vanishes identically,
 \begin{equation*}d_\rho\circ d_\rho=0\,.\end{equation*}

\end{prop}

\begin{proof}
{\rm 1.} Let $V_0=V,\ V_1,\ \ldots,\ V_p\in A^1(M,E)$,
$\eta\in\Omega^p(M,E)$. Then
 \begin{equation*}
 \begin{split}
 \bigl(i(V)\circ d_\rho \eta\bigr)(V_1,
 &\ldots,V_p)\\
 &=d_\rho\eta(V,V_1,\ldots,V_p)\\
 &=\sum_{i=0}^p(-1)^i{\cal
 L}_\rho(V_i)\bigl(\eta(V_0,\ldots,\widehat{V}_i,\ldots,V_p)\bigr)\\
 &\quad+\sum_{0\leq i<j\leq p}(-1)^{i+j}\eta\bigl(\{V_i,V_j\},V_0,\ldots,
 \widehat{V}_i,\ldots,\widehat{V}_j,\ldots,V_p\bigr)\,,
 \end{split}
 \end{equation*}
and
 \begin{equation*}
 \begin{split}
 \bigl(d_\rho\circ i(V) \eta\bigr)(V_1,
 &\ldots,V_p)\\
 &=\sum_{i=1}^p(-1)^{i-1}{\cal
 L}_\rho(V_i)\bigl(\eta(V_0,\ldots,\widehat{V_i},\ldots,V_p)\bigr)\\
 &\quad+\sum_{1\leq i<j\leq
 p}(-1)^{i+j}\eta\bigl(V_0,\{V_i,V_j\},V_1,\ldots,
 \widehat{V}_i,\ldots,\widehat{V}_j,\ldots,V_p\bigr)\,.
 \end{split}
 \end{equation*}
Let us add these two equalities. Several terms cancel, and we
obtain
 \begin{equation*}
 \begin{split}
 \Bigl(\bigl(i(V)\circ d_\rho
 &+d_\rho\circ i(V)\bigr)\eta\Bigr)(V_1,\ldots,V_p)\\
 &={\cal L}_\rho(V_0)\bigl(\eta(V_1,\ldots,V_p)\bigr)
 +\sum_{j=1}^p(-1)^j\eta\bigl(\{V_0,V_j\},V_1,
 \ldots,\widehat{V}_j,\ldots,
 V_p\bigr)\,.
 \end{split}
 \end{equation*}
When we shift, in the last term of the right hand side, the
argument $\{V_0,V_j\}$ to the $j$-th position, we obtain
 \begin{equation*}
 \begin{split}
 \Bigl(\bigl(i(V)&\circ d_\rho+d_\rho\circ i(V)\bigr)
 \eta\Bigr)(V_1,\ldots,V_p)\\
 &={\cal L}_\rho(V_0)\bigl(\eta(V_1,\ldots,V_p)\bigr)
 +\sum_{j=1}^p\eta\bigl(V_1,\ldots,V_{j-1},\{V_0,V_j\},V_{j+1},\ldots,
 V_p\bigr)\\
 &=\bigl({\cal L}_\rho(V_0)\eta\bigr)(V_1,\ldots,V_p)\,.
 \end{split}
 \end{equation*}
\par\smallskip
{\rm 2.} For $\eta=f$ and
$\zeta=g\in\Omega^0(M,E)=C^\infty(M,\RR)$, Property 2 holds since
we have, for alll $V\in A^1(M,E)$,
 \begin{equation*}
 \begin{split}
 \bigl\langle d_\rho(fg),V\bigr\rangle
 =\bigl\langle d(fg),\rho\circ V\bigr\rangle
 &=\langle f\,dg+g\,df,\rho\circ V\rangle\\
 &=\langle f\,d_\rho g+g\,d_\rho f,V\rangle\,.
 \end{split}
 \end{equation*}
Now let $p\geq 0$ and $q\geq 0$ be two integers, $\eta\in
\Omega^p(M,E)$, $\zeta\in\Omega^q(M,E)$. We will prove that
Property 2 holds by induction on $p+q$. Just above, we have seen
that it holds for $p+q=0$. Let us assume that $r$ is an integer
such that Property 2 holds for $p+q\leq r$, and that now
$p+q=r+1$. Let $V\in A^1(M,E)$. We may write
 \begin{equation*}
 \begin{split}
 i(V)d_\rho(\eta\wedge\zeta)
 &={\cal L}_\rho(V)(\eta\wedge\zeta)-d_\rho\circ
 i(V)(\eta\wedge\zeta)\\
 &={\cal L}_\rho(V)\eta\wedge\zeta+\eta\wedge{\cal
 L}_\rho(V)\zeta\\
 &\quad-d_\rho\bigl(i(V)\eta\wedge\zeta+(-1)^p\eta\wedge
 i(V)\zeta\bigr)\,.
 \end{split}
 \end{equation*}
We may now use the induction assumption, since in the last terms
of the right hand side $i(V)\eta\in\Omega^{p-1}(M,E)$ and
$i(V)\zeta\in\Omega^{q-1}(M,E)$. After some rearrangements of the
terms we obtain
 \begin{equation*}i(V)d_\rho(\eta\wedge\zeta)
 =i(V)(d_\rho\eta\wedge\zeta+\eta\wedge
 d_\rho\zeta)\,.
 \end{equation*}
Since that result holds for all $V\in A^1(M,E)$, Property 2 holds
for $p+q=r+1$, and therefore for all $p$ and $q\in\ZZ$.
\par\smallskip
{\rm 3.} Let $V\in A^1(M,E)$. We know 
(Proposition~\ref{Lie derivative by section of algebroid 3}) 
that ${\cal L}_\rho(V)$ is a derivation of degree $0$ of the exterior
algebra $\Omega(M,E)$, and we have just seen (Property 2) that
$d_\rho$ is a derivation of degree $1$ of that algebra. Therefore,
by~\ref{properties of derivations}, their graded bracket
 \begin{equation*}\bigl[{\cal L}_\rho(V),d_\rho\bigr]={\cal L}_\rho(V)\circ
d_\rho
 -d_\rho\circ{\cal L}_\rho(V)\end{equation*}
is a derivation of degree $1$ of $\Omega(M,E)$. In order to prove
that that derivation is equal to $0$, it is enough to prove that
it vanishes on $\Omega^0(M,E)$ and on $\Omega^1(M,E)$. We have
already proven that it vanishes on $\Omega^0(M,E)$ (Property~1 of
\ref{Lie derivative by section of algebroid 3}). Let $\eta\in\Omega^1(M,E)$ and
$W\in A^1(M,E)$. By using
Property~1 of this Proposition and Property~2 
of~\ref{Lie derivative by section of algebroid 3}, we may write
 \begin{equation*}
 \begin{split}
 i(W)\circ\bigl({\cal L}_\rho(V)\circ d_\rho
 &-d_\rho\circ {\cal
 L}_\rho(V)\bigr)\eta\\
 &={\cal L}_\rho(V)\circ i(W)\circ d_\rho\eta-i\bigl(\{V,W\}\bigr)\circ
 d_\rho\eta\\
 &\quad-{\cal L}_\rho(W)\circ {\cal L}_\rho(V)\eta
 +d_\rho\circ i(W)\circ {\cal L}_\rho(V)\eta\,.
 \end{split}
 \end{equation*}
By rearrangement of the terms, we obtain
 \begin{equation*}
 \begin{split}
 i(W)\circ\bigl({\cal L}_\rho(V)\circ d_\rho
 &-d_\rho\circ {\cal
 L}_\rho(V)\bigr)\eta\\
 &=\bigl({\cal L}_\rho(V)\circ{\cal L}_\rho(W)-{\cal L}_\rho(W)\circ
 {\cal L}_\rho(V)-{\cal L}_\rho(\{V,W\})\bigr)\eta\\
 &\quad+d_\rho\circ i(\{V,W\})\eta-d_\rho\circ i(\{V,W\})\eta\\
 &\quad-\bigl({\cal L}_\rho(V)\circ d_\rho
 -d_\rho\circ{\cal L}_\rho(V)\bigr)\bigl(i(W)\eta\bigr)\\
 &=0\,,
 \end{split}
 \end{equation*}
since $i(W)\eta\in\Omega^0(M,E)$, which implies that the last term
vanishes.
\par\smallskip
{\rm 4.} Property 2 shows that $d_\rho$ is a derivation of degree
$1$ of $\Omega(M,E)$. We know 
(\ref{properties of derivations}) 
that $[d_\rho,d_\rho]=2d_\rho\circ d_\rho$ is a derivation of degree
$2$ of $\Omega(M,\RR)$. In order to prove that $d_\rho\circ
d_\rho=0$, it is enough to prove that it vanishes on
$\Omega^0(M,E)$ and on $\Omega^1(M,E)$.
\par
Let $f\in\Omega^0(M,E)=C^\infty(M,\RR)$, $V$ and $W\in A^1(M,E)$.
Then
 \begin{equation*}
 \begin{split}
 (d_\rho\circ d_\rho f)(V,W)
 &={\cal L}_\rho(V)\bigl(d_\rho f(W)\bigr)
 -{\cal L}_\rho(W)\bigl(d_\rho f(V)\bigr)
 -d_\rho f\bigl(\{V,W\}\bigr)\\
 &=\Bigl({\cal L}_\rho(V)\circ{\cal L}_\rho(W)
   -{\cal L}_\rho(W)\circ{\cal L}_\rho(V)
   -{\cal L}_\rho\bigl(\{V,W\}\bigr)\Bigr)f\\
 &=0\,,
 \end{split}
 \end{equation*}
where we have used Property~4 of 
Proposition~\ref{Lie derivative by section of algebroid 3}. 
We have shown
that $d_\rho\circ d_\rho$ vanishes on $\Omega^0(M,E)$.
\par
Now let $\eta\in \Omega^1(M,E)$, $V_0$, $V_1$ and $V_2\in
A^1(M,E)$. Using Property 1, we may write
 \begin{equation*}
 \begin{split}
 (d_\rho\circ d_\rho \eta)(V_0,V_1,V_2)
 &=\Bigl(\bigl(i(V_0)\circ
 d_\rho\bigr)(d_\rho\eta)\Bigr)(V_1,V_2)\\
 &=\Bigl(\bigl({\cal L}_\rho(V_0)\circ d_\rho-d_\rho\circ
 i(V_0)\circ d_\rho\bigr)\eta\Bigr)(V_1,V_2)\,.
 \end{split}
 \end{equation*}
The last term in the right hand side may be transformed, by using
again Property 1:
 \begin{equation*}
 \begin{split}
 d_\rho\circ i(V_0)\circ d_\rho (\eta)
 &=d_\rho\circ {\cal L}_\rho(V_0)
 \eta-d_\rho\circ d_\rho\bigl(i(V_0)\eta\bigr)\\
 &=d_\rho\circ {\cal L}_\rho(V_0)\eta\,,
 \end{split}
 \end{equation*}
since, as $i(V_0)\eta\in \Omega^0(M,E)$, we have $d_\rho\circ
d_\rho\bigl(i(V_0)\eta\bigr)=0$. So we obtain
 \begin{equation*}(d_\rho\circ d_\rho \eta)(V_0,V_1,V_2)=
 \Bigl(\bigl({\cal L}_\rho(V_0)\circ d_\rho-d_\rho\circ
 {\cal L}_\rho(V_0)\bigr)\eta\Bigr)(V_1,V_2)\,.\end{equation*}
But Property 3 shows that
 \begin{equation*}\bigl({\cal L}_\rho(V_0)\circ d_\rho-d_\rho\circ
 {\cal L}_\rho(V_0)\bigr)\eta=0\,,\end{equation*}
so we have
 \begin{equation*}(d_\rho\circ d_\rho \eta)(V_0,V_1,V_2)=0\,,\end{equation*}
and our proof is complete.
\end{proof}

\section{Defining a Lie algebroid by properties of its dual}
Let $(E,\tau,M)$ be a vector bundle and $(E^*,\pi,M)$ its dual
bundle. We have seen in \ref{exterior derivative} that 
when $(E,\tau,M)$ has a Lie
algebroid structure whose anchor is denoted by $\rho$, this
structure determines, on the graded algebra $\Omega(M,E)$ of
sections of the exterior powers of the dual bundle
$(E^*,\pi,M)$, a graded derivation $d_\rho$, of degree $1$,
which satisfies $d_\rho^2=d_\rho\circ d_\rho=0$. Now we are going
to prove a converse of this property: when a graded derivation of
degree $1$, whose square vanishes, is given on $\Omega(M,E)$, it
determines a Lie algebroid structure on $(E,\tau,M)$. This property
will be used later to prove that the cotangent bundle of a Poisson
manifold has a natural Lie algebroid structure.
\par\smallskip
We will need the following lemmas.

\begin{lemma}
\label{deltabracket}
Let $(E,\tau,M)$ be a vector bundle and $(E^*,\pi,M)$ its dual
bundle. Let $\delta$ be a graded derivation of degree $1$ of the
exterior algebra $\Omega(M,E)$ (notations defined
in~\ref{The exterior algebra of sections}). 
For
each pair $(X,Y)$ of smooth sections of $(E,\tau,M)$ there exists a
unique smooth section $[X,Y]_\delta$ of $(E,\tau,M)$, called the
{\it $\delta$-bracket\/} of $X$ and $Y$, such that
 \begin{equation*}i\bigl([X,Y]_\delta\bigr)
 =\bigl[[i(X),\delta],i(Y)\bigr]\,.
 \end{equation*}
\end{lemma}

\begin{proof}
The map defined by the right hand side of the above equality,
  \begin{equation*}D:\eta\mapsto
D(\eta)=\bigl[[i(X),\delta],i(Y)\bigr]\end{equation*} is a
derivation of degree $-1$ of $\Omega(M,E)$, since it is obtained
by repeated application of the graded bracket to derivations
(property~\ref{properties of derivations}~(ii)). Therefore, it vanishes on
$\Omega^0(M,E)=C^\infty(M,\RR)$. As a consequence, that map is
$C^\infty(M,\RR)$-linear; we have indeed, for each $f\in
C^\infty(M,\RR)$ and $\eta\in\Omega(M,R)$,
 \begin{equation*}D(f\eta)=D(f)\wedge\eta+fD(\eta)=fD(\eta)\,.\end{equation*}
Therefore, there exists a unique smooth section $[X,Y]_\delta$ of
$(E,\tau,M)$ such that, for each $\eta\in \Omega^1(M,E)$,
 \begin{equation*}\bigl\langle
\eta,[X,Y]_\delta\bigr\rangle=D(\eta)\,.\end{equation*} Now the
maps
 \begin{equation*}i\bigl([X,Y]_\delta\bigr)\quad\hbox{and}\quad
   \bigl[[i(X),\delta],i(Y)\bigr]\end{equation*}
are both derivations of degree $-1$ of $\Omega(M,E)$, which
coincide on $\Omega^0(M,E)$ and $\Omega^1(M,E)$. Since derivations
are local, and since $\Omega(M,E)$ is locally generated by its
elements of degrees $0$ and $1$, these two derivations are equal.
\end{proof}

\begin{lemma}
\label{deltabracket 2}
Under the same assumptions as those of 
Lemma~\ref{deltabracket}, let us set,
for each smooth section $X$ of $(E,\tau,M)$,
 \begin{equation*}{\cal L}_\delta(X)
 =\bigl[i(X),\delta\bigr]\,.
 \end{equation*}
Then, for each smooth section $X$ of $(E,\tau,M)$, we have
 \begin{equation*}\bigl[{\cal
L}_\delta(X),\delta\bigr]=\bigl[i(X),\delta^2\bigr]\,;\end{equation*}
for each pair $(X,Y)$ of smooth sections of $(E,\tau,M)$, we have
 \begin{equation*}\bigl[{\cal L}_\delta(X),{\cal L}_\delta(Y)\bigr]
 -{\cal L}_\delta\bigl([X,Y]_\delta\bigr)
 =\Bigl[\bigl[i(X),\delta^2\bigr],i(Y)\Bigr]\end{equation*}
and, for each triple $(X,Y,Z)$ of smooth sections of $(E,\tau,M)$,
we have
 \begin{equation*}i\Bigl(\bigl[X,[Y,Z]_\delta\bigr]_\delta+
          \bigl[Y,[Z,X]_\delta\bigr]_\delta+
          \bigl[Z,[X,Y]_\delta\bigr]_\delta\Bigr)
    =\Bigl[\bigl[[i(X),\delta^2\bigr],
    i(Y)\bigr],i(Z)\Bigr]\,.\end{equation*}
\end{lemma}

\begin{proof}
Let us use the graded Jacobi identity. We may write
 \begin{equation*}\bigl[{\cal L}_\delta(X),\delta\bigr]
   =\bigl[[i(X),\delta],\delta\bigr]
   =-\bigl[[\delta,\delta],i(X)\bigr]
    -\bigl[[\delta,i(X)],\delta\bigr]\,.\end{equation*}
Since $[\delta,\delta]=2\delta^2$, we obtain
 \begin{equation*}2\bigl[{\cal
L}_\delta(X),\delta\bigr]=-2\bigl[\delta^2,i(X)\bigr]
    =2\bigl[i(X),\delta^2\bigr]\,,\end{equation*}
which proves the first equality. Similarly, by using again the
graded Jacobi identity and the equality just proven,
  \begin{equation*}
  \begin{split}
  \bigl[{\cal L}_\delta(X),{\cal L}_\delta(Y)\bigr]
  &=\bigl[{\cal L}_\delta(X),[i(Y),\delta]\bigr]\\
  &=-\bigl[i(Y),[\delta,{\cal L}_\delta(X)]\bigr]
      +\bigl[\delta,[{\cal L}_\delta(X),i(Y)]\bigr]\\
  &=-\bigl[[\delta,{\cal L}_\delta(X)],i(Y)\bigr]
      +\bigl[[{\cal L}_\delta(X),i(Y)],\delta\bigr]\\
  &=\bigl[[{\cal L}_\delta(X),\delta],i(Y)\bigr]
      +\bigl[i\bigl([X,Y]_\delta\bigr),\delta\bigr]\\
  &=\bigl[[i(X),\delta^2],i(Y)\bigr]
      +{\cal L}_\delta\bigl([X,Y]_\delta\bigr)\,.
 \end{split}
  \end{equation*}
The second formula is proven. Finally,
  \begin{equation*}
  \begin{split}
  i\Bigl(\bigl[X,[Y,Z]_\delta\bigr]_\delta\Bigr)
  &=\bigl[{\cal L}_\delta(X),i\bigl([Y,Z]_\delta\bigr)\bigr]\\
  &=\bigl[{\cal L}_\delta(X),[{\cal L}_\delta(Y),i(Z)]\bigr]\\
  &=-\bigl[{\cal L}_\delta(Y),[i(Z),{\cal L}_\delta(X)]\bigr]
    -\bigl[i(Z),[{\cal L}_\delta(X),{\cal L}_\delta(Y)]\bigr]\\
  &=\bigl[{\cal L}_\delta(Y),[{\cal L}_\delta(X),i(Z)]\bigr]
    -\bigl[i(Z),{\cal L}_\delta\bigl([X,Y]_\delta\bigr)\bigr]\\
  &\quad-\Bigl[i(Z),\bigl[[i(X),\delta^2],i(Y)\bigr]\Bigr]\\
  &=i\bigl(\bigl[Y,[X,Z]_\delta\bigr]_\delta\bigr)
     +i\bigl(\bigl[[X,Y]_\delta,Z\bigr]_\delta\bigr)\\
  &\quad+\Bigl[\bigl[[i(X),\delta^2],
    i(Y)\bigr],i(Z)\Bigr]\,.
 \end{split}
 \end{equation*}
The proof is complete.
\end{proof}

\begin{theorem}
\label{derivation on dual bundle}
Let $(E,\tau,M)$ be a vector bundle and $(E^*,\pi,M)$ its dual
bundle. Let $\delta$ be a graded derivation of degree $1$ of the
exterior algebra $\Omega(M,E)$ (notations defined 
in~\ref{The exterior algebra of sections}), which
satisfies
 \begin{equation*}\delta^2=\delta\circ\delta=0\,.\end{equation*}
Then $\delta$ determines a natural Lie algebroid structure on
$(E,\tau,M)$. For that structure, the anchor map $\rho:E\to TM$ is
the unique vector bundle map such that, for each smooth section
$X$ of $(E,\tau,M)$ and each function $f\in C^\infty(M,\RR)$,
 \begin{equation*}i(\rho\circ X)\,df=\langle\delta f,X\rangle\,.\end{equation*}
The bracket $(X,Y)\mapsto\{X,Y\}$ is the $\delta$-bracket defined
in Lemma~\ref{deltabracket}; it is such that, for each pair $(X,Y)$ of smooth
sections of $(E,\tau,M)$,
 \begin{equation*}i\bigl(\{X,Y\}\bigr)
 =\bigl[[i(X),\delta],i(Y)\bigr]\,.
 \end{equation*}
The $\omega(M,E)$-valued exterior derivative associated to that
Lie algebroid structure 
(propositions \ref{exterior derivative of forms 1} 
and \ref{exterior derivative of forms 2}) is the
given derivation $\delta$.
\end{theorem}

\begin{proof}
Since $\delta^2=0$, lemmas~\ref{deltabracket} 
and~\ref{deltabracket 2} prove that the
$\delta$-bracket satisfies the Jacobi identity. Let $X$ and $Y$ be
two smooth sections of $(E,\tau,M)$ and $f$ a smooth function on
$M$. By using the definition of the $\delta$-bracket we obtain
 \begin{equation*}i\bigl([X,fY]_\delta\bigr)=f\,i\bigl([X,Y]_\delta\bigr)+
 \bigl({\cal L}_\delta(X)f\bigr)\,i(Y)\,.\end{equation*}
But
 \begin{equation*}{\cal L}_\delta(X)f=\bigl[i(X),\delta\bigr]f=\langle\delta
 f,X\rangle\,,\end{equation*}
since $i(X)f=0$. We must prove now that the value of $\delta(f)$
at any point $x\in M$ depends only on the value of the
differential $df$ of $f$ at that point. We first observe that
$\delta$ being a derivation, the values of $\delta(f)$ in some
open subset $U$ of $M$ depend only on the values of $f$ in that
open subset. Moreover, we have
 \begin{equation*}\delta(1\,f)=\delta(f)=\delta(1)\,f+1\,\delta(f)=\delta(1)\,f+
\delta(f)\,,\end{equation*} which proves that $\delta$ vanishes on
constants.
\par
Let $a\in M$. We use a chart of $M$ whose domain $U$ contains $a$,
and whose local coordinates are denoted by $(x^1,\ldots,x^n)$. In
order to calculate $\delta(f)(a)$, the above remarks allow us to
work in that chart. We may write
 \begin{equation*}f(x)=f(a)+\sum_{i=1}^n(x^i-a^i)\varphi_i(x)\,,\quad\hbox{with}
\quad
        \lim_{x\to a}\varphi_i(x)=\frac{\partial f}{\partial
x^i}(a)\,.\end{equation*} Therefore,
 \begin{equation*}(\delta f)(a)=\sum_{i=1}^n\frac{\partial f}{\partial
x^i}(a)\,\delta(x^i)(a)\,.
 \end{equation*}
We have proven that $\delta(f)(a)$ depends only on $df(a)$, and
that we may write
 \begin{equation*}\delta(f)={}^t\!\rho_\delta\circ df\,,\end{equation*}
where ${}^t\!\rho_\delta:T^*M\to E^*$ is a smooth vector bundle
map. Let $\rho_\delta:E\to TM$ be its transpose. We may now write
  \begin{equation*}[X,fY]_\delta=f[X,Y]_\delta+\langle df,\rho_\delta\circ
X\rangle
  Y\,.\end{equation*}
This proves that the vector bundle $(E,\tau,M)$, with the
$\delta$-bracket and the map $\rho_\delta$ as anchor, is a Lie
algebroid. Finally, by using 
Propositions~\ref{exterior derivative of forms 1}
and~\ref{exterior derivative of forms 2}, we see
that the $\Omega(M,E)$-valued exterior derivative associated to
that Lie algeboid structure is the derivation $\delta$.
\end{proof}

\section{The Schouten-Nijenhuis bracket}

In this subsection $(E,\tau,M,\rho)$ is a Lie algebroid. We have
seen (Propositions~\ref{Lie derivative by section of algebroid 3} 
and~\ref{Lie derivative of multivectors 2})
that the composition law wich
associates, to each pair $(V,W)$ of sections of the Lie algebroid
$(E,\tau,M,\rho)$, the bracket $\{V,W\}$, extends into a map
$(V,P)\mapsto{\cal L}_\rho(V) P$, defined on $A^1(M,E)\times
A(M,E)$, with values in $A(M,E)$. Theorem~\ref{Schouten bracket} below will show
that this map extends, in a very natural way, into a composition
law $(P,Q)\mapsto[P,Q]$, defined on $A(M,E)\times A(M,E)$, with
values in $A(M,E)$, called the {\it Schouten-Nijenhuis bracket\/}.
That bracket was discovered by Schouten \cite{Schou} for
multivectors on a manifold, and its properties were further
studied by Nijenhuis \cite{Ni}.
\par\smallskip
The following lemmas will be used in the proof of 
Theorem~\ref{Schouten bracket}.

\begin{lemma} 
\label{Schouten lemma}
Let $(E,\tau,M,\rho)$ be a Lie algebroid, $p$ and
$q\in\ZZ$, $P\in A^p(M,E)$, $Q\in A^q(M,E)$, $f\in
C^\infty(M,\RR)$ and $\eta\in \Omega(M,E)$. Then
 \begin{equation*}
 \begin{split}
 i(P)\bigl(df\wedge i(Q)\eta\bigr)&-(-1)^pdf\wedge\bigl(i(P)\circ
 i(Q)\eta\bigr)\\
 &+(-1)^{(p-1)q}i(Q)\circ i(P)(df\wedge\eta)\\
 &+(-1)^{(p-1)q+p}i(Q)\bigl(df\wedge i(P)\eta\bigr)\\
 &=0\,.
 \end{split}
 \end{equation*}

\end{lemma}

\begin{proof} Let us denote by $E(P,Q,f,\eta)$ the left hand side of
the above equality. We have to prove that $E(P,Q,f,\eta)=0$.
\par
Obviously, $E(P,Q,f,\eta)=0$ when $p<0$, as well as when $q<0$.
When $p=q=0$, we have
 \begin{equation*}E(P,Q,f,\eta)=PQ\,df\wedge\eta-PQ\,df\wedge\eta
 -QP\,df\wedge\eta+QP\,df\wedge\eta=0\,.\end{equation*}
Now we proceed by induction on $p$ and $q$, with the induction
assumption that $E(P,Q,f,\eta)=0$ when $p\leq p_M$ and $q\leq
q_M$, for some integers $p_M$ and $q_M$. Let $P=X\wedge P'$, with
$X\in A^1(M,E)$ and $P'\in A^{p_M}(M,E)$, $Q\in A^q(M,E)$, with
$q\leq q_M$, $f\in C^\infty(M,\RR)$ and $\eta\in\Omega(M,E)$. We
obtain, after some calculations,
 \begin{equation*}
 \begin{split}
 E(P,Q,f,\eta)
 &=E(X\wedge P',Q,f,\eta)\\
 &=(-1)^{p_M+q-1}E\bigl(P',Q,f,i(X)\eta\bigr)\\
 &\quad+(-1)^{p_M}\langle df,X\rangle i(P)\circ i(Q)\eta\\
 &\quad-(-1)^{p_M+p_Mq} \langle df,X\rangle i(Q)\circ i(P)\eta\\
 &=0\,,
 \end{split}
 \end{equation*}
since, by the induction assumption,
$E\bigl(P',Q,f,i(X)\eta\bigr)=0$.
\par
Since every $P\in A^{p_M+1}(M,E)$ is the sum of terms of the form
$X\wedge P'$, with $X\in A^1(M,E)$ and $P'\in A^{p_M}(M,E)$, we
see that $E(P,Q,f,\eta)=0$ for all $p\leq p_M+1$, $q\leq q_M$,
$P\in A^{p_M+1}(M,E)$ and $Q\in A^{q_M}(M,E)$.
\par
Moreover, $P$ and $Q$ play similar parts in $E(P,Q,f,\eta)$, since
we have
 \begin{equation*}E(P,Q,f,\eta)=(-1)^{pq+p+q}E(Q,P,f,\eta)\,.\end{equation*}
Therefore  $E(P,Q,f,\eta)=0$ for all $p\leq p_M+1$, $q\leq q_M+1$,
$P\in A^{p}(M,E)$ and $Q\in A^{q}(M,E)$. By induction we conclude
that  $E(P,Q,f,\eta)=0$ for all $p$ and $q\in\ZZ$, $P\in
A^{p}(M,E)$ and $Q\in A^{q}(M,E)$.
\end{proof}

\begin{lemma} 
\label{Schouten lemma 2}
Let $(E,\tau,M,\rho)$ be a Lie algebroid, $p$, $q$ and
$r\in\ZZ$, $P\in A^p(M,E)$, $Q\in A^q(M,E)$ and $R\in A^r(M,E)$.
Then
 \begin{equation*}i(R)\circ\bigl[[i(P),d_\rho],i(Q)\bigr]=(-1)^{(p+q-1)r}
 \bigl[[i(P),d_\rho],i(Q)\bigr]\circ i(R)\,.\end{equation*}

\end{lemma}

\begin{proof} Let us first assume that $R=V\in A^1(M,E)$. We may write
 \begin{equation*}
 \begin{split}
 i(V)\circ\bigl[[i(P),d_\rho],i(Q)\bigr]
 &=i(V)\circ i(P)\circ d_\rho\circ i(Q)\\
 &\quad-(-1)^pi(V)\circ d_\rho\circ
 i(P)\circ i(Q)\\
 &\quad-(-1)^{(p-1)q}i(V)\circ i(Q)\circ i(P)\circ d_\rho\\
 &\quad+(-1)^{(p-1)q+p} i(V)\circ i(Q)\circ d_\rho\circ i(P)\,.
 \end{split}
 \end{equation*}
We transform the right hand side by pushing the operator $i(V)$
towards the right, using the formulae (proven 
in~\ref{operations in exterior powers}~(v) and in
Property~1 of~\ref{exterior derivative of forms 2})
 \begin{equation*}i(V)\circ i(P)=(-1)^pi(P)\circ i(V)\quad\hbox{\rm and}\quad
 i(V)\circ d_\rho={\cal L}_\rho(V)-d_\rho\circ i(V)\,.\end{equation*}
We obtain, after rearrangement of the terms,
 \begin{equation*}
 \begin{split}
 i(V)\circ\bigl[[i(P),d_\rho],i(Q)\bigr]
 &=(-1)^{p+q-1} \bigl[[i(P),d_\rho],i(Q)\bigr]\circ i(V)\\
 &\quad+(-1)^pi(P)\circ{\cal L}_\rho(V)\circ i(Q)\\
 &\quad-(-1)^p{\cal L}_\rho(V)\circ i(P)\circ i(Q)\\
 &\quad-(-1)^{(p-1)q+p+q}i(Q)\circ i(P)\circ {\cal L}_\rho(V)\\
 &\quad+(-1)^{(p-1)q+p+q}i(Q)\circ{\cal L}_\rho(V)\circ i(P)\,.
 \end{split}
 \end{equation*}
Now we transform the last four terms of the right hand side by
pushing to the right the operator ${\cal L}_\rho(V)$, using
formulae, proven in~\ref{operations in exterior powers}~(v) and in Property~4
of~\ref{Lie derivative of multivectors 2}, of the
type
 \begin{equation*}i(P)\circ i(Q)=i(P\wedge Q)\quad\hbox{\rm and}\quad
 {\cal L}_\rho(V)\circ i(P)=i(P)\circ{\cal L}_\rho(V)+i\bigl({\cal
 L}_\rho(V)P\bigr)\,.\end{equation*}
The terms containing ${\cal L}_\rho(V)$ become
 \begin{equation*}(-1)^pi\Bigl(P\wedge {\cal L}_\rho(V)Q+\bigl({\cal
 L}_\rho(V)P\bigr)\wedge Q-{\cal L}_\rho(V)(P\wedge Q)\Bigr)\,,\end{equation*}
so they vanish, by Property~3 of~\ref{Lie derivative of multivectors 2}.
So we have
\begin{equation*}i(V)\circ\bigl[[i(P),d_\rho],i(Q)\bigr]=(-1)^{(p+q-1)}
 \bigl[[i(P),d_\rho],i(Q)\bigr]\circ i(V)\,.\end{equation*}
Now let $R=V_1\wedge \cdots\wedge V_r$ be a decomposable element
in $A^r(M,E)$. Since
 \begin{equation*}i(R)=i(V_1)\circ\cdots\circ i(V_r)\,,\end{equation*}
by using $r$ times the above result, we obtain
\begin{equation*}i(R)\circ\bigl[[i(P),d_\rho],i(Q)\bigr]=(-1)^{(p+q-1)r}
 \bigl[[i(P),d_\rho],i(Q)\bigr]\circ i(R)\,.\end{equation*}
Finally the same result holds for all $R\in A^r(M,E)$ by
linearity. \end{proof}

\begin{theorem}
\label{Schouten bracket} 
Let $(E,\tau,M,\rho)$ be a Lie algebroid. Let $p$ and
$q\in\ZZ$, and let $P\in A^p(M,E)$, $Q\in A^q(M,E)$. There exists
a unique element in $A^{p+q-1}(M,E)$, called the {\it
Schouten-Nijenhuis bracket\/} of $P$ and $Q$, and denoted by
$[P,Q]$, such that the interior product $i\bigl([P,Q]\bigr)$,
considered as a graded endomorphism of degree $p+q-1$ of the
exterior algebra $\Omega(M,E)$, is given by the formula
 \begin{equation*}i\bigl([P,Q]\bigr)
 =\bigl[[i(P),d_\rho],i(Q)\bigr]\,,
 \end{equation*}
the brackets in the right hand side being the graded brackets of
graded endomorphism (Definition~\ref{graded bracket}).
\end{theorem}

\begin{proof} We observe that for all $r\in\ZZ$, the map
\begin{equation*}\eta\mapsto
\bigl[[i(P),d_\rho],i(Q)\bigr]\eta\,,\end{equation*} defined on
$\Omega^r(M,E)$, with values in $\Omega^{r-p-q+1}(M,E)$, is
$\RR$-linear. Let us prove that it is in fact
$C^\infty(M,\RR)$-linear. Let $f\in C^\infty(M,\RR)$. By
developing the double graded bracket of endomorphisms, we obtain
after some calculations
 \begin{equation*}
 \begin{split}
 \bigl[[i(P),d_\rho],i(Q)\bigr](f\eta)
 &=f\bigl[[i(P),d_\rho],i(Q)\bigr]\eta\\
 &\quad+i(P)\bigl(df\wedge i(Q)\eta\bigr)-(-1)^pdf\wedge\bigl(i(P)\circ
 i(Q)\eta\bigr)\\
 &\quad+(-1)^{(p-1)q}i(Q)\circ i(P)(df\wedge\eta)\\
 &\quad+(-1)^{(p-1)q+p}i(Q)\bigl(df\wedge i(P)\eta\bigr)\,.
 \end{split}
 \end{equation*}
Lemma~\ref{Schouten lemma} 
shows that the sum of the last four terms of the right
hand side vanishes, so we obtain
 \begin{equation*}\bigl[[i(P),d_\rho],i(Q)\bigr](f\eta)
 =f\bigl[[i(P),d_\rho],i(Q)\bigr]\eta\,.\end{equation*}
Let us take $r=p+q-1$, and $\eta\in \Omega^{p+q-1}(M,E)$. The map
 \begin{equation*}\eta\mapsto
\bigl[[i(P),d_\rho],i(Q)\bigr]\eta\,,\end{equation*} defined on
$\Omega^{p+q-1}(M,E)$, takes its values in
$\Omega^0(M,E)=C^\infty(M,\RR)$, and is $C^\infty(M,\RR)$-linear.
This proves the existence of a unique element $[P,Q]$ in
$\Omega^{p+q-1}(M,E)$ such that, for all $\eta\in
\Omega^{p+q-1}(M,E)$,
 \begin{equation*}\bigl[[i(P),d_\rho],i(Q)\bigr]\eta
 =i\bigl([P,Q]\bigr)\eta\,.
 \end{equation*}
We still have to prove that the same formula holds for all
$r\in\ZZ$ and all $\eta\in\Omega^r(M,E)$. The formula holds
trivially when $r<p+q-1$, so let us assume that $r>p+q-1$. Let
$\eta\in\Omega^r(M,E)$ and $R\in A^{r-p-q+1}(M,E)$. By using
Lemma~\ref{Schouten lemma 2}, we may write
 \begin{equation*}
 \begin{split}
 i(R)\circ\bigl[[i(P),d_\rho],i(Q)\bigr](\eta)
 &=(-1)^{(p+q-1)(r-p-q+1)}\bigl[[i(P),d_\rho],i(Q)\bigr]\bigl(i(R)\eta\bigr)\\
 &=(-1)^{(p+q-1)(r-p-q+1)}i\bigl([P,Q]\bigr)\bigl(i(R)\eta)\,,
 \end{split}
 \end{equation*}
since $i(R)\eta\in\Omega^{p+q-1}(E)$. Therefore
 \begin{equation*}
 \begin{split}
 i(R)\circ\bigl[[i(P),d_\rho],i(Q)\bigr](\eta)
 &=(-1)^{(p+q-1)(r-p-q+1)}i\bigl([P,Q]\bigr)\circ i(R)\eta)\\
 &=i(R)\circ i\bigl([P,Q]\bigr)\eta\,.
 \end{split}
 \end{equation*}
Since that equality holds for all $\eta\in\Omega^r(M,E)$ and all
$R\in A^{r-p-q+1}(M,E)$, we may conclude that
 \begin{equation*}\bigl[[i(P),d_\rho],i(Q)\bigr]
 =i\bigl([P,Q]\bigr)\,,
 \end{equation*}
and the proof is complete.
\end{proof}
\par\medskip
In Proposition~\ref{Lie derivative by section of algebroid 2}, 
we introduced the Lie derivative with
respect to a section of the Lie algebroid $(E,\tau,M,\rho)$. Now we
define, for all $p\in\ZZ$ and $P\in A^p(M,E)$, the Lie derivative
with respect to $P$. The reader will observe that Property~1 of
Proposition~\ref{exterior derivative of forms 2} shows that for $p=1$, the
following definition is in agreement with the definition of the Lie derivative
with respect to an element in $A^1(M,E)$ given 
in~\ref{Lie derivative by section of algebroid 2}.

\begin{defi}
\label{Lie derivative by multivector}
Let $(E,\tau,M,\rho)$ be a Lie algebroid, $p\in\ZZ$ and $P\in
A^p(M,E)$. The {\it Lie derivative\/} with respect to $P$ is the
graded endomorphism of $\Omega(M,P)$, of degree $1-p$, denoted by
${\cal L}_\rho(P)$,
 \begin{equation*}{\cal L}_\rho(P)=\bigl[i(P),d_\rho\bigr]=i(P)\circ d_\rho
 -(-1)^pd_\rho\circ i(P)\,.\end{equation*}

\end{defi}

\begin{rmk} 
\label{Schouten remark}
Under the assumptions of Theorem~\ref{Schouten bracket}, the above
Definition allows us to write
  \begin{equation*}i\bigl([P,Q]\bigr)=\bigl[{\cal L}_\rho(P),i(Q)\bigr]
  ={\cal L}_\rho(P)\circ i(Q)-(-1)^{(p-1)q}i(Q)\circ{\cal
  L}_\rho(P)\,.\end{equation*}
For $p=1$ and $P=V\in A^1(M,E)$,  this formula is simply Property
4 of Proposition~\ref{Lie derivative of multivectors 2}, as shown by the
following Proposition.
\end{rmk}

\begin{prop} 
\label{Schouten bracket and Lie derivative}
Under the assumptions of Theorem~\ref{Schouten bracket}, let
$p=1$, $P=V\in A^1(M,E)$ and $Q\in A^q(M,E)$. The
Schouten-Nijenhuis bracket $[V,Q]$ is simply the Lie derivative of
$Q$ with respect to $V$, as defined in 
Proposition~\ref{Lie derivative of multivectors}:
 \begin{equation*}[V,Q]={\cal L}_\rho(V)Q\,.\end{equation*}
\end{prop}

\begin{proof} As seen in Remark~\ref{Schouten remark}, we may write
 \begin{equation*}i\bigl([V,Q]\bigr)=\bigl[{\cal L}_\rho(V),i(Q)\bigr]
  ={\cal L}_\rho(V)\circ i(Q)-i(Q)\circ{\cal
  L}_\rho(V)\,.\end{equation*}
Property~4 of Proposition~ \ref{Lie derivative of multivectors 2} shows that
 \begin{equation*}i\bigl({\cal L}_\rho(V)Q\bigr)={\cal L}_\rho(V)\circ i(Q)
 -i(Q)\circ{\cal L}_\rho(V)\,.\end{equation*}
Therefore,
 \begin{equation*}i\bigl([V,Q]\bigr)=i\bigl({\cal
L}_\rho(V)Q\bigr)\,,\end{equation*} and finally
 \begin{equation*}[V,Q]={\cal L}_\rho(V)Q\,,\end{equation*}
which ends the proof.
\end{proof}

\begin{rmks}\hfill
\par\nobreak\smallskip\noindent
{\rm(i)}\quad{\it The Lie derivative of elements in
$A^p(M,E)$\/.}\quad One may think to extend the range of
application of the Lie derivative with respect to a multivector
$P\in A^p(M,E)$ by setting, for all $q\in\ZZ$ and $Q\in A^q(M,E)$,
 \begin{equation*}{\cal L}_\rho(P)Q=[P,Q]\,,\end{equation*}
the bracket in the right hand side being the Schouten-Nijenhuis
bracket. However, we will avoid the use of that notation because
it may lead to confusions: for $p>1$, $P\in A^p(M,E)$, $q=0$ and
$Q=f\in A^0(M,E)=C^\infty(M,\RR)$, the Schouten-Nijenhuis bracket
$[P,f]$ is an element in $A^{p-1}(M,E)$ which does not vanish in
general. But $f$ can be considered also as an element in
$\Omega^0(M,E)$, and the Lie derivative of $f$ with respect to
$P$, in the sense of Definition~\ref{Lie derivative by multivector}, is an
element in
$\Omega^{-(p-1)}(M,E)$, therefore vanishes identically. So it
would not be a good idea to write ${\cal L}_\rho(P)f=[P,f]$.
\par\smallskip\noindent
{\rm(ii)}\quad{Lie derivatives and derivations\/.}\quad We have
seen (Property~3 of~\ref{Lie derivative by section of algebroid 3}) that the Lie
derivative ${\cal
L}_\rho(V)$ with respect to a section $V\in A^1(M,R)$ of the Lie
algebroid $(E,\tau,M,\rho)$ is a derivation of degree $0$ of the
exterior algebra $\Omega(M,E)$. For $p>1$ and $P\in A^p(M,E)$, the
Lie derivative ${\cal L}_\rho(P)$ is a graded endomorphism of
degree $-(p-1)$ of $\Omega(M,E)$. Therefore, it vanishes
identically on $\Omega^0(M,E)$ and on $\Omega^1(M,E)$. Unless it
vanishes identically, ${\cal L}_\rho(P)$ is not a derivation of
$\Omega(M,E)$.
\end{rmks}

\begin{prop} 
\label{Lie derivative by multivector 2}
Let $(E,\tau,M,\rho)$ be a Lie algebroid,
$p$ and $q\in\ZZ$, $P\in A^p(M,E)$, $Q\in A^q(M,E)$.
\par\nobreak\smallskip
{\rm 1.} The graded bracket of the Lie derivative ${\cal
L}_\rho(P)$ and the exterior differential $d_\rho$ vanishes
identically:
 \begin{equation*}\bigl[{\cal L}_\rho(P),d_\rho\bigr]={\cal L}_\rho(P)\circ
 d_\rho-(-1)^{p-1}d_\rho\circ{\cal L}_\rho(P)=0\,.\end{equation*}
\par\smallskip
{\rm 2.} The graded bracket of the Lie derivatives ${\cal
L}_\rho(P)$ and ${\cal L}_\rho(Q)$ is equal to the Lie derivative
${\cal L}_\rho\bigl([P,Q]\bigr)$:
 \begin{equation*}\bigl[{\cal L}_\rho(P),{\cal L}_\rho(Q)\bigr]=
 {\cal L}_\rho(P)\circ {\cal L}_\rho(Q)-(-1)^{(p-1)(q-1)}
 {\cal L}_\rho(Q)\circ{\cal L}_\rho(P)={\cal
 L}_\rho\bigl([P,Q]\bigr)\,.\end{equation*}

\end{prop}

\begin{proof}
 {\rm 1.} We have seen~(\ref{examples graded Lie algebras}~(ii)) that the space
of graded
endomorphisms of $\Omega(M,E)$, with the graded bracket as
composition law, is a graded Lie algebra. By using the graded
Jacobi identity, we may write
 \begin{equation*}(-1)^p\bigl[[i(P),d_\rho],d_\rho\bigr]
  +(-1)^p\bigl[[d_\rho,d_\rho],i(P)\bigr]
  -\bigl[[d_\rho,i(P)],d_\rho\bigr]=0\,.\end{equation*}
But
 \begin{equation*}[d_\rho,d_\rho]=2 d_\rho\circ d_\rho=0\quad
 \hbox{\rm and}\quad
\bigl[i(P),d_\rho\bigr]=-(-1)^p\bigl[d_\rho,i(P)\bigr]\,.\end{equation*}
So we obtain
 \begin{equation*}2\bigl[[i(P),d_\rho],d_\rho\bigr]=2\bigl[{\cal
 L}_\rho(P),d_\rho\bigr]=0\,.\end{equation*}
\par\smallskip
{\rm 2.} We have
 \begin{equation*}{\cal
 L}_\rho\bigl([P,Q]\bigr)=\bigl[i\bigl([P,Q]\bigr),d_\rho\bigr]
 =\bigl[[{\cal L}_\rho(P),i(Q)],d_\rho\bigr]\,.\end{equation*}
Using the graded Jacobi identity, we may write
 \begin{equation*}
 \begin{split}
 (-1)^{p-1}\bigl[[{\cal L}_\rho(P),i(Q)],d_\rho\bigr]
 &+(-1)^{q(p-1)}\bigl[[i(Q),d_\rho],{\cal L}_\rho(P)\bigr]\\
 &+(-1)^q\bigl[[d_\rho,{\cal L}_\rho(P)],i(Q)\bigr]=0\,.
 \end{split}
 \end{equation*}
But, according to~\ref{Lie derivative by multivector} and Property~1 above,
 \begin{equation*}\bigl[i(Q),d_\rho\bigr]={\cal L}_\rho(Q)\quad\hbox{\rm
 and}\quad\bigl[d_\rho,{\cal L}_\rho(P)\bigr]=0\,.\end{equation*}
So we obtain
 \begin{equation*}{\cal L}_\rho\bigl([P,Q]\bigr)=-(-1)^{(p-1)(q-1)}
 \bigl[{\cal L}_\rho(Q),{\cal L}_\rho(P)\bigr]
 =\bigl[{\cal L}_\rho(P),{\cal L}_\rho(Q)\bigr]\,,\end{equation*}
as announced.
\end{proof}

\begin{prop}
\label{Schouten bracket 2}
 Under the same assumptions as those of
Theorem~\ref{Schouten bracket}, the Schouten-Nijenhuis bracket has the following
properties.
 \par\nobreak\smallskip

{\rm 1.} For $f$ and $g\in A^0(M,E)=C^\infty(M,\RR)$,
 \begin{equation*}[f,g]=0\,.\end{equation*}
\par\smallskip

{\rm 2.} For $V\in A^1(M,E)$, $q\in\ZZ$ and $Q\in A^q(M,E)$,
 \begin{equation*}[V,Q]={\cal L}_\rho(V)Q\,.\end{equation*}
\par\smallskip

{\rm 3.} For $V$ and $W \in A^1(M,E)$,
 \begin{equation*}[V,W]=\{V,W\}\,,\end{equation*}
the bracket in the right hand side being the Lie algebroid
bracket.

{\rm 4.} For all $p$ and $q\in\ZZ$, $P\in A^p(M,E)$, $Q\in
A^q(M,E)$,
 \begin{equation*}[P,Q]=-(-1)^{(p-1)(q-1)}[Q,P]\,.\end{equation*}
\par\smallskip

{\rm 5.} Let $p\in\ZZ$, $P\in A^p(M,E)$. The map $Q\mapsto [P,Q]$
is a derivation of degree $p-1$ of the graded exterior algebra
$A(M,E)$. In other words, for $q_1$ and $q_2\in\ZZ$, $Q_1\in
A^{q_1}(M,E)$ and $Q_2\in A^{q_2}(M,E)$,
 \begin{equation*}[P,Q_1\wedge Q_2]=[P,Q_1]\wedge
 Q_2+(-1)^{(p-1)q_1}Q_1\wedge[P,Q_2]\,.\end{equation*}
\par\smallskip

{\rm 6.} Let $p$, $q$ and $r\in\ZZ$, $P\in A^p(M,E)$, $Q\in
A^q(M,E)$ and $R\in A^r(M,E)$. The Schouten-Nijenhuis bracket
satisfies the graded Jacobi identity:
 \begin{equation*}
 \begin{split}
 (-1)^{(p-1)(r-1)}\bigl[[P,Q],R\bigr]
 &+(-1)^{(q-1)(p-1)}\bigl[[Q,R],P\bigr]\\
 &+(-1)^{(r-1)(q-1)}\bigl[[R,P],Q\bigr]\\
 &=0\,.
 \end{split}
 \end{equation*}

\end{prop}

\begin{proof}
 {\rm 1.} Let $f$ and $g\in A^0(M,E)$. Then $[f,g]\in
 A^{-1}(M,E)=\{0\}$, therefore $[f,g]=0$.
\par\smallskip

{\rm 2.} See Proposition~\ref{Schouten bracket and Lie derivative}.
\par\smallskip

{\rm 3.} See Property~1 of Proposition~\ref{Lie derivative of multivectors 2}.
\par\smallskip

{\rm 4.}  Let $p$ and $q\in\ZZ$, $P\in A^p(M,E)$, $Q\in A^q(M,E)$.
By using the graded Jacobi identity for graded endomorphisms of
$\Omega(M,E)$, we may write
 \begin{equation*}(-1)^{pq}\bigl[[i(P),d_\rho],i(Q)\bigr]
 +(-1)^{p}\bigl[[d_\rho,i(Q)],i(P)\bigr]
 +(-1)^{q}\bigl[[i(Q),i(P)],d_\rho\bigr]
 =0\,.\end{equation*}
By using
 \begin{equation*}\bigl[i(Q),i(P)\bigr]=i(Q\wedge P)-i(Q\wedge
 P)=0\quad\hbox{and}\quad
 \bigl[d_\rho,i(Q)\bigr]=-(-1)^{q}\bigl[i(Q),d_\rho\bigr]\,,\end{equation*}
we obtain
 \begin{equation*}(-1)^{pq}\bigl[[i(P),d_\rho],i(Q)\bigr]+(-1)^{p+q-1}
 \bigl[[i(Q),d_\rho],i(P)\bigr]=0\,,\end{equation*}
so the result follows immediately.
\par\smallskip

{\rm 5.} Let $p$, $q_1$ and $q_2\in\ZZ$, $P\in A^p(M,E)$, $Q_1\in
A^{q_1}(M,E)$ and $Q_2\in A^{q_2}(M,E)$. We may write
 \begin{equation*}
 \begin{split}
 i\bigl([P,Q_1\wedge Q_2]\bigr)
 &=\bigl[{\cal L}_\rho(P),i(Q_1\wedge Q_2)\bigr]\\
 &={\cal L}_\rho(P)\circ i(Q_1\wedge Q_2)\\
 &\quad-(-1)^{(p-1)(q_1+q_2)}
 i(Q_1\wedge Q_2)\circ{\cal L}_\rho(P)\,.
 \end{split}
 \end{equation*}
We add and substract $(-1)^{(p-1)q_1}i(Q_1)\circ {\cal
L}_\rho(P)\circ i(Q_1)$ from the last expression, and replace
$i(Q_1\wedge Q_2)$ by $i(Q_1)\circ i(Q_2)$. We obtain
 \begin{equation*}i\bigl([P,Q_1\wedge Q_2]\bigr)=\bigl[{\cal
 L}_\rho(P),i(Q_1)\bigr]\circ i(Q_2)
 +(-1)^{(p-1)q_1}i(Q_1)\circ\bigl[{\cal
 L}_\rho(P),i(Q_2)\bigr]\,.\end{equation*}
The result follows immediately.
\par\smallskip

{\rm 6.}
 Let $p$, $q$ and $r\in\ZZ$, $P\in A^p(M,E)$, $Q\in
A^q(M,E)$ and $R\in A^r(M,E)$. By using Property~2 of
Proposition~\ref{Lie derivative by multivector 2}, we may write
 \begin{equation*}
 \begin{split}
 i\bigl(\bigl[[P,Q],R\bigr]\bigr)
 &=\bigl[{\cal L}_\rho([P,Q]),i(R)\bigr]\\
 &=\bigl[[{\cal L}_\rho(P),{\cal L}_\rho(Q)],i(R)\bigr]\,.
 \end{split}
 \end{equation*}
Using the graded Jacobi identity, we obtain
 \begin{equation*}
 \begin{split}
 (-1)^{(p-1)r}\bigl[[{\cal L}_\rho(P),{\cal L}_\rho(Q)],i(R)\bigr]
  &+(-1)^{(q-1)(p-1)}\bigl[[{\cal L}_\rho(Q),i(R)],{\cal
  L}_\rho(P)\bigr]\\
  &+(-1)^{r(q-1)}\bigl[[i(R),{\cal L}_\rho(P)],{\cal L}_\rho(Q)\bigr]
  =0\,.
  \end{split}
  \end{equation*}
But
 \begin{equation*}
 \begin{split}
 \bigl[[{\cal L}_\rho(Q),i(R)],{\cal L}_\rho(P)\bigr]
 &=\bigl[i\bigl([Q,R]\bigr),{\cal L}_\rho(P)\bigr]\\
 &=-(-1)^{(q+r-1)(p-1)}\,\bigl[{\cal
 L}_\rho(P),i\bigl([Q,R]\bigr)\bigr]\\
 &=-(-1)^{(q+r-1)(p-1)}\,i\bigl(\bigl[P,[Q,R]\bigr]\bigr)\\
 &=(-1)^{(q+r-1)(p-1)+(p-1)(q+r-2)}\,i\bigl(\bigl[[Q,R],P\bigr]\bigr)\\
 &=(-1)^{p-1}\,i\bigl(\bigl[[Q,R],P\bigr]\bigr)\,.
 \end{split}
 \end{equation*}
Similarly,
 \begin{equation*}
 \begin{split}
 \bigl[[i(R),{\cal L}_\rho(P)],{\cal L}_\rho(Q)\bigr]
 &=-(-1)^{(p-1)r}\,\bigl[[{\cal L}_\rho(P),i(R)],{\cal
 L}_\rho(Q)\bigr]\\
 &=-(-1)^{(p-1)r}\,\bigl[i\bigl([P,R]\bigr),{\cal
 L}_\rho(Q)\bigr]\\
 &=(-1)^{(p-1)r+(p+r-1)(q-1)}\,\bigl[{\cal
 L}_\rho(Q),i\bigl([P,R]\bigr)\bigr]\\
 &=(-1)^{(p-1)(r+q-1)+r(q-1)}\,i\bigl(\bigl[Q,[P,R]\bigr]\bigr)\\
 &=-(-1)^{(p-1)q+(q-1)(p-2)}\,
 i\bigl(\bigl[[R,P],Q\bigr]\bigr)\\
 &=-(-1)^{p+q}\,i\bigl(\bigl[[R,P],Q\bigr]\bigr)\,.
 \end{split}
 \end{equation*}
Using the above equalities, we obtain
 \begin{equation*}
 \begin{split}
 (-1)^{(p-1)(r-1)}\,i\bigl(\bigl[[P,Q]),R\bigr]\bigr)
 &+(-1)^{(q-1)(p-1)}\,i\bigl(\bigl[[Q,R],P\bigr]\bigr)\\
 &+(-1)^{(r-1)(q-1)}\,i\bigl(\bigl[[R,P],Q\bigr]]\bigr)\\
 &=0\,.
 \end{split}
 \end{equation*}
 The proof is complete.
 \end{proof}

\begin{rmks}
\label{remarks about degrees}
 Let $(E,\tau,M,\rho)$ be a Lie algebroid.
\par\nobreak\smallskip\noindent
{\rm(i)}\quad{\it Degrees for the two algebra structures of
$A(M,E)$\/}.\quad The algebra
$A(M,E)=\bigoplus_{p\in\ZZ}A^p(M,E)$ of sections of the exterior
powers $(\bigwedge^pE,\tau,M)$, with the exterior product as
composition law, is a graded associative algebra; for that
structure, the space of homogeneous elements of degree $p$ is
$A^p(M,E)$. Proposition~\ref{Schouten bracket 2} 
shows that $A(M,E)$, with the
Schouten-Nijenhuis bracket as composition law, is a graded Lie
algebra; for that structure, the space of homogeneous elements of
degree $p$ is not $A^p(M,E)$, but rather $A^{p+1}(M,E)$. For
homogeneous elements in $A(M,E)$, one should therefore make a
distinction between the degree for the graded associative algebra
structure and the degree for the graded Lie algebra structure; an
element in $A^p(M,E)$ has degree $p$ for the graded associative
algebra structure, and degree $p-1$ for the graded Lie algebra
structure.
\par\nobreak\smallskip\noindent
{\rm(ii)}\quad{\it The anchor as a graded Lie algebras
homomorphism\/}.\quad The anchor $\rho:E\to TM$ allows us to
associate to each smooth section $X\in A^1(M,E)$ a smooth vector
field $\rho\circ X$ on $M$; according to 
Definition~\ref{Lie algebroid}, that
correspondence is a Lie algebras homomorphism. We can extend that
map, for all $p\geq 1$, to the space $A^p(M,E)$ of smooth sections
of the $p$-th external power $(\bigwedge^pE,\tau,M)$. First, for a
decomposable element $X_1\wedge\cdots\wedge X_p$, with $X_i\in
A^1(M,E)$, we set
 \begin{equation*}\rho\circ(X_1\wedge\cdots\wedge X_p)=(\rho\circ
 X_1)\wedge\cdots\wedge(\rho\circ X_p)\,.\end{equation*}
For $p=0$, $ f\in A^0(M,E)=C^\infty(M,\RR)$, we set, as a
convention,
  \begin{equation*}\rho\circ f=f\,.\end{equation*}
Then we extend that correspondence to all elements in $A(M,E)$ by
$C^\infty(M,\RR)$-linearity. The map $P\mapsto \rho\circ P$
obtained in that way is a homomorphism from $A(M,E)$ into
$A(M,TM)$, both for their graded associative algebras structures
(with the exterior products as composition laws) and their graded
Lie algebras structures (with the Schouten-Nijenhuis brackets,
associated to the Lie algebroid structure of $(E,\tau,M,\rho)$ and
to the Lie algebroid structure of the tangent bundle
$(TM,\tau_M,M,\id_{TM})$ as composition laws).
\par\smallskip
In~\ref{Poisson manifolds properties}~(iii), we will see that when the Lie
algebroid under
consideration is the cotangent bundle to a Poisson manifold, the
anchor map has still an additional property: it induces a
cohomology anti-homomorphism.
\end{rmks}

\chapter{Poisson manifolds and Lie algebroids}
\label{Poisson} 
In this final section we will show that there
exist very close links between Poisson manifolds and Lie
algebroids.

\section{Poisson manifolds}
Poisson manifolds were introduced by A.~Lichnerowicz in the very
important paper \cite{Lich}.
Their importance was soon recognized,
and their properties were investigated in depth by A.~Weinstein
\cite{Wein1}. Let us recall briefly their definition and some of
their properties. The reader is referred to \cite{Lich, Wein1,
Vai} for the proofs of these properties.

\begin{defi}
Let $M$ be a smooth manifold. We assume that the space
$C^\infty(M,\RR)$ of smooth functions on $M$ is endowed with a
composition law, denoted by $(f,g)\mapsto\{f,g\}$, for which
$C^\infty(M,\RR)$ is a Lie algebra, which moreover satisfies the
Leibniz-type formula
  \begin{equation*}\{f,gh\}=\{f,g\}h+g\{f,h\}\,.\end{equation*}
We say that the structure defined on $M$ by such a composition law
is a {\it Poisson structure\/}, and that the manifold $M$,
equipped with that structure, is a {\it Poisson manifold\/}.
\end{defi}

The following Proposition is due to A.~Lichnerowicz \cite{Lich}.
Independently, A.~Kirillov \cite{Kir} introduced local Lie
algebras (which include both Poisson manifolds and Jacobi
manifolds, which were introduced too by A.~Lichnerowicz
\cite{Lich2}) and obtained, without using the Schouten-Nijenhuis
bracket, an equivalent result and its generalization for Jacobi
manifolds.

\begin{prop}
On a Poisson manifold $M$, there exists a unique smooth section of the
bundle of bivectors, 
$\Lambda\in A^2(M,TM)$, called the {\it Poisson bivector\/}, which
satisfies
 \begin{equation*}[\Lambda,\Lambda]=0\,,\eqno(*)\end{equation*}
such that for any $f$ and $g\in C^\infty(M,\RR)$,
 \begin{equation*}\{f,g\}=\Lambda(df, dg)\,.\eqno(**)\end{equation*}
The bracket in the left hand side of $(*)$ is the
Schouten-Nijenhuis bracket of multivectors on $M$, for the
canonical Lie algebroid structure of $(TM,\tau_M,M)$ (with
$\id_{TM}$ as anchor).
\par
Conversely, let $\Lambda$ be a smooth section of $A^2(TM,M)$. We
use formula $(**)$ to define a composition law on
$C^\infty(M,\RR)$. The structure defined on $M$ by that
composition law is a Poisson structure if and only if $\Lambda$
satisfies formula $(*)$.
\end{prop}
\par\smallskip
In what follows, we will denote by $(M,\Lambda)$ a manifold $M$
equipped with a Poisson structure whose Poisson bivector is
$\Lambda$.

\section{The Lie algebroid structure on the cotangent bundle of
a Poisson manifold}

The next theorem shows that the cotangent bundle of a Poisson
manifold has a canonical structure of Lie algebroid. That property
was discovered by Dazord and Sondaz \cite{DazSon}.

\begin{theorem}
\label{cotangent of Poisson}
Let $(M,\Lambda)$ be a Poisson manifold. The
cotangent bundle $(T^*M,\pi_M,M)$ has a canonical structure of Lie
algebroid characterized by the following properties:
\par\nobreak\smallskip
\begin{description}
\item{\rm(i)}\quad the bracket $[\eta,\zeta]$ of two sections $\eta$ and
$\zeta$ of $(T^*M,\pi_M,M)$, {\it i.e.}, of two Pfaff forms on
$M$, is given by the formula
  \begin{equation*}\bigl\langle[\eta,\zeta],X\bigr\rangle
  =\bigl\langle\eta,\bigl[\Lambda,\langle \zeta,X\rangle\bigr]\bigr\rangle
  -\bigl\langle\zeta,\bigl[\Lambda,\langle \eta,X\rangle\bigr]\bigr\rangle
  -[\Lambda,X](\eta,\zeta)\,,\end{equation*}
where $X$ is any smooth vector field on $M$; the bracket in the
right hand side of that formula is the Schouten-Nijenhuis bracket
of multivectors on $M$;

\item{\rm(ii)}\quad the anchor is the vector bundle map
$\Lambda^\sharp:T^*M\to TM$ such that, for each $x\in M$, $t$
and $s\in T^*_xM$,

\begin{equation*}\langle \beta,\Lambda^\sharp
\alpha\rangle=\Lambda(\alpha,\beta)\
,.\end{equation*}
\end{description}
\end{theorem}

\begin{proof}
We define a linear endomorphism $\delta_\Lambda$ of $A(M,TM)$ by
setting, for each $P\in A(M,TM)$,
 \begin{equation*}\delta_\Lambda(P)=[\Lambda,P]\,,\end{equation*}
where the bracket in the right hand side is the Schouten-Nijenhuis
bracket of multivectors on $M$, {\it i.e.}, the Schouten-Nijenhuis
bracket for the canonical Lie algebroid structure of
$(TM,\tau_M,M)$ (with $\id_{TM}$ as anchor map). When $P$ is in
$A^p(M,TM)$, $\delta_\Lambda(P)$ is in $A^{p+1}(M,TM)$, therefore
$\delta_\Lambda$ is homogeneous of degree $1$. For each $P\in
A^p(M,TM)$ and $Q\in A^q(M,TM)$, we have
 \begin{equation*}
 \begin{split}
 \delta_\Lambda(P\wedge Q)
 &=[\Lambda,P\wedge Q]\\
 &=[\Lambda,P]\wedge Q+(-1)^pP\wedge[\Lambda,Q]\\
 &=\delta_\Lambda(P)\wedge Q+P\wedge\delta_\Lambda(Q)\,.
 \end{split}
 \end{equation*}
This proves that $\delta_\Lambda$ is a graded derivation of degree
$1$ of the exterior algebra $A(M,TM)$.
\par
Moreover, for each $P\in A^p(M,TM)$ we obtain, by using the graded
Jacobi identity,
 \begin{equation*}
 \begin{split}
 \delta_\Lambda\circ \delta_\Lambda(P)
 &=\bigl[\Lambda,[\Lambda,P]\bigr]\\
 &=(-1)^{p-1}\bigl[\Lambda,[P,\Lambda]\bigr]-\bigl[P,[\Lambda,\Lambda]\bigr]\\
 &=-[\Lambda,[\Lambda,P]\bigr]-\bigl[P,[\Lambda,\Lambda]\bigr]\\
 &=-\delta_\Lambda\circ
 \delta_\Lambda(P)-\bigl[P,[\Lambda,\Lambda]\bigr]\,.
 \end{split}
 \end{equation*}
Therefore
 \begin{equation*}2 \delta_\Lambda\circ
 \delta_\Lambda(P)=-\bigl[P,[\Lambda,\Lambda]\bigr]=0\,,\end{equation*}
since $[\Lambda,\Lambda]=0$. We have proven that the graded
derivation $\delta_\Lambda$, of degree $1$, satisfies
 \begin{equation*}\delta_\Lambda^2=\delta_\Lambda\circ
\delta_\Lambda=0\,.\end{equation*} Now we observe that the tangent
bundle $(TM,\tau_M,M)$ can be considered as the dual bundle of the
cotangent bundle $(T^*M,\pi_M,M)$. Therefore, we may apply
Theorem~\ref{derivation on dual bundle}, which shows that there exists on
$(T^*M,\pi_M,M)$ a
Lie algebroid structure for which $\delta_M$ is the associated
derivation on the space $\Omega(M,T^*M)=A(M,TM)$ (with the
notations defined 
in~\ref{The exterior algebra of sections}). 
That theorem also shows that the
bracket of two smooth sections of $(T^*M,\pi_M,M)$, {\it i.e.}, of
two Pfaff forms $\eta$ and $\zeta$ on $M$, is given by the
formula, where $X$ is any smooth vector field on $M$,
 \begin{equation*}\bigl\langle[\eta,\zeta],X\bigr\rangle
 =\bigl\langle\eta,\bigl[\Lambda,\langle \zeta,X\rangle\bigr]\bigr\rangle
 -\bigl\langle\zeta,\bigl[\Lambda,\langle \eta,X\rangle\bigr]\bigr\rangle
 -[\Lambda,X](\eta,\zeta)\,.\end{equation*}
The anchor map $\rho$ is such that, for each $\eta\in
\Omega^1(M,TM)$ and each $f\in C^\infty(M,\RR)$,
 \begin{equation*}i(\rho\circ\eta)\,df
 =\bigl\langle\eta,[\Lambda,f]\bigr\rangle\,
.\end{equation*} 
The bracket which appears in the right hand sides
of these two formulae is the Schouten-Nijenhuis bracket of
multivectors on $M$. By using Theorem~\ref{Schouten bracket}, we see that
 \begin{equation*}[\Lambda,f]=-\Lambda^\sharp(df)\,.\end{equation*}
Therefore,
 \begin{equation*}\langle df,\rho\circ\eta\rangle=i(\rho\circ\eta)\,df
 =\bigl\langle\eta,-\Lambda^\sharp(df)\bigr\rangle=\bigl\langle
 df,\Lambda^\sharp(\eta)\bigr\rangle\,.\end{equation*}
So we have $\rho=\Lambda^\sharp$.
\end{proof}

\begin{rmks}
\label{Poisson manifolds properties}
Let $(M,\Lambda)$ be a Poisson manifold.
\par\nobreak\smallskip\noindent
{\rm(i)}\quad{\it The bracket of forms of any degrees on
$M$\/}.\quad Since, by Theorem~\ref{cotangent of Poisson},
$(T^*M,\pi_M,M,\Lambda^\sharp)$ is a Lie algebroid, we can define
a composition law in the space $A(M,T^*M)=\Omega(M,\RR)$ of smooth
differential forms of all degrees on $M$: the Schouten-Nijenhuis
bracket for the Lie algebroid structure of $(T^*M,\pi_M,M)$, with
$\Lambda^\sharp$ as anchor. With that composition law, denoted by
$(\eta,\zeta)\mapsto[\eta,\zeta]$, $\Omega(M,\RR)$ is a graded Lie
algebra. Observe that a form $\eta\in\Omega^p(M,\RR)$, of degree
$p$ for the graded associative algebra structure whose composition
law is the exterior product, has degree $p-1$ for the graded Lie
algebra structure.
\par
The bracket of differential forms on a Poisson manifold was first
discovered for Pfaff forms by Magri and Morosi \cite{MaMo}. It is
related to the Poisson bracket of functions by the formula
 \begin{equation*}[df,dg]=d\{f,g\}\,,\quad\hbox{with}\quad f\ \hbox{and}\ g\in
 C^\infty(M,\RR)\,.\end{equation*}
That bracket was extended to forms of all degrees by Koszul
\cite{Ko2}, and rediscovered, with the Lie algebroid structure of
$T^*M$, by Dazord and Sondaz \cite{DazSon}.
\par\nobreak\smallskip\noindent
{\rm(ii)}\quad{\it The Lichnerowicz-Poisson cohomology\/}.\quad
The derivation $\delta_\Lambda$,
 \begin{equation*}P\mapsto\delta_\Lambda(P)=[\Lambda,P]\,,\quad P\in
A(M,TM)\,,\end{equation*} used in the proof 
of~\ref{cotangent of Poisson}, was
first introduced by A.~Lichnerowicz \cite{Lich}, who observed that
it may be used to define a cohomology with elements in $A(M,TM)$
as cochains. He began the study of that cohomology, often called
the Poisson cohomology (but which should be called the
Lichnerowicz-Poisson cohomology). The study of that cohomology was
carried on by Vaisman \cite{Vai}, Huebschmann \cite{Huebs}, Xu
\cite{Xu} and many other authors.
\par\nobreak\smallskip\noindent
{\rm(iii)}\quad{\it The map $\Lambda^\sharp$ as a cohomology
anti-homomorphism\/}.\quad In~\ref{remarks about degrees}~(ii), 
we have seen that the
anchor map $\rho$ of a Lie algebroid $(E,\tau,M,\rho)$ yields a map
$P\mapsto \rho\circ P$ from $A(M,E)$ into $A(M,TM)$, which is both
a homomorphism of graded associative algebras (the composition
laws being the exterior products) and a homomorphism of graded Lie
algebras (the composition laws being the Schouten brackets). When
applied to the Lie algebroid $(T^*M,\pi_M,M,\Lambda^\sharp)$, that
property shows that the map $\eta\mapsto\Lambda^\sharp\circ\eta$
is a homomorphism from the space of differential forms
$\Omega(M,\RR)$ into the space of multivectors $A(M,\RR)$, both
for their structures of graded associative algebras and their
structures of graded Lie algebras. As observed by A.~Lichnerowicz
\cite{Lich}, that map exchanges the exterior derivation $d$ of
differential forms and the derivation $\delta_\Lambda$ of
multivectors (with a sign change, under our sign conventions), in
the following sense: for any $\eta\in \Omega^p(M,\RR)$, we have
 \begin{equation*}\Lambda^\sharp(d\eta)=
 -\delta_\Lambda\bigl(\Lambda^\sharp(\eta)\bigr)
 =-\bigl[\Lambda,\Lambda^\sharp(\eta)\bigr]\,.\end{equation*}
That property is an easy consequence of the formula, valid for any
smooth function $f\in C^\infty(M,\RR)$, which can be derived from
Theorem~\ref{Schouten bracket},
 \begin{equation*}\Lambda^\sharp(df)=-[\Lambda,f]\,.\end{equation*}
The map $\Lambda^\sharp$ therefore induces an anti-homomorphism
from the Lichnerowicz-Poisson cohomology of the Poisson manifold
$(M,\Lambda)$, into its De~Rham cohomology.
\par\nobreak\smallskip\noindent
{\rm(iv)}\quad{\it Lie bialgebroids\/}.\quad Given a Poisson
manifold $(M,\Lambda)$, we have Lie algebroid structures both on
the tangent bundle $(TM,\tau_M,M)$ and on the cotangent bundle
$(T^*M,\pi_M,M)$, with $\id_{TM}:TM\to TM$ and
$\Lambda^\sharp:T^*M\to TM$ as their respective anchor maps.
Moreover, these two Lie algebroid structures are compatible in the
following sense: the derivation $\delta_\Lambda:P\mapsto [\Lambda,
P]$ of the graded associative algebra $A(M,TM)$ (the composition
law being the exterior product) determined by the Lie algebroid
structure of $(T^*M,\pi_M,M)$ is also a derivation for the graded
Lie algebra structure of $A(M,E)$ (the composition law being now
the Schouten-Nijenhuis bracket). We have indeed, as an easy
consequence of the graded Jacobi identity, for $P\in A^p(M,TM)$
and $Q\in A^q(M,TM)$,
 \begin{equation*}
 \begin{split}
 \delta_\Lambda\bigl([P,Q]\bigr)
 &=\bigl[\Lambda,[P,Q]\bigr)
 =\bigl[[\Lambda,P],Q\bigr]+(-1)^{p-1}\bigl[P,[\Lambda,Q]]\\
 &=[\delta_\Lambda P,Q]+(-1)^{p-1}[P,\delta_\Lambda Q]\,.
 \end{split}
 \end{equation*}
When two Lie algebroid structures on two vector bundles in duality
satisfy such a compatibility condition, it is said that that pair
of Lie algebroids is a {\it Lie bialgebroid\/}. The very important
notion of a Lie bialgebroid is dut to K.~Mackenzie and P.~Xu
\cite{MackXu}. Its study was developed by Y.~Kosmann-Schwarzbach
\cite{Kosm} and her student \cite{BanKo} and many other authors.
D.~Iglesias and J.C.~Marrero have introduced a
generalization of that notion in relation with Jacobi manifolds
\cite{IgMar}.
\end{rmks}

\section{The Poisson structure on the dual bundle of a Lie
algebroid} 
We will now prove that there is a $1$--$1$ correspondence between
Lie algebroid structures on a vector bundle $(E,\tau,M)$ and
homogeneous Poisson structures on the total space of the dual bundle
$(E^*,\pi,M)$. This will allow us to recover
well known results (Remarks~\ref{known examples of Poisson structures}).
\par\smallskip
We will use the following definition.

\begin{defi}
 \label{vertical functions}
Let $(E,\tau,M)$ be a vector bundle and $(E^*,\pi,M)$ its
dual bundle.
To each smooth section $X\in A^1(M,E)$, we associate
the smooth function $\Phi_X$ defined on $E^*$ by
 \begin{equation*}
 \Phi_X(\xi)=\bigl\langle\xi,X\circ\pi(\xi)\bigr\rangle\,,\quad \xi\in E^*\,.
 \end{equation*}
We will say that $\Phi_X$ is the \emph{vertical function} on $E^*$ associated
to the smooth section $X$. 
\end{defi}

\begin{lemma}
\label{lemma poisson on dual bundle}
Let $(E,\tau,M)$ be a vector bundle and $(E^*,\pi,M)$ its
dual bundle. 
\par\smallskip\noindent
{\rm 1.\quad} If, for some smooth section $X\in A^1(M,E)$, some smooth
function $f\in C^\infty(M,\RR)$ and some $\xi\in E^*$,
$d(\Phi_X+f\circ\pi)(\xi)=0$, where $\Phi_X$ is the vertical function
associated to $X$ (Definition~\ref{vertical functions}), then
$X\bigl(\pi(\xi)\bigr)=0$.
\par\smallskip\noindent
{\rm 2.\quad}  For each $\xi\in E^*$ and each 
$\eta\in T^*_\xi E^*$, there exists a smooth section $X\in A^1(M, E)$ and a
smooth function $f\in C^\infty(M,\RR)$ such that 
$d(\Phi_X+f\circ\pi)(\xi)=\eta$.
\end{lemma}
\begin{proof} These properties being local we may work in an open subset $U$ of
$M$ on which there exists a system of local coordinates $(x^1,\ldots,x^n)$ and
smooth sections $(s_1,\ldots,s_k)$ of $\tau$, such that for each $x\in U$,
$(s_1(x),\ldots,s_k(x))$ is a basis of $E_x$. 
A smooth section $X$ of $\tau$ defined on $U$ can be written
 $$X=\sum_{r=1}^k X^rs_r\,,$$
where the $X^r$ are smooth functions on $U$. We will denote by the
same letters $X^r$ the expression of these functions in local
coordinates $(x^1,\ldots,x^n)$. Similarly we will denote by $f$
both a smooth function in $C^\infty(M,\RR)$ and its expression in local
coordinates. The vertical function, defined on $\pi^{-1}(U)$, which
corresponds to $X$ is
 $$\Phi_X(\xi)= \sum_{r=1}^k\xi_rX^r\bigl(\pi(\xi)\bigr)\,,\quad
 \xi\in\pi^{-1}(U)\,,\quad\hbox{where}\ \xi_r=\bigl\langle\xi,
 s_r\bigl(\pi(\xi)\bigr)\bigr\rangle\,.$$
On $\pi^{-1}(U)$, $(x^1,\ldots,x^n,\allowbreak \xi_1,\ldots,\xi_k)$ is a smooth
system of local coordinates, in which
 $$d(\Phi_X+f\circ\pi)(\xi)=\sum_{r=1}^k X^r(x^1,\ldots,x^n)\,d\xi_r
 +\sum_{j=1}^n\left(\sum_{r=1}^k \xi_r
 \frac{\partial X^r(x^1,\ldots,x^n)}{\partial x^j}
 +\frac{\partial f(x^1,\ldots,x^n)}{\partial x^j}\right)\,dx^j\,.$$
This result shows that if $d(\Phi_X+f\circ\pi)(\xi)=0$, then
$X\bigl(\pi(\xi)\bigr)=0$.
\par\smallskip 
Let $\xi\in E^*$ and $\eta\in T^*_\xi E^*$ be given. The above formula
shows that if $\xi\neq 0$, we can take $f=0$ and choose $X$ such that
$d\Phi_X(\xi)=\eta$. If $\xi=0$, we can take $X=0$,
and $f$ such that $d(f\circ\pi)(\xi)=\eta$.
\end{proof}

\begin{defi}
\label{homogeneous Poisson structure}
Let $(E,\tau,M)$ be a vector bundle and $(E^*,\pi,M)$ its
dual bundle. A Poisson structure on $E^*$ is said to be \emph{homogeneous}
if the Poisson bracket of two vertical functions 
(Definition~\ref{vertical functions}) is vertical.
\end{defi}

\begin{prop}
\label{properties of homogeneous Poisson structures}
Let $(E,\tau,M)$ be a vector bundle, $(E^*,\pi,M)$ its
dual bundle, and $\Lambda$ be a Poisson structure on $E^*$. The following
properties are equivalent.

\par\smallskip\noindent
{\rm 1.\quad} There exists a dense subset $U$ of $E^*$ and a subset $\cal F$ of
the set of vertical functions on $E^*$ whose differentials $df(\xi)$,
$f\in{\cal F}$, span the cotangent space $T^*_\xi E^*$, for all $\xi\in U$,
such that the Poisson bracket of two functions in $\cal F$ is vertical.

\par\smallskip\noindent
{\rm 2.\quad} Let $Z_{E^*}$ be the vector field on $E^*$ whose
flow generates homotheties in the fibres. 
We recall that its value at $\xi\in E^*$ is
$\displaystyle Z_{E^*}(\xi)=\frac{d\bigl(\exp(t)\xi\bigr)}{dt}\bigm|_{t=0}$.
The Poisson structure on $E^*$  satisfies
 \begin{equation*}
 [Z_{E^*},\Lambda]=-\Lambda\,.
 \end{equation*}

\par\smallskip\noindent
{\rm 3.\quad} The Poisson structure $\Lambda$ is homogeneous.
\end{prop}

\begin{proof}
The reduced flow of the
vector field $Z_{E^*}$ is the one parameter
group of homotheties in the fibres
$(t,\xi)\mapsto H_t(\xi)=\exp(t)\xi$, with $t\in\RR$, $\xi\in E^*$.
For any smooth section $X\in A^1(M,E)$ and any $t\in\RR$, we have
 \begin{equation*}
 (H_t^*\Phi_X)(\xi)=\Phi_X\circ
H_t(\xi)=\Phi_X\bigl(\exp(t)\xi\bigr)=\exp(t)\Phi_X(\xi)\,, 
 \end{equation*}
therefore
 \begin{equation*}
 H^*_t\Phi_X=\exp(t)\Phi_X\,,\quad
 {\cal L}(Z_{E^*})\Phi_X=
 \frac{dH_t^*\Phi_X}{dt}\big|_{t=0}=\Phi_X\,. 
 \end{equation*}
Let us assume that~1 is true. Let $X$ and $Y\in A^1(M,E)$ be such that
$\Phi_X$ and $\Phi_Y$ are in the subset $\cal F$. Then $\{\Phi_X,\Phi_Y\}$ is
vertical, so for all $t\in\RR$ we have
 \begin{equation*}
 H_t^*\bigl(\Lambda(d\Phi_X,d\Phi_Y)\bigr)
 =\{\Phi_X,\Phi_Y\}\circ H_t=\exp(t)\{\Phi_X,\Phi_Y\}\,.
 \end{equation*}
But we may also write
 \begin{equation*}
 H_t^*\bigl(\Lambda(d\Phi_X,d\Phi_Y)\bigr)
 =(H_t^*\Lambda)\bigl(H_t^*d\Phi_X,H_t^*d\Phi_Y\bigr)
 =\exp(2t)(H_t^*\Lambda)\bigl(d\Phi_X,d\Phi_Y\bigr)\,.
 \end{equation*}
Since for each $\xi\in U$, the differentials at $\xi$ of functions in $\cal F$
generate $T_\xi^*E^*$, this result proves that in $U$
 \begin{equation*}
 H_t^*(\Lambda)=\exp(-t)\Lambda\,.
 \end{equation*}
Since $U$ is dense in $E^*$ that equality holds everywhere on $E^*$,
therefore
\begin{equation*}
 [Z_{E^*},\Lambda]={\cal L}(Z_{E^*})\Lambda=
\frac{dH_t^*\Lambda}{dt}\big|_{t=0}=-\Lambda\,.
\end{equation*}
We have proven that~1 implies~2. Let us now assume that~2 is true. 
For all $X$ and $Y\in A^1(M,E)$,
 \begin{equation*}
\begin{split}
 {\cal L}(Z_{E^*})\bigl(\{\Phi_X,\Phi_Y\}\bigr)
 &={\cal L}(Z_{E^*})\bigl(\Lambda(d\Phi_X,d\Phi_Y)\bigr)\\
 &=\bigl({\cal L}(Z_{E^*})\Lambda\bigr)(d\Phi_X,d\Phi_Y)
 +\Lambda\bigl({\cal L}(Z_{E^*})\Phi_X,\Phi_Y) 
 +\Lambda\bigl(\Phi_X,{\cal L}(Z_{E^*})\Phi_Y\bigr)\\
 &=\Lambda(\Phi_X,\Phi_Y)=\{\Phi_X,\Phi_Y\}\,.
 \end{split}
 \end{equation*}
Since $\{\Phi_X,\Phi_Y\}$ is smooth on $E^*$,
including on the zero section, this function is linear on each fibre of $E^*$,
in other words it is  vertical, and we have proven that~2 implies~3.

Finally,~3 inplies of course~1, and our proof is complete.
\end{proof}

\begin{theorem}
\label{theorem poisson on dual bundle}
Let $(E,\tau,M)$ be a vector bundle and $(E^*,\pi,M)$ its
dual bundle. There is a $1$--$1$ correspondence between Lie algebroid
structures on $(E,\tau,M)$ and homogeneous Poisson structures on $E^*$
(Definition~\ref{homogeneous Poisson structure})
such that, for each pair $(X,Y)$ of smooth sections of $\tau$, $\Phi_X$ and
$\Phi_Y$ being the corresponding vertical functions on $E^*$
(Definition~\ref{vertical functions}),
 \begin{equation*}
 \{\Phi_X,\Phi_Y\}=\Phi_{\{X,Y\}}\,,
 \end{equation*}
the bracket in the left hand side being the Poisson bracket of
functions on $E^*$, and the bracket in the right hand side the
bracket of sections for the corresponding Lie algebroid structure on
$(E,\tau,M)$. 
\end{theorem}

\begin{proof}
First let $\Lambda$ be an homogeneous Poisson structure
on $E^*$. Let $(X,Y)$ be a pair of smooth sections of $\tau$,
$\Phi_X$ and
$\Phi_Y$ the corresponding vertical functions on $E^*$. Since $\Lambda$ is
homogeneous, there exists a unique smooth section of $\tau$ whose corresponding
vertical function on $E^*$ is $\{\Phi_X,\Phi_Y\}$. We define $\{X,Y\}$ as being
that section. So we have a composition law on the space $A^1(M,E)$ of smooth
sections of $\tau$, which is bilinear and satisfies the Jacobi identity, and
therefore is a Lie algebra bracket. Now let $f$ be a smooth function on $M$.
Then
\begin{equation*}
 \begin{split}
 \{\Phi_X,\Phi_{fY}\}&=\{\Phi_X,(f\circ\pi)\Phi_Y\}=(f\circ\pi)
 \{\Phi_X,\Phi_Y\}+\{\Phi_X,f\circ\pi\}\Phi_Y\\
 &=(f\circ\pi)\Phi_{\{X,Y\}}+
  \{\Phi_X,f\circ\pi\}\Phi_Y
 =(f\circ\pi)\Phi_{\{X,Y\}}
 +\Bigl(i\bigl(\Lambda^\sharp(d\Phi_X)\bigr)d(f\circ\pi)\Bigr)s\Phi_Y\,.
 \end{split}
\end{equation*}
The term $(f\circ\pi)\Phi_{\{X,Y\}}$ is the vertical function which
corresponds to the smooth section $f\{X,Y\}$. Therefore the other term of
the last side,
$\Bigl(i\bigl(\Lambda^\sharp(d\Phi_X)\bigr)d(f\circ\pi)\Bigr)\Phi_Y$, must be a
vertical function. But $\Phi_Y$ is vertical, so
$\Bigl(i\bigl(\Lambda^\sharp(d\Phi_X)\bigr)d(f\circ\pi)\Bigr)\Phi_Y$
is vertical for all $Y\in A^1(M,E)$ if and only if 
the function $i\bigl(\Lambda^\sharp(d\Phi_X)\bigr)d(f\circ\pi)$ is constant on
each fibre $\pi^{-1}(x)$, $x\in M$.
This happens
for any function $f\in C^\infty(M,\RR)$ if and only if
for each $x\in M$,
$T_\xi\pi\bigl(\Lambda^\sharp(d\Phi_X)\bigr)$ does not depend on
$\xi\in\pi^{-1}(x)$. In other words, for any $X\in A^1(M,E)$ the vector field
$\Lambda^\sharp(d\Phi_X)$ must be projectable by $\pi$ on $M$. We may take
$\xi=0$ (the origin of the fibre $E^*_x$), use the formula for $d\Phi_X$
given in the proof of Lemma~\ref{lemma poisson on dual bundle} and we obtain for
that projection in local coordinates
  $$\sum_{r=1}^kX^r(x^1,\ldots,x^n)T\pi\bigl(\Lambda^\sharp(d\xi_r)\bigr)\,.
  $$ 
The value of that vector field at a point $x\in M$ only depends of $X(x)$, and
that dependence is linear. So there exists a smooth vector bundle map $
\rho:E\to TM$ with all the properties of an anchor map.
The Lie algebra structure we have defined on $A^1(M,E)$ is a Lie algebroid
bracket.
\par\smallskip
Conversely, let us assume that we have on $(E,\pi,M)$ a Lie algebroid structure
with anchor $\rho$. We must prove that there exists a Poisson structure on $E^*$
such that for each pair $(X,Y)$ of smooth sections of $\tau$, $\Phi_X$ and
$\Phi_Y$ being the corresponding vertical functions on $E^*$
(Definition~\ref{vertical functions}),
$\{\Phi_X,\Phi_Y\}=\Phi_{\{X,Y\}}$.
More generally, let $(X,Y)$ be a pair of smooth sections of $\tau$, $f$ and
$g$ two smooth functions. Let us write that
$\{\Phi_{fX},\Phi_{gY}\}=\bigl\{f\circ\pi)\Phi_X,(g\circ\pi)\Phi_Y\bigr\}
=\Phi_{\{ fX , gY\} } $. We use the property 
of the Lie
algebroid bracket
 \begin{equation*}
 \{fX,gY\}=fg\{X,Y\}+\bigl(f{\cal L}(\rho\circ X)g\bigr)Y
                    -\bigl(g{\cal L}(\rho\circ Y)f\bigr)X\,,
 \end{equation*}
which implies
 \begin{equation*}
 \Phi_{\{fX,gY\}}=(fg\circ\pi)\Phi_{\{X,Y\}}
                  +\bigl(f{\cal L}(\rho\circ X)\circ\pi\bigr)\,\Phi_Y
                  -\bigl(g{\cal L}(\rho\circ
                  Y)f\circ\pi\bigr)\,\Phi_X\,.
 \end{equation*}
This calculation shows that if such a Poisson structure on $E^*$ exists, it
must be such that
\begin{equation*}
 \{\Phi_X,g\circ\pi\}=\bigl({\cal L}(\rho\circ
 X)g\bigr)\circ\pi\,,\qquad
 \{f\circ\pi,g\circ\pi\}=0\,.
 \end{equation*}
For each $\xi\in E^*$, $\eta$ and $\zeta\in
T^*_\xi E^*$, point 2 of Lemma~\ref{lemma poisson on dual bundle} shows that
there exists a (non unique) pair $(X,Y)$ of sections of $\tau$ and a (non
unique) pair $(f,g)$ of smooth functions on $M$ such that
 $\eta=d(\Phi_X+f\circ\pi)(\xi)$,
 $\zeta=d(\Phi_Y+g\circ\pi)(\xi)$.
Our Poisson bivector $\Lambda$ is therefore
 \begin{equation*}
 \Lambda(\xi)(\eta,\zeta)
 =\{\Phi_X+f\circ\pi,\Phi_Y+g\circ\pi\}(\xi)\,.
 \end{equation*}
This proves that if such a Poisson structure exists, it is unique.
By point 1 of 
Lemma~\ref{lemma poisson on dual bundle}, the
right hand side of the above formula depends only on $\eta$ and
$\zeta$, not on the particular choices we have made for
$(X,f)$ and $(Y,g)$. Moreover, it is smooth, bilinear and skew-symetric with
respect to the pair
$\bigl((X,f),(Y,g)\bigr)$, so $\Lambda$ is a smooth bivector.
\par\smallskip
When restricted to vertical functions on $E^*$, the bracket defined by
$\Lambda$ satisfies the Jacobi
identity. Therefore, for
each $\xi\in E^*\,\backslash\{0\}$, $\eta$, $\zeta$ and 
$\theta\in T^*_\xi E^*$
which are the differentials, at $\xi$, of vertical functions, the Schouten
bracket $[\Lambda,\Lambda]$
satisfies $[\Lambda,\Lambda](\xi)(\eta,\zeta,\theta)=0$. Point 2
of~Lemma~\ref{lemma poisson on dual bundle} proves that
$[\Lambda,\Lambda]$ vanishes identically on $E^*\,\backslash\{0\}$. By
continuity, it vanishes everywhere on $E^*$. So $\Lambda$ is a Poisson
structure on $E^*$ with all the stated properties.
\end{proof}

\begin{prop}
\label{poisson properties of dual bundle}
Let $(E,\tau,M,\rho)$ be a Lie algebroid and $(E^*,\pi,M)$ its
dual bundle. The Poisson structure on $E^*$ defined in
Theorem~\ref{theorem poisson on dual bundle} has the following properties:
\par\smallskip\noindent
{\rm 1.}\quad For any $X\in A^1(M,E)$, $f$ and $g\in
C^\infty(M,\RR)$,
 \begin{equation*}
 \{\Phi_X,g\circ\pi\}=\bigl({\cal L}(\rho\circ
 X)g\bigr)\circ\pi\,,\qquad
 \{f\circ\pi,g\circ\pi\}=0\,,
 \end{equation*}
where $\Phi_X$ is the function on $M$ associated to the
section $X$ as indicated in 
Theorem~\ref{theorem poisson on dual bundle}.
\par\smallskip\noindent
{\rm 2.}\quad The transpose $\,{}^t\!\rho:T^*M\to E^*$ of the
anchor map $\rho:E\to TM$ is a Poisson map (the cotangent bundle
being endowed with the Poisson structure associated to its
canonical symplectic structure).
\end{prop}

\begin{proof}
We have proven Properties~1 in the proof of 
Theorem~\ref{theorem poisson on dual bundle}. In
order to prove Property~2, we must prove that for all pairs
$(h_1,h_2)$ of smooth functions on $E^*$,
 \begin{equation*}
 \{h_1\circ\,{}^t\!\rho,h_2\circ\,{}^t\!\rho\}
 =\{h_1,h_2\}\circ\,{}^t\!\rho\,,
 \end{equation*}
the bracket in the left hand side being the Poisson bracket of
functions on $T^*M$, and the bracket in the right hand side the
Poisson bracket of functions on $E^*$. It is enough to check that
property when $h_1$ and $h_2$ are of the type $\Phi_X$, where
$X\in A^1(M,E)$, or of the type $f\circ\pi$, with $f\in
C^\infty(M,\RR)$, since the differentials of functions of these
two types generate $T^*E^*$. For $h_1=\Phi_X$ and $h_2=\Phi_Y$,
with $X$ and $Y\in A^1(M,E)$, and $\zeta\in T^*M$, we have
 \begin{equation*}
 \begin{split}
 \{\Phi_X,\Phi_Y\}\circ\,{}^t\!\rho(\zeta)
 &=\Phi_{\{X,Y\}}\circ\,{}^t\!\rho(\zeta)\\
 &=\bigl\langle\,{}^t\!\rho(\zeta),\{X,Y\}
   \circ\pi\circ\,{}^t\!\rho(\zeta)\bigr\rangle\\
 &=\bigl\langle\zeta,\rho\circ\{X,Y\}\circ\pi_M(\zeta\bigr\rangle\\
 &=\bigl\langle\zeta,[\rho\circ X,\rho\circ Y]\circ\pi_M(\zeta\bigr\rangle\,,
 \end{split}
 \end{equation*}
since the canonical projection $\pi_M:T^*M\to M$ satisfies
$\pi\circ\,{}^t\!\rho=\pi_M$. But let us recall a well known
property of the Poisson bracket of functions on $T^*M$
(\cite{LibMa}, exercise 17.5 page 182). To any vector field
$\widehat X$ on $M$, we associate the function $\Psi_{\widehat X}$
on $T^*M$ by setting, for each $\zeta\in T^*M$,
 \begin{equation*}
 \Psi_{\widehat X}(\zeta)=\bigl\langle\zeta,\widehat
 X\circ\pi_M(\zeta)\bigr\rangle\,.
 \end{equation*}
Then, for any pair $(\widehat X,\widehat Y)$ of vector fields on
$M$,
 \begin{equation*}
 \{\Psi_{\widehat X},\Psi_{\widehat Y}\}=\Psi_{[\widehat X,\widehat Y]}\,.
 \end{equation*}
By using $\pi_M=\pi\circ\,{}^t\!\rho$, we easily see that for
each $X\in A^1(M,E)$,
 \begin{equation*}
 \Psi_{\rho\circ X}=\Phi_X\circ\,{}^t\!\rho\,.
 \end{equation*}
Returning to our pair of sections $X$ and $Y\in A^1(M,E)$, we see
that
 \begin{equation*}
 \{\Phi_X,\Phi_Y\}\circ\,{}^t\!\rho(\zeta)=\Psi_{[\rho\circ
 X,\rho\circ Y]}(\zeta)=\{\Psi_{\rho\circ X},\Psi_{\rho\circ
 Y}\}(\zeta)
 =\{\Phi_X\circ\,{}^t\!\rho,\Phi_Y\circ\,{}^t\!\rho\}(\zeta)\,.
 \end{equation*}
Now for $h_1=\Phi_X$ and $h_2=f\circ\pi$, with $X\in A^1(M,E)$
and $f\in C^\infty(M,\RR)$, we have
 \begin{equation*}
 \begin{split}
 \{\Phi_X\circ\,{}^t\!\rho,f\circ\pi\circ\,{}^t\!\rho\}
 &=\{\Psi_{\rho\circ X},f\circ\pi_M\}
 ={\cal L}(\rho\circ X)f\circ\pi_M\\
 &={\cal L}(\rho\circ X)f\circ\pi\circ\,{}^t\!\rho
 =\{\Phi_X,f\circ\pi\}\circ\,{}^t\!\rho\,.
 \end{split}
 \end{equation*}
Similarly, for $h_1=f\circ\pi$ and $h_2=g\circ\pi$, we have
 \begin{equation*}
 \{f\circ\pi,g\circ\pi\}\circ\,{}^t\!\rho=0=\{f\circ\pi_M,g\circ\pi_M\}
 =\{f\circ\pi\circ\,{}^t\!\rho,g\circ\pi\circ\,{}^t\!\rho\}\,.
 \end{equation*}
Property 2 is proven, and our proof is complete.
\end{proof}

\begin{rmks}
\label{known examples of Poisson structures}
 \par\nobreak\smallskip\noindent
{\rm(i)}\quad{\it The symplectic structure of a cotangent bundle.}\quad
Let us take as Lie algebroid the tangent bundle $TM,\tau_M,M)$, with $\id_{TM}$
as anchor. Its dual bundle is the cotangent bundle $(T^*M,\pi_M,M)$. The
transpose of the anchor map being $\id_{T^*M}$, 
Proposition~\ref{poisson properties of dual bundle} shows that the Poisson
structure on $T^*M$ given by Theorem~\ref{theorem poisson on dual bundle} is
the structure associated to its canonical symplectic $2$-form.
 \par\nobreak\smallskip\noindent
{\rm(i)}\quad{\it The symplectic structure on the dual of a Lie algebra.}\quad
Now we take as Lie algebroid a finite dimensional Lie algebra $\cal G$. The
Poisson structure on its dual vector space ${\cal G}^*$ given 
by Theorem~\ref{theorem poisson on dual bundle} is the well
known
Kirillov-Kostant-Souriau Poisson structure \cite{Kir, Kos, Sou}.
\end{rmks}

\section{Tangent lifts}
\label{tangent lifts}
G.~Sanchez~de~Alvarez~\cite{Sanchez} discovered the lift of a Poisson structure
on a manifold $P$ to the tangent bundle $TP$. We show below
that its existence and properties can be easily deduced
from Theorems~\ref{cotangent of Poisson} 
and~\ref{theorem poisson on dual bundle}.
The reader will find many other properties of tangent and cotangent lifts of
Poisson an Lie algebroid structures in~\cite{GraUrb1, GraUrb2}.

\begin{theorem}
 \label{tangent lift of Poisson}
Let $(P,\Lambda)$ be a Poisson manifold. There exists on
its tangent bundle $TP$ a Poisson structure, determined by that of $P$
and called its \emph{tangent lift}. It is such that, if $f$ and 
$g\in C^\infty(P,\RR)$ are two smooth functions on $P$,
 $$\{df, dg\}_{TP}=d\{f,g\}_P\,,\quad \{f\circ\tau_P,g\circ
\tau_P\}_{TP}=0\,,\quad
\{df,g\circ\tau_P\}_{TP}=\{f,g\}_P\circ\tau_P\,.$$
In these formulae we have denoted by $\{\ ,\ \}_{P}$ and $\{\ ,\, \}_{TP}$ the
Poisson brackets of functions on $P$ and $TP$, respectively, and we have
considered $df$, $dg$ and $d\{f,g\}_P$ as vertical functions on $TP$.
\end{theorem}

\begin{proof}
The cotangent bundle $(T^*P, \pi_P,P)$ has a Lie algebroid structure, with
$\Lambda^\sharp: T^*P\to TP$ as anchor 
(Theorem~\ref{cotangent of Poisson}). Its dual is the tangent bundle 
$(TP,\tau_P,P)$, and by Theorem\ref{theorem poisson on dual bundle}, there
exists on its total space $TP$ a Poisson structure such that, for each pair
$(\eta,\zeta)$ of sections of $\pi_P$,
 $$\{\Phi_\eta,\Phi_\zeta\}_{TP}=\Phi_{[\eta,\zeta]}\,.$$
We have denoted by $\Phi_\eta$ and $\Phi_\zeta$ the vertical functions on $TP$
associated to the sections $\eta$ and $\zeta$ of $\pi_P$ 
(\ref{vertical functions}) and denoted by $[\eta,\zeta]$ the bracket of the
Pfaff forms $\eta$ and $\zeta$ on the Poisson manifold $(P,\Lambda)$
(\ref{cotangent of Poisson}). When $\eta=df$ and $\zeta=dg$, we have
$[df,dg]=d\{f,g\}_P$. The properties of  the Poisson bracket on $TP$ follow from
Proposition~\ref{poisson properties of dual bundle}.
\end{proof}

\begin{ex}
Let us assume that $P$ is of even dimension $2m$ and that its Poisson
structure is associated to a symplectic $2$-form $\omega_P$. In local Darboux
coordinates $(x^1\ldots,x^m,\allowbreak y_1,\ldots,y_m)$ we have
  $$\omega_P=\sum_{i=1}^m dy_i\wedge dx^i\,,\quad 
    \Lambda_P=\sum_{i=1}^m \frac{\partial}{\partial
    y_i}\wedge\frac{\partial}{\partial x^i}\,.
 $$
Let $(x^1,\ldots,x^m,y_1,\ldots,y_m,\dot x^1,\ldots,\dot x^m,\allowbreak
\dot y_ 1,\ldots,\dot y_m)$ be the local coordinates on $TP$ naturally
associated to the local coordinates $(x^1,\ldots,x^m,y_1,\ldots,y_m)$ on $P$. We
easily see that the lift to $TP$ of the Poisson structure on $P$ is associated
to a symplectic structure $\omega_{TP}$, and that the expressions of
$\Lambda_{TP}$ and $\omega_{TP}$ in local coordinates are
\begin{equation*}
 \omega_{TP}=\sum_{i=1}^m (d\dot y_i\wedge dx^i+ dy_i\wedge d\dot x^i)\,,\quad 
    \Lambda_{TP}=\sum_{i=1}^m \left(\frac{\partial}{\partial
    \dot y_i}\wedge\frac{\partial}{\partial x^i}+
   \frac{\partial}{\partial
    y_i}\wedge\frac{\partial}{\partial \dot x^i}\right)\,.
\end{equation*}
The symplectic form $\omega_{TP}$ was defined and used by
W.M.~Tulczyjew~\cite{Tul1, Tul2}, mainly when $P$ is  a cotangent bundle. 
It can be defined by several other
methods. For example, since $\Lambda_P$ is associated to a symplectic
structure, $\Lambda^\sharp_P$ is a fibre bundle isomorphim from $T^*P$ onto
$TP$. There is on $T^*P$ a canonical symplectic form $\omega_{T^*P}$ (the
exterior differential of its Liouville $1$-form). With our sign conventions,
$\omega_{TP}=-\bigl((\Lambda^\sharp_P) ^{-1}\bigr)^*(\omega_{T^*P})$. The $-$
sign is in agreement with Point 2 of \ref{poisson properties of dual bundle},
since the transpose on $\Lambda_P^\sharp$ is $-\Lambda_P^\sharp$.
\end{ex}

For Lie algebroid structures, there is an even richer
notion of lift: the next proposition
shows that a Lie algebroid structure on a vector bundle $(E,\tau,M)$ gives rise
to Lie algebroid structures on two vector bundles: $(T^*E^*,\pi_{E^*}, E^*)$
and $(TE, T\tau, TM)$. Formulae in local coordinates for these
algebroid structures are given in \cite{GraUrb2}, and other propertie of
these lifts can be found in \cite{Mack2}.

\begin{prop}
 \label{tangent lift of Lie algebroid}
Let $(E,\tau,M,\rho)$ be a Lie algebroid. Let $\Lambda_{E^*}$ be the
associated Poisson structure on the total space of the dual bundle $(E^*,\pi,M)$
(\ref{theorem poisson on dual bundle}) and $\Lambda_{TE^*}$ its lift to $TE$
(\ref{tangent lift of Poisson}).

\par\smallskip\noindent
{\rm 1.\quad} The Poisson structure $\Lambda_{TE^*}$ is homogeneous
(Definition~\ref{homogeneous Poisson structure})
for each of the two vector fibrations
$(TE^*,\tau_{E^*},E^*)$ and
$(TE^*,T\pi,TM)$.

\par\smallskip\noindent
{\rm 2.\quad} The vector bundle dual to $(TE^*,T\pi,TM)$. is
$(TE,T\tau,TM)$, and the Lie algebroid structure on that dual associated
to the homogeneous Poisson structure
$\Lambda_{TE^*}$ on the total space of $(TE^*,T\pi, TM)$ 
(\ref{theorem poisson on dual bundle})
is such that for each pair $(X,Y)$ of smooth sections of $\tau$, the
bracket $\{TX,TY\}$ is equal to $T\{X,Y\}$.

\par\smallskip\noindent
{\rm 3.\quad} The Lie algebroid structure on the vector bundle
$(T^*E^*,\pi_{E^*},E^*)$ associated to the homogeneous Poisson structure
$\Lambda_{TE^*}$ on the total space of its dual bundle
$(TE^*,\tau_{E^*},E^*)$ (\ref{theorem poisson on dual bundle}) is the same as
the Lie algebroid structure on the cotangent bundle
to the Poisson manifold $E^*$ (\ref{cotangent of Poisson}).

\end{prop}

\begin{proof}
Since $E^*$ is the total space of a vector bundle, $TE^*$ is a \emph{double
vector bundle} (\cite{Mack, Mack2, KonUrb}), {\it i.e.}, 
it is the total space of
two different vector fibrations: the tangent fibration $\tau_{E^*}:TE^*\to E^*$,
and the tangent lift $T\pi:TE^*\to TM$ of the vector fibration $\pi:E^*\to M$.
As a consequence of its definition, the Poisson structure
$\Lambda_{TE^*}$ is homogeneous with respect to the first vector fibration
$\tau_{E^*}:TE^*\to TM$. Let us
prove that it is homogeneous also with respect to the second. That Poisson
structure is characterized by the following
properties:  for each pair $(f,g)$ of smooth functions on $E^*$,
 $$\{df,dg\}_{TE^*}=d\{f,g\}_{E^*}\,,\quad
    \{f\circ\tau_{E^*},g\circ\tau_{E^*}\}_{TE^*}=0\,,\quad
    \{df,g\circ\tau_{E^*}\}_{TE^*}=\{f,g\}_{E^*}\circ\tau_{E^* }\,.
 $$
We need to prove first a part of Point~2: the duality between the vector
bundles $(TE^*,T\pi, TM)$ and $(TE,T\tau, TM)$. It is obtained by tangent lift
of the duality between $(E^*,\pi,M)$ and $(E,\tau,M)$. Let $Z\in TE$ and
$\Xi\in TE^*$ be such that $T\tau(Z)=T\pi(\Xi)$. There exist smooth curves
$t\mapsto \varphi(t)$  and $t\mapsto\psi(t)$, defined on an open interval $I$
containing $0$, with
values in $E$ and in $E^*$, respectively, such that 
$\displaystyle\frac{d\varphi(t)}{dt}\bigm|_{t=0}=Z$ and
$\displaystyle\frac{d\psi(t)}{dt}\bigm|_{t=0}=\Xi$. We may choose $\varphi$ and
$\psi$ such that $\tau\circ\varphi=\pi\circ\psi$, so
$\bigl\langle\psi(t),\varphi(t)\bigr\rangle$ is well defined for all $t\in I$.
We define
 $$\langle\Xi,Z\rangle=\frac{d\bigl\langle\psi(t),\varphi(t)\bigr\rangle}{dt}
 \bigm|_{t=0}\,.
 $$
The left hand side does not depend on the choices of $\varphi$ and
$\psi$, so it is a legitimate definition of$\langle\Xi,Z\rangle$. The vector
bundles $(TE^*,T\pi, TM)$ and
$(TE,T\tau, TM)$ are in duality.
\par\smallskip
For each smooth section $X:M\to E$ of $\tau$, $TX:TM\to TE$ is a smooth
section of $T\tau$. Let $\Phi_X:E^*\to\RR$ and $\Psi_{TX}:TE^*\to\RR$ be the
associated vertical functions (\ref{vertical functions}), defined, respectively
on the total spaces of the vector bundles $(E^*,\pi,M)$ and $(TE^*,T\pi,TM)$. 
For $\Xi\in TE^*$, let us calculate $\Psi_{TX}(\Xi)$. We take a smooth curve
$t\mapsto \psi(t)$ in $E^*$ such that
$\displaystyle\frac{d\psi(t)}{dt}\bigm|_{t=0}=\Xi$. The smooth curve in $E$,
$\varphi=X\circ\pi\circ\psi$, is such that $\tau\circ\varphi=\pi\circ \psi$
and
$\displaystyle\frac{d\varphi(t)}{dt}\bigm|_{t=0}=TX\bigl(T\pi(\Xi)\bigr)$, and
 $$\Psi_{TX}(\Xi)=\frac{d\bigl\langle\psi(t),\varphi(t)\bigr\rangle}{dt}
 \bigm|_{t=0}
 =\frac{d\Phi_X\bigl(\psi(t)\bigr)}{dt}\bigm|_{t=0}=d\Phi_X(\Xi)\,.
 $$
We have proven that $\Psi_{TX}=d\Phi_X$.
Now if $Y:M\to E$ is another smooth section of $\tau$, we have
 $$\{\Psi_{TX},\Psi_{TY}\}_{TE^*}=\{d\Phi_X,d\Phi_Y\}_{TE^*}
 =d\{\Phi_X,\Phi_Y\}_{E^*}=d\Phi_{\{X,Y\}}\,,$$
the last bracket $\{X,Y\}$ being the Lie algebroid bracket of sections $X$ and
$Y$ of $(E,\tau,M,\rho)$. These equalities prove that for each pair $(X,Y)$ of
smooth sections of $\tau$, the Poisson bracket $\{\Psi_{TX},\Psi_{TY}\}_{TE^*}$
is a vertical function. 
Proposition~\ref{properties of homogeneous Poisson structures} shows that
$\Lambda_{TE^*}$ is homogeneous with respect to the vector fibration
$(TE^*,T\pi,TM)$. Point~1 is
proven.
\par\smallskip
Theorem~\ref{theorem poisson on dual bundle} shows that associated
to the Poisson structure $\Lambda_{TE^*}$, we have Lie algebroid structures on
the dual bundles of $(TE^*,\tau_{E^*}, E^*)$ and $(TE^*,T\pi, TM)$. 
\par\smallskip
The last formula also proves that for each pair $(X,Y)$ of smooth sections of
$\tau$, the bracket $\{TX,TY\}$ is equal to $T\{X,Y\}$. So Point~2 is proven.
\par\smallskip
For the Lie algebroid structure on $(T^*E^*,\pi_{E^*}, E^*)$ considered as
cotangent bundle to the Poisson manifold $E^*,\Lambda_{E^*})$ 
(\ref{cotangent of Poisson}), the bracket of two sections of $\pi_{E^*}$, that
means the bracket of two Pfaff forms on $E^*$, is the bracket defined in
Remark~\ref{Poisson manifolds properties}~(i). Point~3 follows from the
properties of that bracket, as shown by Grabowski and
Urba\'nski \cite{GraUrb2}.
\end{proof}

%% ----------------------- REFERENCES -------------------------------%%

\end{document}